\documentclass[11pt]{amsart}
\usepackage{txfonts}
\usepackage{amsmath,amsthm,amscd,amssymb,mathtools}
\usepackage{enumitem}
\usepackage{tikz}
\usetikzlibrary{cd}
\usepackage[margin=1.07in]{geometry}
\usepackage[colorlinks=true,linkcolor=blue,citecolor=blue,urlcolor=blue,backref=page]{hyperref}
\renewcommand*{\backrefalt}[4]{%
  \ifcase #1 \or $\uparrow$\,#2\else $\uparrow$\,#2\fi}

\usepackage{orcidlink}

\theoremstyle{plain}
\newtheorem{theorem}{Theorem}[section]
\newtheorem{lemma}[theorem]{Lemma}
\newtheorem{proposition}[theorem]{Proposition}
\newtheorem{corollary}[theorem]{Corollary}

\theoremstyle{definition}
\newtheorem{definition}[theorem]{Definition}
\newtheorem{question}[theorem]{Question}

\newtheorem{remark}[theorem]{Remark}
\newtheorem{example}[theorem]{Example}

\theoremstyle{plain}
\newtheorem{theoremx}{Theorem}

\renewcommand{\k}{\Bbbk}
\newcommand{\N}{\mathbb{N}}
\newcommand{\Z}{\mathbb{Z}}
\newcommand{\Q}{\mathbb{Q}}
\newcommand{\R}{\mathbb{R}}
\newcommand{\C}{\mathbb{C}}
\newcommand{\F}{\mathbb{F}}
\newcommand{\A}{\mathcal{A}}
\newcommand{\CC}{\mathcal{C}}
\newcommand{\DD}{\mathcal{D}}
\newcommand{\E}{\mathcal{E}}
\newcommand{\FF}{\mathcal{F}}
\newcommand{\GG}{\mathcal{G}}
\newcommand{\RR}{\mathcal{R}}
\newcommand{\VV}{\mathcal{V}}
\newcommand{\ZZ}{\mathcal{Z}}
\newcommand{\one}{\mathbf{1}}
\newcommand{\g}{\mathfrak{g}}
\newcommand{\z}{\mathfrak{z}}

\DeclareMathOperator{\Spec}{Spec}
\DeclareMathOperator{\coker}{coker}

\DeclareMathOperator{\ab}{ab}
\DeclareMathOperator{\Hom}{Hom}
\DeclareMathOperator{\Ext}{Ext}

\DeclareMathOperator{\Mod}{Mod}
\DeclareMathOperator{\Sp}{Sp}
\DeclareMathOperator{\PSL}{PSL}
\DeclareMathOperator{\SL}{SL}
\DeclareMathOperator{\ad}{ad}

\DeclareMathOperator{\codim}{codim}
\newcommand{\id}{\mathrm{id}}

\newcommand{\LA}{\Lambda}
\newcommand{\zps}{\zeta_{p^s}}
\newcommand{\etaform}{$\eta$-formal}
\newcommand{\bwedge}{\mbox{$\bigwedge$}}

\newcommand{\surj}{\twoheadrightarrow}
\newcommand{\inj}{\hookrightarrow}
\newcommand\isom{\xrightarrow{ \,\smash{\raisebox{-0.5ex}{\ensuremath{\scriptstyle\simeq}}}\,}}

\newtheorem{propositionx}[theoremx]{Proposition}
\newtheorem{corollaryx}{Corollary}

\begin{document}

\title[Formality of finite abelian covers and Jordan blocks]%
{Formality of finite abelian covers\\ and Jordan blocks of the monodromy}

\author[Alexander~I.~Suciu]{Alexander~I.~Suciu\,\orcidlink{0000-0002-5060-7754}}
\address{Department of Mathematics, Northeastern University,
Boston, MA 02115, USA}
\email{\href{mailto:a.suciu@northeastern.edu}{a.suciu@northeastern.edu}}
\urladdr{\href{https://suciu.sites.northeastern.edu}%
{suciu.sites.northeastern.edu}}



\date{November 24, 2026}

\begin{abstract}
Let $X$ be a connected finite-type CW-complex with a map to the circle that is surjective on
fundamental groups. If $X$ is formal, the algebraic monodromy of the associated infinite cyclic cover
has no Jordan block of size greater than one at the eigenvalue $1$, over the rationals; formality of
the $d$-fold cyclic cover rules out such blocks at the $d$-th roots of unity. We study how formality of
finite abelian covers and Jordan blocks of the monodromy constrain each other. We call $X$
\emph{Alexander formal} if there are no such blocks at any root of unity, rationally, nor at the
eigenvalue $1$ modulo any prime; this is a formality property of cochains, as DG modules over the
cohomology of the circle. The complement of a central hyperplane arrangement, with its Milnor
fibration, is Alexander formal exactly when the numbers of hyperplanes in its irreducible factors have
no common divisor greater than $1$.

We also show that a CDGA graded by an abelian group, over an algebraically closed field of
characteristic zero, is formal if and only if it is formal as a graded CDGA. Hence the strong formality
of Kohno and Pajitnov implies formality of every finite abelian cover of $X$, but not conversely.
Modulo $p$, eigenvalues of $p$-power order collide with $1$, producing Jordan blocks that strong
formality allows; the two conditions are incomparable, and all four resulting classes contain
arrangement complements.
\end{abstract}

\maketitle
\setcounter{tocdepth}{1}
\tableofcontents

\section{Introduction}
\label{sec:intro}

Formality of a space, in the sense of Sullivan, has a known consequence for spaces that map to the
circle: the \emph{algebraic monodromy} of the associated infinite cyclic cover, that is, the
automorphism of its rational homology induced by the deck transformation, has no Jordan block of size
greater than one at the eigenvalue $1$. By the monodromy we always mean this operator, 
not the geometric monodromy that induces it. The $d$-fold cyclic cover of the space along the map has
the same infinite cyclic cover, with the $d$-th power of the monodromy, so a Jordan block of this kind
at a $d$-th root of unity obstructs the formality of the $d$-fold cover, even when the space itself is
formal. This paper studies how the formality of finite abelian covers and the Jordan blocks of the
monodromy constrain each other, rationally and modulo primes. Two conditions organize the answer.
\emph{Alexander formality}, introduced here, asks for the absence of such blocks at every root of
unity, rationally, and at the eigenvalue $1$ modulo every prime; it concerns a single map to the
circle, but integrally. The \emph{strong formality} of Kohno and Pajitnov concerns all rank-one local
systems at once, but only over $\C$; it implies the formality of every finite abelian cover
(Corollary~\ref{corx:cover}), by a result of independent interest about algebras graded by an abelian
group (Theorem~\ref{thmx:cover}). For the complement of a hyperplane arrangement, Alexander formality
is decided by the sizes of the irreducible factors (Theorem~\ref{thmx:arr}); in general, the two
conditions are incomparable (Theorem~\ref{thmx:incomp}).

\subsection{Formality and monodromy}
\label{subsec:setting}
A space $X$ is \emph{formal} if its algebra of rational polynomial forms is connected to its
cohomology algebra, with zero differential, by a chain of quasi-isomorphisms of CDGAs
\cite{Sullivan-ihes77}. Compact
K\"ahler manifolds \cite{DGMS-inv75} and, by a theorem of Brieskorn \cite{Brieskorn-71}, complements
of complex hyperplane arrangements are formal, while nonvanishing Massey products obstruct formality. 
Formality of $X$ implies \emph{$1$-formality} of $\pi_1(X)$, which means that the Malcev Lie algebra of 
the group is isomorphic, as a filtered Lie algebra, to the degree completion of its holonomy Lie algebra;
see \cite{Suciu-Wang-forum19} for a systematic treatment and \cite{Suciu-ems23} for a survey.

Let $X$ be a connected CW-complex of finite type and $\nu\colon\pi_1(X)\surj\Z$ an epimorphism;
equivalently, $\nu$ is induced by a map from $X$ to the circle. Let $X^\nu\to X$ be the infinite cyclic 
cover defined by $\ker\nu$, with deck transformation $t$. The automorphism $h_\ast$ induced by $t$ 
on homology is the \emph{algebraic monodromy}, and it makes
\[
A_\nu=H_\ast(X^\nu;\Z)\cong H_\ast\big(X;\Z[t^{\pm1}]_\nu\big)
\]
into a module over $\LA=\Z[t^{\pm1}]$, finitely generated in each degree: the \emph{integral Alexander
module} of $(X,\nu)$. Likewise, for a field $\k$, each $H_q(X^\nu;\k)$ is a finitely generated module
over the principal ideal domain $\k[t^{\pm1}]$. Throughout, the \emph{torsion} of a finitely generated
$\k[t^{\pm1}]$-module $M$ means its $\k[t^{\pm1}]$-torsion submodule, the set of elements annihilated
by a nonzero Laurent polynomial; it is finite-dimensional over $\k$ and contains every generalized
eigenspace of $t$, so the Jordan blocks of $t$ on $M$ are those of its restriction to the torsion.
For abelian groups, and for $\LA$-modules such as $A_\nu$, the terms ``$p$-torsion'' and
``torsion-free'' keep their usual meaning and refer to the underlying abelian group.

Formality and monodromy are linked by a dictionary between Jordan blocks and Massey products. For a
field $\k$, let $\eta\in H^1(X;\k)$ be the class of $\nu$, so that $\eta\cup\eta=0$. The dictionary has
its origins in the study of mapping tori. Fern\'andez, Gray, and Morgan used a Jordan block at the
eigenvalue~$1$ on the first cohomology of the fiber to show that certain mapping tori, and the circle
bundles over them, are not formal, so that such bundles, when symplectic, carry no K\"ahler metric
\cite[Prop.~9 and Note~10]{Fernandez-Gray-Morgan-mmj91}. Papadima and Suciu showed, using 
resonance rather than Massey products, that such a block obstructs even the $1$-formality of any 
bundle over the mapping torus with closed connected fiber \cite[Prop.~4.1]{Papadima-Suciu-forum10}; 
more generally, if $\pi_1(X)$ is $1$-formal and $H_1(X^\nu;\Q)$ is finite-dimensional, then $h_\ast$ 
has no Jordan block of size greater than one at the eigenvalue $1$ on $H_1(X^\nu;\Q)$
\cite[Prop.~3.2]{Papadima-Suciu-forum10}. Bazzoni, Fern\'andez, and Mu\~noz showed that some triple
product $\langle\eta,\eta,x\rangle$ does not vanish whenever, in some degree, the largest Jordan block at the
eigenvalue~$1$ has size two \cite[Thm.~13]{Bazzoni-Fernandez-Munoz-tams15}, and Kohno and Pajitnov
that, if $X$ is formal, then $h_\ast$ has no Jordan block of size greater than one at the eigenvalue $1$
on $H_\ast(X^\nu;\Q)$ \cite[Thm.~6.1]{Kohno-Pajitnov-cejm14}. In its general form, the dictionary is a 
theorem of Liu, Maxim, and Wang \cite[Thm.~1.1]{Liu-Maxim-Wang-ejm26}, valid over any field: the 
largest size of a Jordan block of $h_\ast$ for the eigenvalue $1$ on $H_q(X^\nu;\k)$ is one less
than the largest length of a nonvanishing iterated Massey product $\langle\eta,\dots,\eta,x\rangle$ 
with $x\in H^q(X;\k)$, so that a Jordan block of size two corresponds to a nonvanishing triple product
$\langle\eta,\eta,x\rangle$. In positive characteristic the products are computed in
the cochain algebra $C^\bullet(X;\k)$, from a fixed cocycle representing $\eta$
(Remark~\ref{rem:restricted}). Over $\C$, and for compact manifolds, Pajitnov had established the same
correspondence at every eigenvalue $\lambda$, for products with coefficients twisted by the character
$g\mapsto\lambda^{\nu(g)}$ \cite[Thm.~5.1]{Pajitnov-ejm17}, and later extended it to coefficients of
higher rank \cite[Thm.~C]{Pajitnov-topappl19}. The basic example is the Heisenberg nilmanifold, 
a torus bundle over the circle with a Jordan block of size two at the eigenvalue $1$ and a nonvanishing 
triple Massey product \cite[Exs.~2.4 and~9.5]{Papadima-Suciu-bullroum09}.

For $d\ge1$, let $X^{(d)}\to X$ be the $d$-fold cyclic cover defined by $\nu^{-1}(d\Z)$, taken with
the epimorphism $\nu^{(d)}=\nu/d$. Its infinite cyclic cover is again $X^\nu$, with deck transformation
$t^d$, so the theorem of Kohno and Pajitnov, applied to $X^{(d)}$, excludes Jordan blocks of size
greater than one at the $d$-th roots of unity when $X^{(d)}$ is formal (Remark~\ref{rem:intrinsic}).
For instance, let $X$ be the mapping torus of a diffeomorphism of a closed orientable surface $\Sigma$
of positive genus acting on $H_1(\Sigma;\Q)$ by Jordan blocks of size two at the eigenvalue $-1$, as in
Pajitnov's example \cite[\S4.2]{Pajitnov-ejm17}. Then $X$ is formal (Example~\ref{ex:paj}), but its
double cover $X^{(2)}$, the mapping torus of the square of the diffeomorphism, has unipotent monodromy
with Jordan blocks of size two, and so is not formal.

In positive characteristic the dictionary persists, but formality no longer controls it, even over the
field of coefficients. Three lines through the origin of $\C^2$, defined by $f=zw(z-w)$, already show
this. Their complement $M$ is formal, and over a field $\k$ the first homology of the Milnor fiber $F$,
a torus with three punctures, is $\k[t^{\pm1}]/(t^3-1)\oplus\k[t^{\pm1}]/(t-1)$
(Example~\ref{ex:toric-pencil}). This module is semisimple over $\C$, but over $\F_3$, where
$t^3-1=(t-1)^3$, the eigenvalues $\zeta_3$ and $\bar\zeta_3$ collide with $1$, and the monodromy has a
Jordan block of size three at the eigenvalue $1$. By the dictionary, $M$ carries a nonvanishing
quadruple product $\langle\eta,\eta,\eta,x\rangle$ over $\F_3$, with $\eta$ the class of the Milnor
fibration $f\colon M\to\C^\ast$. Yet $M\simeq S^1\times(S^1\vee S^1)$ and $F\simeq\bigvee^4S^1$ are
formal even over $\Z$, so the usual Massey products of $M$ over $\F_3$ all contain zero. The products of
the dictionary are more special, with defining systems built from a fixed cocycle pulled back from the
circle; they detect formality of $C^\bullet(M;\F_3)$ as a DG module over $H^\bullet(S^1;\F_3)$, not as a
DG algebra (Proposition~\ref{prop:circle} and Example~\ref{ex:three-lines}). Moreover, every finite
abelian cover of $M$ is formal (Theorem~\ref{thmx:incomp}\ref{x1} and Corollary~\ref{corx:cover}):
unlike the rational blocks above, this one obstructs the formality of none of them.

\subsection{Alexander formality}
\label{subsec:definition}
We base the definition on the monodromy itself, rationally at every root of unity and 
modulo every prime. The pair $(X,\nu)$, or simply $X$, is \emph{Alexander formal} if the following hold.
\begin{enumerate}[label=\textup{(\alph*)},leftmargin=2.2em,itemsep=2pt]
\item \label{AF1}
For every root of unity $\zeta$, the monodromy has no Jordan block of size greater than one for the
eigenvalue $\zeta$ on the torsion of $H_\ast(X^\nu;\C)$.
\item \label{AF2}
For every prime $p$, it has no Jordan block of size greater than one for the eigenvalue $1$ on the
torsion of $H_\ast(X^\nu;\F_p)$.
\end{enumerate}
The name refers to the Alexander module, but the definition is a formality condition in the literal
sense, imposed on the finite cyclic covers of $X$ along $\nu$ and on $X$ modulo every prime: by
Proposition~\ref{prop:circle} and Remark~\ref{rem:intrinsic}, $X$ is Alexander formal if and only if,
for every $d\ge1$, the rational cochains of $X^{(d)}$, and, for every prime $p$, the mod-$p$ cochains of
$X$, are formal as DG modules over the exterior algebra $H^\bullet(S^1)$, whose generator acts by cup
product with a cocycle pulled back from the circle.

When \ref{AF1} holds at $\zeta$ we say that $X$ is \emph{\etaform{} over $\Q(\zeta)$ at $\zeta$}, or
simply \emph{\etaform{} over $\Q$} when $\zeta=1$, and when \ref{AF2} holds at $p$ we say that $X$ is
\emph{\etaform{} over $\F_p$}. By the dictionary, these are vanishing conditions on the Massey
products along $\eta$, twisted by $g\mapsto\zeta^{\nu(g)}$ in the first case and computed modulo $p$
in the second; each is the degeneration at $E_2$ of the equivariant spectral sequence of $(X,\eta)$
\cite{Papadima-Suciu-tams10}, over the appropriate field and after the appropriate twist
(Definitions~\ref{def:etak} and~\ref{def:Aformal}).

Formality of $X$, or of $X^{(d)}$, implies \ref{AF1} at $\zeta=1$, or at the $d$-th roots of unity,
but formality of $X$ implies neither the other clauses of \ref{AF1} nor \ref{AF2}, by the examples of
\S\ref{subsec:setting}. 
The asymmetry between \ref{AF1} and \ref{AF2} is
deliberate: modulo $p$, every root of unity of $p$-power order reduces to $1$, since
$\Phi_{p^s}\equiv(t-1)^{\varphi(p^s)}\bmod p$, so \ref{AF2} at $p$ constrains all these eigenvalues at
once, together with the $p$-torsion of $A_\nu$ (Lemma~\ref{lem:collision} and
Remarks~\ref{rem:scope} and~\ref{rem:allzeta}). When $H_\ast(X^\nu;\Z)$ is finitely generated, each
clause reduces to the acyclicity of an Aomoto complex, $\big(H^\bullet(X;\F_p),\eta\cup-\big)$ or its
analogue with twisted complex coefficients (Corollary~\ref{cor:finite-type}).

\subsection{Arrangement complements}
\label{subsec:diagonal}
Hyperplane arrangements are where Alexander formality can be decided completely. Let $\A$ be a
central arrangement of $n$ hyperplanes in $\C^\ell$, with complement $M(\A)$, and let $\nu$ be the
class of the Milnor fibration $M(\A)\to\C^\ast$ given by the defining polynomial. The infinite cyclic
cover is then the Milnor fiber $F(\A)$, a finite complex, and the monodromy has order $n$, so
condition~\ref{AF1} holds automatically. Condition~\ref{AF2} does not.

\begin{theoremx}[see Theorem~\ref{thm:arr}]
\label{thmx:arr}
Let $\A=\A_1\times\dots\times\A_k$ be the decomposition of a central arrangement into irreducible
arrangements, and let $n_i=\lvert\A_i\rvert$. Then $M(\A)$, with the class of its Milnor fibration,
is Alexander formal if and only if, for every prime $p$ dividing $n=n_1+\dots+n_k$, some $n_i$ is prime
to $p$; equivalently, if and only if $\gcd(n_1,\dots,n_k)=1$. In particular, if $\A$ is irreducible and
$n\ge2$, then $M(\A)$ is not $\eta$-formal over $\F_p$ for any prime $p\mid n$, and so is not
Alexander formal.
\end{theoremx}

\noindent By Proposition~\ref{prop:circle}, for $p\mid n$ the mod-$p$ cochains of the complement of an
irreducible arrangement are not formal as a DG module over the cohomology of the circle, although its
rational cochains are formal as a DG algebra. The reason is an Euler characteristic. The Aomoto complex
$\big(H^\bullet(M(\A);\k),\nu\cup-\big)$ is acyclic over every field in which $n$ is invertible
(Lemma~\ref{lem:homotopy}). If $\A$ is irreducible and $p\mid n$, the mod-$p$ Aomoto complex is the
tensor product of an exterior algebra on one generator, with zero differential, and a complex whose
Euler characteristic is that of the projectivized complement; up to sign, this is Crapo's beta invariant
of the underlying matroid, which is nonzero for irreducible $\A$ \cite{Crapo-jct67}. In terms of the
monodromy, every $n$-th root of unity is then an eigenvalue in some degree, and, when $p^s>2$ divides
$n$, an eigenvalue of order $p^s$ collides modulo $p$ with the eigenvalue $1$, as for the three lines
(Remark~\ref{rem:arr-collision}). Alexander formality of an arrangement complement is thus determined
by the intersection lattice, and cannot distinguish the members of a Zariski pair. Nor need the failure
be visible in degree one: for many arrangements and primes $p\mid n$, the mod-$p$ Aomoto complex is
exact in degree one (Example~\ref{ex:betti}).

Condition~\ref{AF1}, automatic for complex arrangements, is a genuine constraint for arrangements of
real planes: for many horizontal arrangements of $n$ transverse planes in $\R^4$
\cite{Matei-Suciu-top00}, the monodromy on $H_1(X^\nu;\Q)$, for $X$ the complement and $\nu$ sending
each meridian to $1$, has Jordan blocks of size greater than one at roots of unity $\zeta\ne1$, already
for $n=4$ at $\zeta=-1$ \cite[Prop.~10.10]{Suciu-planes-r4}.

\subsection{Strong formality and finite abelian covers}
\label{subsec:graded}
Strong formality, in the sense of Kohno and Pajitnov, controls the finite abelian covers through a
result about graded algebras. This result mentions neither $X$ nor $A_\nu$, and we prove it in that
generality. A CDGA $\bar A$ over a field is \emph{graded by an abelian group} $\CC$ if it is given a
decomposition $\bar A=\bigoplus_{\rho\in\CC}\bar A_\rho$ into subcomplexes such that $\bar
A_\rho\cdot\bar A_\tau\subseteq\bar A_{\rho\tau}$. Such an algebra is \emph{formal} if it can be joined
to its cohomology by a chain of quasi-isomorphisms of CDGAs, and \emph{$\CC$-graded formal} if that
chain can be chosen among $\CC$-graded CDGAs and grading-preserving quasi-isomorphisms. 
The second condition implies the first.

\begin{theoremx}[see Theorem~\ref{thm:G}]
\label{thmx:cover}
Let $\CC$ be an abelian group and let $\bar A$ be a CDGA graded by $\CC$ over a field $\k$
of characteristic zero such that the characters $\CC\to\k^\ast$ separate the points of $\CC$
\textup{(}for instance, $\k$ algebraically closed\textup{)}. If $\bar A$ is formal, then it is formal
as a $\CC$-graded CDGA. In particular the sub-CDGA $\bigoplus_{\rho\in\DD}\bar A_\rho$ is
formal for \emph{every} subgroup $\DD\subseteq\CC$, of any cardinality and any index.
\end{theoremx}

The key point is that a grading by an abelian group amounts to an action of its character
group by diagonalizable automorphisms. The proof passes to the minimal $C_\infty$-model of $\bar A$,
where formality becomes the gauge equivalence of two invariant Maurer--Cartan elements, and makes the
gauge invariant as well by a descent argument whose only arithmetic ingredient is Artin's independence
of characters (Lemmas~\ref{lem:charcocycle} and~\ref{lem:descent}); it needs characteristic zero
(Remark~\ref{rem:Gproof}).

The algebra to which we apply Theorem~\ref{thmx:cover} was introduced by Kohno and Pajitnov
\cite[\S6]{Kohno-Pajitnov-cejm14}. Let $\CC=H^1(X;\C^\ast)$ be the group of rank-one characters of
$\pi_1(X)$, and for $\rho\in\CC$ let $A^\bullet(X,E_\rho)$ be the complex of differential forms with
values in the flat line bundle $E_\rho$ with monodromy $\rho$, whose cohomology is
$H^\bullet(X;\C_\rho)$; for a finite complex one works on a manifold of the same homotopy type
(\S\ref{subsec:KP}). The isomorphisms $E_\rho\otimes E_\tau\cong E_{\rho\tau}$ make
\[
\bar A^\bullet(X)=\bigoplus_{\rho\in\CC}A^\bullet(X,E_\rho)
\]
a CDGA graded by $\CC$, and $X$ is \emph{strongly formal} if $\bar A^\bullet(X)$ is formal
\cite[Def.~6.2]{Kohno-Pajitnov-cejm14}. This is not the bigraded notion of Milivojevi\'c and Stelzig
\cite{Milivojevic-Stelzig-cag24}, which bears the same name; the two are separated by a curve 
of genus two (Remark~\ref{rem:clash}). Compact K\"ahler manifolds are strongly formal
\cite[Thm.~6.5]{Kohno-Pajitnov-cejm14}, and strong formality excludes Jordan blocks of size greater
than one at every eigenvalue, for every $\nu$ (Proposition~\ref{prop:sf-char0}). The link with
Theorem~\ref{thmx:cover} is that a connected finite abelian cover $Y\to X$ with deck group $G$ has
$A^\bullet(Y)\cong\bigoplus_{\rho\in\widehat G}A^\bullet(X,E_\rho)$: the CDGA of the cover is the part
of $\bar A^\bullet(X)$ lying over the subgroup $\widehat G=\Hom(G,\C^\ast)$ of $\CC$
(Lemma~\ref{lem:graded}).

\begin{corollaryx}[see Proposition~\ref{prop:cover}, Corollary~\ref{cor:allD}]
\label{corx:cover}
If $X$ is strongly formal, then every connected finite abelian cover of $X$ is formal.
More generally, for every subgroup $\DD\subseteq\CC$ the CDGA
$\bigoplus_{\rho\in\DD}A^\bullet(X,E_\rho)$ is formal, so all of its Massey products,
of all lengths, vanish.
\end{corollaryx}

Already the trivial subgroup is informative: taking $\DD=\{\one\}$ shows that strong formality implies
formality of $X$, which does not follow from the definition alone, since the projection of
$\bar A^\bullet(X)$ onto its summand $A^\bullet(X)$ is not multiplicative. The abelian hypothesis cannot
be dropped: the Brieskorn sphere $\Sigma(2,3,7)$ is strongly formal, and its $\PSL(2,7)$-cover, a circle
bundle over the Klein quartic, is not $1$-formal (Example~\ref{ex:nonab-cover}).

For every epimorphism $\nu\colon\pi_1(X)\surj\Z$, the results in characteristic zero thus form a chain
of implications,
\[
X\ \text{strongly formal}\ \Longrightarrow\ \text{every finite abelian cover of $X$ formal}\
\Longrightarrow\ (X,\nu)\ \text{satisfies \ref{AF1}},
\]
the second arrow passing through the cyclic covers $X^{(d)}$. Neither arrow can be reversed in general.
The first fails for a surface bundle over the circle with monodromy $J_2(\lambda)\oplus
J_2(\lambda^{-1})$ on $H_1$, where $\lambda$ is not a root of unity, whose finite abelian covers are all
formal (Proposition~\ref{prop:converse}): its obstruction sits at a character of infinite order, which
no finite cover involves; for $b_1(X)=0$ the first arrow is an equivalence. The second fails for the
complement of the Ceva arrangement, with its Milnor fibration: the complement is formal and its
monodromy has finite order, so \ref{AF1} holds, but the $9$-fold cyclic cover along $\nu$, the product of
the Milnor fiber with $\C^\ast$, is not formal \cite{Zuber-blms10}, although its own monodromy is
trivial (Remark~\ref{rem:line}). Condition~\ref{AF2} lies off the chain: modulo a prime $p$, eigenvalues
of $p$-power order collide with $1$, and the resulting Jordan blocks are compatible with strong
formality (Theorem~\ref{thmx:incomp}\ref{x1}).

For a finite grading group, Theorem~\ref{thmx:cover} reduces to the classical descent of formality
along finite covers. The case needed here is the opposite one: for the projectivized complement of an
arrangement of $n$ hyperplanes, $\CC=(\C^\ast)^{n-1}$ is divisible, and the finite subgroup $\mu_n$ that
classifies the Milnor fiber is neither a direct summand of $\CC$ nor detected by any homomorphism from
$\CC$ to a finite group. With Corollary~\ref{corx:cover}, if the Milnor fiber of an arrangement, or any
other connected finite abelian cover of the projectivized complement, is not formal, as for the Ceva
arrangement, then the complement is formal but not strongly formal (Corollary~\ref{cor:milnor-strong}).
Strong formality is also the rank-one end of Pajitnov's scale of $\FF$-formality
\cite[Def.~3.9]{Pajitnov-topappl19}, which Theorem~\ref{thmx:cover} reaches, and its proof a little
beyond (\S\ref{subsec:Fformal}).

Theorem~\ref{thmx:cover} also controls the cohomology jump loci of a strongly formal space, at every
character and not only at the trivial one. Let $\VV^i_k(X)\subseteq\CC$ be the set of characters
$\rho$ with $\dim_\C H^i(X;\C_\rho)\ge k$, and let $\RR^i_k(X,\rho)\subseteq H^1(X;\C)$ be the
set of classes $\omega$ for which the complex $\big(H^\bullet(X;\C_\rho),\omega\cup-\big)$ has
$i$-th cohomology of dimension at least $k$. Libgober showed that the tangent cone to
$\VV^i_k(X)$ at $\rho$ always lies in $\RR^i_k(X,\rho)$ \cite{Libgober-topappl02}; for
formal $X$ and $\rho=\one$, the two agree as germs 
\cite{Dimca-Papadima-Suciu-duke09,Dimca-Papadima-ccm14}.

\begin{theoremx}[see Theorem~\ref{thm:jump}]
\label{thmx:jump}
Let $X$ be a strongly formal space with the homotopy type of a finite CW-complex, and let
$\rho\in\CC$. For all $i$ and $k$, the map $\omega\mapsto\rho\exp(\omega)$ identifies the germ of
$\RR^i_k(X,\rho)$ at $0$ with the germ of $\VV^i_k(X)$ at $\rho$. Consequently, every
irreducible component of $\RR^i_k(X,\rho)$ is a linear subspace defined over $\Q$, and every
irreducible component of $\VV^i_k(X)$ is a translated subtorus of $\CC$.
\end{theoremx}

\noindent The proof restricts the $\CC$-graded formality of Theorem~\ref{thmx:cover} to the degrees
$\one$ and $\rho$, and applies results of Budur and Wang \cite{Budur-Wang-comp15a,Budur-Wang-pamq20}.
The last assertion gives an obstruction to strong formality that involves neither finite covers nor
Massey products, only the shape of a characteristic variety away from the trivial character, where
formality says nothing: the complement of the link $6^2_3$ is formal, but one component of its first
characteristic variety is a curve that is not a translated subtorus, so it is not strongly formal
(Example~\ref{ex:link623}). For smooth quasi-projective varieties, whose jump loci are unions of
torsion-translated subtori \cite{Budur-Wang-ens15}, only the germ statement carries information
(Remark~\ref{rem:jump}). Theorem~\ref{thmx:jump} cannot detect every failure: the complement of the
braid arrangement is not strongly formal, yet the conclusions of the theorem hold for it at every
character (see the discussion after Proposition~\ref{prop:A3-tc}).

\subsection{Alexander formality versus strong formality}
\label{subsec:results}
\begin{theoremx}[see Theorem~\ref{thm:incomp} and Corollary~\ref{cor:arr-corners}]
\label{thmx:incomp}
Alexander formality and Kohno--Pajitnov strong formality are incomparable.
\begin{enumerate}[label=\textup{(\roman*)},leftmargin=2.2em,itemsep=2pt]
\item \label{x1} Strong formality excludes Jordan blocks of size greater than one at every
eigenvalue over $\C$, so a strongly formal space can fail to be Alexander formal only at a prime, and
it does: bielliptic surfaces of types $\Z/3$ and $\Z/4$, with a suitable class $\nu$, and
complements of pencils of at least three lines, are strongly formal and not Alexander formal. The
obstruction is the congruence $\Phi_{p^s}\equiv(t-1)^{\varphi(p^s)}\bmod p$, which turns a
semisimple rational eigenvalue of order $p^s>2$ into a Jordan block of size greater than one at $t=1$
modulo $p$ (Corollary~\ref{cor:collision}).
\item \label{x2} Alexander formality does not imply strong formality: some fibered $3$-manifold with
$b_1=1$ whose monodromy has a Jordan block of size two at an eigenvalue of infinite order
is Alexander formal and
not strongly formal, and so is the complement of the product arrangement $\A_3\times\A_1$.
\item \label{x3} All four classes cut out by the two conditions are nonempty, and each contains an
arrangement complement.
\end{enumerate}

\end{theoremx}

For part~\ref{x2}, note that Alexander formality constrains only the Massey products along $\eta$: a
class $\nu$ that splits off a circle factor, $X=S^1\times Y$, makes $X$ Alexander formal whatever $Y$ is
(Proposition~\ref{prop:dir}), and the product arrangement is of this kind. The
mapping torus shows that the gap persists when $\nu$ is unique up to sign. Each of the four 
classes contains the following examples.

\vspace*{2pt}
\begin{center}
\small
\renewcommand{\arraystretch}{1.3}
\begin{tabular}{l|l|l}
\hline
& \textit{strongly formal} & \textit{not strongly formal}\\
\hline
\textit{Alexander formal}
& complex torus; bielliptic of types $\Z/2$, $\Z/6$;
& mapping torus (Prop.~\ref{prop:converse});\\
& Boolean arrangement
& product $\A_3\times\A_1$\\
& (\S\ref{subsec:corners}, Cors.~\ref{cor:bielliptic}, \ref{cor:arr-corners})
& (Prop.~\ref{prop:dir}, Cor.~\ref{cor:arr-corners})\\
\hline
\textit{not Alexander formal}
& bielliptic of types $\Z/3$, $\Z/4$;
& mapping tori (Ex.~\ref{ex:paj}, Ex.~\ref{ex:phi6});\\
& generic arrangements; $2$-trees
& Ceva $\A(3,3,3)$; braid $\A_m$ with $m\ge3$\\
& (Cor.~\ref{cor:bielliptic}, Cor.~\ref{cor:generic}, Ex.~\ref{ex:2tree})
& (Cor.~\ref{cor:arr-corners}, Prop.~\ref{prop:nark}, Cor.~\ref{cor:A3-loc})\\
\hline
\end{tabular}
\end{center}

\vspace*{2pt}
The two arrangements in the lower right cell fail to be strongly formal for different reasons: the
Ceva arrangement through the non-formal finite cyclic cover described above, and the braid arrangement
through a nonvanishing triple Massey product over a subgroup $\DD\cong\Z^2$ of characters, found by
Narkawicz \cite{Narkawicz-thesis07} (Example~\ref{ex:nark}). Since $\DD$ is infinite and of infinite
index, no finite-cover argument applies, and the conclusion does not depend on the open question of
whether the Milnor fiber of $\A_3$ is $1$-formal. Both arrangements fail Alexander formality for the
arithmetic reason of Theorem~\ref{thmx:arr}; the product $\A_3\times\A_1$ removes that reason and keeps
Narkawicz's obstruction. Whether the complement of Yoshinaga's icosidodecahedral arrangement, which is
not Alexander formal, is strongly formal we do not know (Remark~\ref{rem:corners}). Toric complexes,
which are all formal, occupy the four classes as well (\S\ref{subsec:toric}), and the classifying space
of the right-angled Artin group of a graph $\Gamma$ is strongly formal if and only if $\Gamma$ has no
induced path on four vertices (Corollary~\ref{cor:toric-P4}). For arrangements, localizations play the
role of full subcomplexes: the failure for $\A_3$ propagates to every arrangement with a localization
isomorphic to $\A_3$, such as the braid arrangements $\A_m$ with $m\ge3$ and the graphic arrangements of
graphs containing the complete graph $K_4$ (Corollary~\ref{cor:A3-loc}), and the graphic arrangement of
a chordal graph $\Gamma$ has strongly formal complement if and only if $\Gamma$ does not contain $K_4$
(Theorem~\ref{thm:chordal}).

\subsection{The classical constraints}
\label{subsec:question}
Alexander formality was designed to assemble several known constraints on the algebraic monodromy,
recalled in \S\ref{subsec:constraints}: rationally and at
the eigenvalue $1$, the theorem on $1$-formality quoted above \cite{Papadima-Suciu-forum10}; modularly,
the bounds on the twisted Betti numbers at characters of $p$-power order by the mod-$p$ Aomoto--Betti
numbers \cite[Thm.~C]{Papadima-Suciu-tams10}, which can be strict because of torsion in the homology of
the cover \cite{Yoshinaga-ejm20}; and, over $\C$ and at all characters at once, strong formality. Each is
an instance of the dictionary of \S\ref{subsec:setting}, over a different field or after a different
twist, and together they raise the question that motivated the present work.

\begin{question}
\label{q:main}
Does the pair of conditions
\begin{center}
\emph{$1$-formality at the eigenvalue $\lambda=1$, together with $\eta$-formality at every
prime-power character}
\end{center}
assemble into a single formality property of the $\Z[t^{\pm1}]$-module $A_\nu$? And if it does, is that
property equivalent to the all-character strong formality of \cite{Kohno-Pajitnov-cejm14}, or do the
higher-rank, cross-character Massey products of \cite{Kohno-Pajitnov-cejm14,Narkawicz-thesis07} carry
strictly more information?
\end{question}

The answer to the first clause is affirmative, after two adjustments, and the resulting notion is
Alexander formality.

\begin{propositionx}[see Proposition~\ref{thm:assembly} and Corollary~\ref{cor:assembly}]
\label{thmx:assembly}
The clause of Alexander formality at $\zeta=1$, namely $\eta$-formality over $\Q$, is implied by
formality of $X$, and in degree one by $1$-formality of $\pi_1(X)$ when $H_1(X^\nu;\Q)$ is
finite-dimensional. When the $p$-power cyclic covers of $X$ along $\nu$ have no $p$-torsion in their
homology, the clause at a prime $p$, namely $\eta$-formality over $\F_p$, is equivalent to sharpness of
the modular Aomoto--Betti bounds at every character of $p$-power order and in every degree. Consequently,
reading the second condition of Question~\ref{q:main} both in characteristic zero and modulo $p$,
Alexander formality is the conjunction of its conditions after two changes:
$1$-formality of $\pi_1(X)$ is replaced by $\eta$-formality over $\Q$ in every degree, and
$\eta$-formality in characteristic zero is required at every root of unity, not only at those of
prime-power order.
\end{propositionx}

\noindent Proposition~\ref{thmx:assembly} reorganizes results already in the literature, and neither
change can be omitted (Corollary~\ref{cor:assembly}): $S^1\times H$, with $H$ the Heisenberg
nilmanifold, is Alexander formal, although its fundamental group is not $1$-formal; and the fibered
$3$-manifold of Example~\ref{ex:phi6} satisfies both conditions of the question, but has a rational
Jordan block of size two at a primitive sixth root of unity.

The answer to the second clause is negative, by Theorem~\ref{thmx:incomp}, and the difference does
\emph{not} lie in the rank of the coefficients: by a theorem of Pajitnov
\cite[Thm.~C]{Pajitnov-topappl19}, the products $\langle\eta,\dots,\eta,x\rangle$ with coefficients of
any rank are again governed by the Jordan blocks of twisted Alexander modules
(Remark~\ref{rem:higher-rank}). What lies beyond the reach of $A_\nu$ is the rest of the multiplicative
structure of $\bar A^\bullet(X)$: on the line of characters that factor through $\nu$, as for the Ceva
arrangement (Remark~\ref{rem:line}), or in two independent degrees of the $\CC$-grading, as for
Narkawicz's product, which is defined on no sub-CDGA attached to a single epimorphism
$\nu\colon\pi_1(X)\surj\Z$ (Proposition~\ref{prop:nark}).

\subsection{Organization and notation}
\label{subsec:organization}
Sections~\ref{sec:setup}--\ref{sec:first} set up the equivariant spectral sequence and the
dictionary, give the definitions, and assemble the classical constraints; \S\ref{sec:covers} proves
Theorem~\ref{thmx:cover}, \S\ref{sec:strong} derives Corollary~\ref{corx:cover} and
Theorem~\ref{thmx:jump}, \S\ref{sec:second} proves Theorem~\ref{thmx:incomp}, and
\S\ref{sec:functorial} treats arrangement complements. Each section opens with a summary of its
contents. 

Throughout, $\LA=\Z[t^{\pm1}]$, and $A_\nu=H_\ast(X^\nu;\Z)$ is the integral Alexander module of
$(X,\nu)$, the subscript recording its dependence on $\nu$. With a bullet, $A^\bullet$ is
the Orlik--Solomon algebra $A^\bullet(\A)$, a twisted de Rham complex $A^\bullet(X,E_\rho)$, or
the Kohno--Pajitnov algebra $\bar A^\bullet(X)$. For a central arrangement $\A$ of $n$ hyperplanes in
$\C^\ell$ we write $M(\A)$ for its complement, $U(\A)=\mathrm{P}(M(\A))$ for the projectivized
complement, $F(\A)$ for the Milnor fiber of the defining polynomial, a regular $\mu_n$-cover 
of $U(\A)$, and $\nu=\omega_1+\dots+\omega_n\in H^1(M(\A);\Z)$ for the \emph{diagonal class}, 
the class of the Milnor fibration, where $\omega_1,\dots,\omega_n$ are the classes of the 
logarithmic forms; then $M(\A)\cong U(\A)\times\C^\ast$. The Aomoto--Betti numbers $\beta^p_q(X,\eta)$ and
the twisted Betti numbers $b_q(X,\zeta)$ are defined in \S\ref{subsec:ss}.

\section{Setup and the spectral sequence}
\label{sec:setup}

Throughout, $X$ is a connected CW-complex of finite type, $\nu\colon\pi_1(X)\surj\Z$ is an
epimorphism, $A_\nu=H_\ast(X;\LA_\nu)\cong H_\ast(X^\nu;\Z)$ is the integral Alexander module, and
$h_\ast$ is the algebraic monodromy. For a field $\k$ we put $\eta=\nu_\k\in H^1(X;\k)$; since $\eta$ 
is pulled back from $H^1(S^1;\k)$, we have $\eta\cup\eta=0$. The results of \cite{Liu-Maxim-Wang-ejm26} quoted 
below are stated for finite complexes, and those of \cite{Pajitnov-ejm17,Pajitnov-topappl19} for compact
manifolds; each concerns a single degree, depends only on a finite skeleton, and so applies to
complexes of finite type, a finite complex being homotopy equivalent to a compact manifold with
boundary, to which Pajitnov's arguments apply unchanged: they use only the finite generation of the
twisted Alexander modules and the spectral sequences of \cite{Kohno-Pajitnov-cejm14}, which are set up
for arbitrary smooth manifolds.

This section recalls the equivariant spectral sequence of $(X,\eta)$ and the dictionary that reads the
Jordan blocks of $h_\ast$ off its differentials, which are Massey products, and shows that degeneration of
the sequence is the formality of the cochains of $X$ as a DG module over the cohomology of the circle
(\S\ref{subsec:ss}). It then explains which part of $A_\nu$ is constrained by degeneration over a given 
field and after a given twist (\S\ref{subsec:translations}), and proves the collision lemma, 
Lemma~\ref{lem:collision}, the arithmetic input of the paper.

\subsection{The equivariant spectral sequence}
\label{subsec:ss}
Let $R=\k[t^{\pm1}]$, let $s=t-1$, and let $\hat R=\k\llbracket s\rrbracket$ be the $(t-1)$-adic completion 
of $R$. Let $L_\nu$ be the local system on $X$ with fiber $R$ on which $\pi_1(X)$ acts through $\nu$, a
generator of $\Z$ acting by multiplication by $t$, so that $H_\ast(X;L_\nu)\cong H_\ast(X^\nu;\k)$.
As in \cite[\S2]{Liu-Maxim-Wang-ejm26}, we take $X$ to be a simplicial complex, or more generally a
$\Delta$-complex, and the map to the circle to be simplicial, and we let the cocycle representing $\eta$ 
be the pullback of the $1$-cochain dual to one
edge of the circle; since the circle has no $2$-simplices, $\eta\cup\eta=0$ already on cochains. The
complex $C^\bullet(X;\k)\otimes_\k\hat R$, with differential $\delta+(\eta\cup-)\otimes s$, then computes
$H^\bullet(X;L_\nu\otimes_R\hat R)$ \cite[Cor.~2.4]{Liu-Maxim-Wang-ejm26}. Filtering it by
the powers of $s$ gives a spectral sequence with
\[
E_1^{i,j}=H^{i+j}(X;\k)\cdot s^i\quad(i\ge0),\qquad d_1=\eta\cup- ,
\]
so that $E_2^{i,j}=H^{i+j}\big(H^\bullet(X;\k),\eta\cup-\big)$ for every $i\ge1$
\cite[\S3]{Liu-Maxim-Wang-ejm26}. The cochain complex $\big(H^\bullet(X;\k),\eta\cup-\big)$ is the
\emph{Aomoto complex} of $\eta$, and we write
\[
\beta^{\k}_q(X,\eta)=\dim_\k H^q\big(H^\bullet(X;\k),\,\eta\cup-\big)
\]
for its \emph{Aomoto--Betti numbers}. Since $H^\bullet(X;\k)=H^\bullet(X;\k_0)\otimes_{\k_0}\k$ for the
prime field $\k_0$ of $\k$, and $\eta$ is defined over $\k_0$, these numbers depend only on the
characteristic $p\ge0$ of $\k$, and we write $\beta^p_q(X,\eta)$ for them; thus $E_2^{i,j}$ has dimension
$\beta^p_{i+j}(X,\eta)$ for $i\ge1$. The sequence is dual to the
$J$-adic spectral sequence of Papadima and Suciu \cite{Papadima-Suciu-tams10}, see
\cite[Rem.~3.2]{Liu-Maxim-Wang-ejm26}, and we call either one the \emph{equivariant spectral
sequence} of $(X,\eta)$. Its higher differentials are Massey products of a special kind, built from
the cocycle $\eta$ \cite[Defs.~3.3 and~3.4]{Liu-Maxim-Wang-ejm26}. As in
\cite[\S3]{Kohno-Pajitnov-cejm14}, an \emph{$r$-chain} starting from $\omega\in H^q(X;\k)$ is a sequence
$\omega_1,\dots,\omega_r$ in $C^q(X;\k)$ with $\delta\omega_1=0$, $[\omega_1]=\omega$, and
$\delta\omega_{j+1}=\eta\cup\omega_j$ for $j<r$; an \emph{infinite chain} is defined in the same way. If
$\omega$ starts an $r$-chain, the $(r+1)$-fold product $\langle\eta,\dots,\eta,\omega\rangle$, with $r$
entries $\eta$, is the set of classes $[\eta\cup\omega_r]\in H^{q+1}(X;\k)$ over all such chains. It is a
coset of a subspace, its indeterminacy, and it \emph{vanishes} if it contains $0$, that is, if $\omega$
starts an $(r+1)$-chain; its \emph{length} is $r+1$, so that $\langle\eta,\omega\rangle=\eta\cup\omega$
has length two. For $r\ge2$ the differential $d_r$ is computed by these products
\cite[\S3]{Liu-Maxim-Wang-ejm26}, and by \cite[Thm.~1.1]{Liu-Maxim-Wang-ejm26}, in each degree $q$ the
maximal size of a Jordan block of the eigenvalue-$1$ part of the torsion of $H_q(X^\nu;\k)$ is one less
than the maximal length of a nonvanishing product $\langle\eta,\dots,\eta,\omega\rangle$ with
$\omega\in H^q(X;\k)$. Consequently the following three conditions are equivalent:
\begin{enumerate}[label=\textup{(\alph*)},leftmargin=2.2em,itemsep=2pt]
\item \label{dg1} The equivariant spectral sequence of $(X,\eta)$ degenerates at $E_2$.
\item \label{dg2} All products $\langle\eta,\dots,\eta,\omega\rangle$ of length at least three vanish.
\item \label{dg3} The $(t-1)$-primary part of the torsion of $H_\ast(X^\nu;\k)$ is semisimple.
\end{enumerate}
When they hold, $\beta^p_q(X,\eta)$ is the rank of $H_q(X^\nu;\k)$ as a module over $\k[t^{\pm1}]$
(for $\k=\C$ this is observed in \cite[\S6]{Kohno-Pajitnov-cejm14}); when $X^\nu$ has
finite-dimensional homology, they hold if and only if the Aomoto complex of $\eta$ is acyclic
(Proposition~\ref{prop:finite-type}). Degree by degree, a product $\langle\eta,\dots,\eta,\omega\rangle$
of length $r+1\ge3$ with $\omega\in H^q(X;\k)$ computes $d_r$ on the class of $\omega$, with values in
a subquotient of $E_2^{r,q+1-r}$, a space of dimension $\beta^p_{q+1}(X,\eta)$; so if
$\beta^p_{q+1}(X,\eta)=0$, all these products vanish, and $H_q(X^\nu;\k)$ has no Jordan block of size
greater than one at the eigenvalue~$1$. In positive characteristic, condition~\ref{dg2} concerns these
products, not the classical Massey products of $C^\bullet(X;\k)$, which may all vanish while
\ref{dg2} fails (Remark~\ref{rem:restricted}).

In characteristic zero, this sequence is the restriction to the line through $\eta$ of the Koszul
spectral sequence of \cite[\S3.2]{Suciu-arx26b}, which treats all classes in $H^1$ at once and is
carried out on algebraic models; its higher differentials are iterated Massey products along the
canonical class \cite[Thm.~5.1]{Suciu-arx26b}, and it collapses at $E_2$ for formal CDGAs
\cite[Cor.~5.9]{Suciu-arx26b}. Along a single class the boundary map has only the two components
$\delta$ and $s(\eta\cup-)$, unlike its multivariable counterpart \cite[Thm.~16.3]{Suciu-arx26b};
this is why the theorem of Liu, Maxim, and Wang holds over every field, as the rank-one instance of
\cite[Thm.~5.1(3)]{Suciu-arx26b} (see \cite[Rem.~16.17]{Suciu-arx26b}).

Over $\C$ the dictionary holds at every eigenvalue, and for an arbitrary epimorphism $\xi\colon\pi_1(X)\surj\Z$ 
in place of $\nu$. For such $\xi$ and $\lambda\in\C^\ast$, let $\rho_\lambda$ be the character
$g\mapsto\lambda^{\xi(g)}$. Then the maximal size of a $\lambda$-Jordan block of the 
monodromy on the torsion of $H_q(X^\xi;\C)$ is one less than the maximal length of a 
nonvanishing product $\langle\xi,\dots,\xi,a\rangle$ with $a\in H^q(X;\C_{\rho_\lambda})$, 
defined by chains as above \cite[Def.~2.1 and Thm.~5.1]{Pajitnov-ejm17}; these
products are the differentials of the equivariant spectral sequence twisted by $\rho_\lambda$, the
sequence $\E_r^\ast(\rho_\lambda,\xi)$ of \cite{Kohno-Pajitnov-cejm14}. (Pajitnov counts only the
entries equal to $\xi$, so his lengths are one less than ours.) Two consequences will be used
repeatedly: a formal $X$ has no Jordan block of size $>1$ at $\lambda=1$, and a strongly formal $X$
has none at any $\lambda$ \cite[Thm.~5.1]{Pajitnov-ejm17}. The dictionary holds as well with $\C_{\rho_\lambda}$
replaced by any finite-dimensional representation and $A_\nu$ by the corresponding twisted Alexander
module \cite[Thm.~C]{Pajitnov-topappl19}; see Proposition~\ref{prop:circle} below and
Remark~\ref{rem:higher-rank}.

The first two pages of $\E_r^\ast(\rho_\lambda,\xi)$ are explicit. Its $E_1$-term is
$H^\bullet(X;\C_{\rho_\lambda})$, and $d_1$ is the cup product with $\xi$ under the pairing
$H^1(X;\C)\otimes H^\bullet(X;\C_{\rho_\lambda})\to H^{\bullet+1}(X;\C_{\rho_\lambda})$
\cite[\S3]{Kohno-Pajitnov-cejm14}. So its $E_2$-term is the cohomology of the \emph{twisted Aomoto
complex} $\big(H^\bullet(X;\C_{\rho_\lambda}),\xi\cup-\big)$, which for $\lambda=1$ and $\xi=\nu$ is the
Aomoto complex of $\eta$ over $\C$. By \cite[Lem.~2.3]{Libgober-topappl02}, the twisted Aomoto complex
is the derivative complex of \cite[Def.~2.1]{Libgober-topappl02} at $\rho_\lambda$ in the direction
$\xi$, and its differential is the connecting homomorphism in cohomology of the extension of local
systems $0\to\C_{\rho_\lambda}\to L\to\C_{\rho_\lambda}\to0$, where $L$ is the local system with fiber
$\C[s]/(s^2)$ on which $g\in\pi_1(X)$ acts by multiplication by $\lambda^{\xi(g)}(1+s)^{\xi(g)}$.

Degeneration of these spectral sequences is a formality property in the literal sense. Let
$\bwedge(\eta)=H^\bullet(S^1;\k)$ be the exterior algebra on a generator $\eta$ of degree one. A \emph{DG
$\bwedge(\eta)$-module} is a cochain complex $(N,d)$ with an operator $\eta$ of degree one such that
$\eta^2=0$ and $d\eta+\eta d=0$; it is \emph{formal} if a chain of quasi-isomorphisms of DG
$\bwedge(\eta)$-modules joins it to $H(N)$, with zero differential and $\eta$ acting by the induced map
$\eta_\ast$. With $\eta$ the cocycle chosen at the beginning of this subsection, cup product with $\eta$
makes $C^\bullet(X;V)$ a DG $\bwedge(\eta)$-module for every finite-dimensional representation $V$ of
$\pi_1(X)$ over $\k$. It is the cocycle that acts, not its class. If $f\colon X\to S^1$ is the
simplicial map used to define $\eta$, then $C^\bullet(X;V)$ is a DG module over the cochain algebra
$C^\bullet(S^1;\k)$ through $f^\ast$, and the map $\bwedge(\eta)\to C^\bullet(S^1;\k)$ sending $\eta$ to
the cochain dual to the chosen edge, a cocycle of square zero, is a quasi-isomorphism of DG algebras;
the DG $\bwedge(\eta)$-module structure on $C^\bullet(X;V)$ is obtained by restriction along it. Cup product 
with another cocycle representing the same class is a different operator, and even when that cocycle has square
zero the resulting DG module need not be quasi-isomorphic to this one (Example~\ref{ex:three-lines}).
Over $\C$, multiplication by a closed $1$-form representing $\eta$ makes the de Rham complex
$A^\bullet(X,E_V)$ of the flat bundle with monodromy $V$ a DG $\bwedge(\eta)$-module, and there the
choice of the form does not matter (Remark~\ref{rem:restricted}). Let $N\llbracket s\rrbracket$ be the
complex of power series with coefficients in $N$ and differential $d+s\eta$, filtered by the powers of
$s$: it is the completed Koszul complex of \cite[\S3.2]{Suciu-arx26b} in one variable, with the DG
module $N$ in place of a CDGA. For $N=C^\bullet(X;\k)$, its spectral sequence is the equivariant 
spectral sequence of $(X,\eta)$, and by \cite[\S2]{Liu-Maxim-Wang-ejm26} the filtered complex
$N\llbracket s\rrbracket$ is isomorphic to $\Hom_R\big(C_\bullet(X^\nu;\k),\hat R\big)$, filtered by the
powers of $s$; so the spectral sequence, and with it the formality of the DG module $C^\bullet(X;\k)$ by
Proposition~\ref{prop:circle} below, depends only on $\nu$, not on the choices made in defining $\eta$. 
Chains in $N$ starting from a class $a\in H(N)$, and the products $\langle\eta,\dots,\eta,a\rangle$, 
are defined as for $C^\bullet(X;\k)$, with $d$ in place of $\delta$.

\begin{proposition}
\label{prop:circle}
Let $N$ be a DG $\bwedge(\eta)$-module over a field $\k$. The following are equivalent:
\begin{enumerate}[label=\textup{(\roman*)},leftmargin=2.2em,itemsep=2pt]
\item \label{pc1} $N$ is formal.
\item \label{pc2} Every class $a\in H(N)$ with $\eta_\ast a=0$ is the start of an infinite chain.
\item \label{pc3} The spectral sequence of the filtered complex $N\llbracket s\rrbracket$ degenerates at
$E_2$.
\end{enumerate}
Now let $V$ be a finite-dimensional representation of $\pi_1(X)$ over $\k$, let $\lambda\in\k^\ast$, and let
$V_\lambda=V\otimes\k_{\rho_\lambda}$ be the twist of $V$ by the character
$\rho_\lambda\colon g\mapsto\lambda^{\nu(g)}$. If $V=\k$ and $\lambda=1$, or if $\k\subseteq\C$, then for
$N=C^\bullet(X;V_\lambda)$ these conditions are also equivalent to:
\begin{enumerate}[label=\textup{(\roman*)},leftmargin=2.2em,itemsep=2pt,resume]
\item \label{pc4} The monodromy has no Jordan block of size greater than one for the eigenvalue $\lambda$
on the torsion of the twisted Alexander module $H_\ast(X^\nu;V^\vee)$, where
$V^\vee$ is the dual representation.
\end{enumerate}
\end{proposition}

\begin{proof}
\ref{pc1}$\Rightarrow$\ref{pc3}. A quasi-isomorphism $N'\to N$ of DG $\bwedge(\eta)$-modules induces a map
of filtered complexes $N'\llbracket s\rrbracket\to N\llbracket s\rrbracket$ which is an isomorphism on the first pages,
$H(N')[s]\to H(N)[s]$, and hence on all later ones; so the spectral sequence depends only on the
quasi-isomorphism type of $N$. For $H(N)$, with zero differential, the differential of $H(N)\llbracket s\rrbracket$ is
$s\eta_\ast$, so $d_1=\eta_\ast$ and $d_r=0$ for $r\ge2$.

\ref{pc3}$\Rightarrow$\ref{pc2}. Let $\eta_\ast a=0$. A cocycle $\omega_1$ representing $a$ is a $1$-chain
starting from $a$, and it extends to a $2$-chain, since $[\eta\,\omega_1]=\eta_\ast a=0$. Let $r\ge2$, and
let $\omega_1,\dots,\omega_r$ be an $r$-chain starting from $a$. The element
$c=\sum_{j=1}^r(-s)^{j-1}\omega_j$ of $N\llbracket s\rrbracket$ has leading term $\omega_1$, and
$(d+s\eta)c=-(-s)^r\eta\,\omega_r$. The class of $a$ survives to $E_r$, by \ref{pc3}, and so $d_r$ applied to
it is the class of $\pm\,\eta\,\omega_r$ in the term of filtration $r$. By \ref{pc3}, this differential
vanishes and $E_r=E_2$, whose term of filtration $r\ge1$ is $\ker\eta_\ast/\eta_\ast H(N)$. Hence
$[\eta\,\omega_r]=\eta_\ast[c']$ for some cocycle $c'$, and replacing $\omega_r$ by $\omega_r-c'$ gives an
$r$-chain with $\eta\,\omega_r$ exact, which therefore extends to an $(r+1)$-chain. Each step changes only
the last term, so induction on $r$ produces an infinite chain starting from $a$.

\ref{pc2}$\Rightarrow$\ref{pc1}. Write $H=H(N)$, and choose, in each degree, a complement $U$ of
$\ker\eta_\ast$ in $H$ and a complement $W$ of $\eta_\ast H$ in $\ker\eta_\ast$, so that
$H=U\oplus\eta_\ast U\oplus W$, with $\eta_\ast\colon U\to\eta_\ast U$ bijective. Let $Z$ be the DG
$\bwedge(\eta)$-module with basis $e_1,e_2,\dots$ in degree $0$ and $\eta e_1,\eta e_2,\dots$ in degree~$1$,
and with $de_1=0$ and $de_{j+1}=\eta e_j$; its cohomology is spanned by the class of $e_1$, on which $\eta$
acts by zero. Let $P$ be the direct sum of a free module $\bwedge(\eta)\,e_u$ on a generator $e_u$ of the
degree of $u$, with $de_u=0$, for each $u$ in a basis of $U$, and of a copy $Z_w$ of $Z$, shifted to the
degree of $w$, for each $w$ in a basis of $W$. Choose a cocycle $\omega_u$ representing $u$ and, by
\ref{pc2}, an infinite chain $(\omega_{w,j})_{j\ge1}$ starting from $w$. Then $e_u\mapsto\omega_u$ and
$e_{w,j}\mapsto\omega_{w,j}$ define a map of DG $\bwedge(\eta)$-modules $P\to N$, and $e_u\mapsto u$,
$e_{w,1}\mapsto w$, $e_{w,j}\mapsto0$ for $j\ge2$ define one $P\to H$. Both induce the isomorphism
$H(P)\cong U\oplus\eta_\ast U\oplus W=H$, so $N\leftarrow P\to H$ is a chain of quasi-isomorphisms.

Finally, let $N=C^\bullet(X;V_\lambda)$. As in \cite[\S2]{Liu-Maxim-Wang-ejm26}, the filtered complex
$N\llbracket s\rrbracket$ computes the cohomology of $X$ with coefficients in the local system
$V_\lambda\otimes_\k(L_\nu\otimes_R\hat R)$, in which $t$ acts on $\hat R$ as $1+s$, filtered by the
powers of $s$. When $V=\k$ and $\lambda=1$, degeneration of its spectral sequence is equivalent
to~\ref{pc4} by \cite[Thm.~1.1]{Liu-Maxim-Wang-ejm26}. When $\k=\C$, the twisted de Rham complex,
deformed by a closed $1$-form representing $\eta$, computes the cohomology with coefficients in the same
local system with $e^s$ in place of $1+s$, filtered in the same way \cite[Thm.~5.1]{Kohno-Pajitnov-cejm14};
the filtered automorphism $s\mapsto e^s-1$ of $\hat R$ identifies the two filtered local systems, so the
two spectral sequences agree from $E_1$ on, and degeneration is equivalent to~\ref{pc4} by
\cite[Thm.~5.1]{Pajitnov-ejm17} and \cite[Thm.~C]{Pajitnov-topappl19}. For $\k\subseteq\C$, both
conditions are invariant under extension of scalars to $\C$.
\end{proof}

\noindent In particular, for $N=C^\bullet(X;\k)$, condition~\ref{pc3} is condition~\ref{dg1} of the list at
the beginning of this subsection, so the conditions \ref{dg1}--\ref{dg3} there are also equivalent to
formality of the DG $\bwedge(\eta)$-module $C^\bullet(X;\k)$; and over $\C$, the twisted sequence
$\E_r^\ast(\rho_\lambda,\xi)$ degenerates at $E_2$ if and only if $A^\bullet(X,E_{\rho_\lambda})$ is 
formal as a DG module over the exterior algebra on $\xi$.

\begin{remark}
\label{rem:restricted}
The products $\langle\eta,\dots,\eta,\omega\rangle$ of this subsection are restricted forms of the
classical Massey products of $C^\bullet(X;\k)$. An $r$-chain starting from $\omega$ is a defining system, in
the classical sense, for the product of $r$ copies of $\eta$ followed by $\omega$, of a special kind: its
entries are the cocycle $\eta$ fixed above, the terms $\omega_j$ of the chain, and zero, the last for every
entry attached to a product of two or more copies of $\eta$, as is possible because $\eta\cup\eta=0$ on
cochains \cite[Def.~3.3]{Liu-Maxim-Wang-ejm26}. So the restricted product is contained in the classical
one, and its indeterminacy is smaller \cite[Rem.~2.2]{Pajitnov-ejm17}; for triple products it is
$\eta\cup H^q(X;\k)$, rather than $\eta\cup H^q(X;\k)+H^1(X;\k)\cup\omega$. The restriction comes from the
construction of the spectral sequence: an element $\sum_{j\ge1}(-s)^{j-1}\omega_j$ of $N\llbracket s\rrbracket$,
for $N=C^\bullet(X;\k)$, has coboundary in $s^rN\llbracket s\rrbracket$ exactly when $\omega_1,\dots,\omega_r$
is an $r$-chain, the cocycle $\eta$ being built into the differential $d+s\eta$. So only these defining
systems compute the differentials, and the ambiguity in their choice dies in $E_r$. Two consequences
point in opposite directions.

First, nonvanishing classical products along $\eta$ obstruct \ref{dg1}--\ref{dg3}. If the classical
product of $r\ge2$ copies of $\eta$ followed by $\omega\in H^q(X;\k)$ is defined and does not contain
$0$, then $\eta\cup\omega=0$, and the restricted product starting from $\omega$ of length $r+1$, or of
length $m+1$ if the longest chains starting from $\omega$ have length $m<r$, is defined and does not vanish. Its
length is at least three, as $\eta\cup\omega=0$, so the torsion of $H_q(X^\nu;\k)$ has a Jordan block
of size at least two at the eigenvalue $1$.

Second, the restricted products are invariants not of the DG algebra $C^\bullet(X;\k)$ and the class
$\eta$, but of the DG $\bwedge(\eta)$-module $C^\bullet(X;\k)$, which is defined through the map $X\to S^1$
and so depends on $\nu$. In characteristic zero the distinction disappears: on a CDGA model $A$ of $X$
\cite[Thms.~5.1 and~5.4]{Kohno-Pajitnov-cejm14}, multiplication by $e^{-sg}$ identifies the filtered
complexes $(A\llbracket s\rrbracket,d+s\theta)$ and $(A\llbracket s\rrbracket,d+s(\theta+dg))$ for every
closed $1$-form $\theta$ and every $g\in A^0$, since $d(e^{-sg})=-s\,dg\,e^{-sg}$; so the spectral sequence
depends only on the class of $\theta$, and it degenerates at $E_2$ when $A$ is formal, as in
\cite[Thm.~6.1]{Kohno-Pajitnov-cejm14}. In positive characteristic there is no such gauge, the exponential
needing denominators: the epimorphisms $\nu$ and $\nu+p\xi$, for $\xi\in H^1(X;\Z)$ such that $\nu+p\xi$ 
is onto, have the same reduction $\eta$ modulo $p$, but they may have different Jordan blocks, 
and the restricted products can fail to vanish when $X$ is
formal over $\F_p$, and even over $\Z$ (Example~\ref{ex:three-lines}). So, in positive characteristic,
the vanishing of the restricted products is a property of the pair $(X,\nu)$, not implied by formality
of the DG algebra $C^\bullet(X;\k)$. In the other direction, formality of the DG module does not imply
formality of the DG algebra, even in characteristic zero (Remark~\ref{rem:line}).
\end{remark}

\begin{example}
\label{ex:three-lines}
Let $M$ be the complement of the three lines of \S\ref{subsec:setting}, defined by $f=zw(z-w)$, and let
$\omega_1,\omega_2,\omega_3\in H^1(M;\Z)$ be the classes of $z$, $w$, and $z-w$, so that
$\nu=\omega_1+\omega_2+\omega_3$. The map $(z,w)\mapsto(z,w/z)$ identifies $M$ with
$\C^\ast\times(\C\smallsetminus\{0,1\})$, as in the proof of Theorem~\ref{thm:arr}: $\omega_1$ becomes the
class of the first factor, and $\bar\omega_2=\omega_2-\omega_1$ and $\bar\omega_3=\omega_3-\omega_1$
those of the two circles of $\C\smallsetminus\{0,1\}\simeq S^1\vee S^1$. So $M\simeq S^1\times(S^1\vee S^1)$,
and $\nu=3\omega_1+\bar\omega_2+\bar\omega_3$ is the epimorphism of Example~\ref{ex:toric-pencil} with
$m=3$. Over $\F_3$ we have $\eta=\bar\omega_2+\bar\omega_3$, and, since
$\bar\omega_2\bar\omega_3=0$, each $\eta\cup\omega_i$ equals $-\omega_1\eta$; so the kernel of $\eta\cup-$
on $H^1(M;\F_3)$ is spanned by $\bar\omega_2$ and $\bar\omega_3$, and $\beta^3_1(M,\eta)=\beta^3_2(M,\eta)=1$.

\emph{Formality.} The space $M$ has the homotopy type of the toric complex of
Example~\ref{ex:toric-pencil} with $m=3$, hence is formal over every commutative ring
\cite[Thm.~4.8 and Rem.~4.9]{Notbohm-Ray-agt05}, and so is $F\simeq\bigvee^4S^1$; thus every classical
Massey product of $M$ over $\F_3$ contains $0$.

\emph{Restricted products.} By Example~\ref{ex:toric-pencil}, $H_1(F;\F_3)$ has a Jordan block 
of size three at the eigenvalue $1$, so, by \cite[Thm.~1.1]{Liu-Maxim-Wang-ejm26}, some product
$\langle\eta,\eta,\eta,x\rangle$ with $x\in H^1(M;\F_3)$ does not vanish. To see it on a cochain model, 
give the circle factor vertices $v_0,v_1,v_2$ and edges $\sigma_1,\sigma_2,\sigma_3$, 
with $\sigma_i$ running from $v_{i-1}$ to $v_i$, indices modulo $3$, and give the wedge one 
vertex and edges $\alpha$ and $\gamma$, representing $\bar\omega_2$ and $\bar\omega_3$; up to
homotopy, $f$ maps every edge once around the circle. For a graph whose edges are mapped once around
the circle, the complex of this subsection is, directly from the definitions, the cochain complex of the
infinite cyclic cover, $\eta$ being the sum of the duals of the edges. The infinite cyclic cover of $M$ is
the quotient of the product of those of the two factors by the antidiagonal action of $\Z$, so its
cochain complex is the tensor product of theirs. Hence the spectral sequence of this subsection, and
with it the products $\langle\eta,\dots,\eta,x\rangle$, may be computed on the tensor product of the
cochain algebras of the factors, with $\eta=\eta_3\otimes1+1\otimes(\alpha^\ast+\gamma^\ast)$, where
$\eta_3=\sigma_1^\ast+\sigma_2^\ast+\sigma_3^\ast$. Over $\F_3$ we have $\eta_3=\delta\phi$ with
$\phi=v_1^\ast+2v_2^\ast$, then $\eta_3\cup\phi=\delta\psi$ with $\psi=v_1^\ast$, and
$\eta_3\cup\psi=\sigma_1^\ast$. Hence $(1\otimes\alpha^\ast,\phi\otimes\alpha^\ast,
\psi\otimes\alpha^\ast)$ is a $3$-chain starting from $\bar\omega_2$, and
$\eta\cup(\psi\otimes\alpha^\ast)=\sigma_1^\ast\otimes\alpha^\ast$ represents 
$\omega_1\bar\omega_2=\omega_1\omega_2$, which is not in $\eta\cup H^1(M;\F_3)$, 
the span of $\omega_1\eta$. This span is the whole indeterminacy of the product, because the triple
products $\langle\eta,\eta,\bar\omega_2\rangle$ and $\langle\eta,\eta,\bar\omega_3\rangle$ contain $0$: the
chain above extends, and so does its analogue with $\gamma$ in place of $\alpha$. So
$\langle\eta,\eta,\eta,\bar\omega_2\rangle$ does not vanish, and 
likewise for every $x$ in the kernel of $\eta\cup-$ that is not a multiple of $\eta$. 
The source is the circle factor, on which $f$ has degree $3$: there $[\eta_3]=0$ in 
$H^1(S^1;\F_3)$, but $\langle\eta_3,\eta_3,\eta_3,1\rangle=[\sigma_1^\ast]\ne0$.

\emph{Dependence on $\nu$.} The epimorphism $\nu-3\omega_1=\bar\omega_2+\bar\omega_3$ has the same
reduction $\eta$ modulo $3$, and its cocycle $\eta'=1\otimes(\alpha^\ast+\gamma^\ast)$, cohomologous to
$\eta$ and also of square zero, has zero products with $1\otimes\alpha^\ast$ and $1\otimes\gamma^\ast$.
So the chains starting from $\bar\omega_2$ and $\bar\omega_3$ continue by zeros, every product of length at
least three vanishes, and $(M,\nu-3\omega_1)$ satisfies \ref{dg1}--\ref{dg3} over $\F_3$; indeed it is
Alexander formal, by Proposition~\ref{prop:toric}\ref{to2}. The structures of DG $\bwedge(\eta)$-module on
the cochains of $M$ over $\F_3$ given by cup product with $\eta$ and with $\eta'$ are therefore not
quasi-isomorphic: the class of $\eta$ alone does not determine the action. Since a classical defining
system may use any cocycle representing $\eta$, the same computation with $\eta'$ also shows directly that
the classical product $\langle\eta,\eta,\eta,\bar\omega_2\rangle$ contains $0$.

The same conclusions hold for the pencil of $n\ge3$ lines through the origin, whose complement has the
homotopy type of $S^1\times\bigvee^{n-1}S^1$, at every prime $p$ dividing $n$, with
$\beta^p_1(M,\eta)=\beta^p_2(M,\eta)=n-2$, with $\nu-n\omega_1$ in place of $\nu-3\omega_1$, and with
Jordan blocks of size $p^s$, where $p^s$ is the exact power of $p$
dividing $n$ (Remark~\ref{rem:arr-collision}).
\end{example}

For a root of unity $\zeta$ and $q\ge0$ we write $b_q(X,\zeta)=\dim_\C H^q(X;\C_{\rho_\zeta})$, where
$\rho_\zeta$, with $\xi=\nu$, is the character obtained from $\nu$ and $\zeta$, so that $b_q(X,1)=b_q(X)$;
Galois conjugate roots of unity give the same value. If $\zeta\ne1$ and $H_1(X^\nu;\C)$ is
finite-dimensional, then $b_1(X,\zeta)$ is the number of Jordan blocks of $h_\ast$ for the
eigenvalue $\zeta$ on $H_1(X^\nu;\C)$, by the exact sequence of the cover, as in
\cite[Lem.~10.4]{Suciu-planes-r4}. So $b_1(X,\zeta)$ is
at most the algebraic multiplicity of $\zeta$ as an eigenvalue of $h_\ast$ on $H_1(X^\nu;\C)$, with equality
if and only if the monodromy is semisimple at $\zeta$, as it is when the monodromy has finite order
(compare \cite[Cor.~6.15]{Dimca-Papadima-Suciu-imrn08}).

\subsection{The two translations}
\label{subsec:translations}
The choice of field and twist selects the part of $A_\nu$ that is constrained.

\emph{At an eigenvalue $\zeta$, in characteristic zero.} Twisting the coefficients by $\rho_\zeta$,
that is, replacing $t$ by $\zeta t$, carries the $(t-\zeta)$-primary part of $A_\nu\otimes\Q(\zeta)$ to
the $(t-1)$-primary part. Degeneration of the twisted sequence over $\Q(\zeta)$ is therefore a
condition at the horizontal prime $(\Phi_d)$ of $\Spec\LA$, where $d$ is the order of $\zeta$.

\emph{At a prime $p$.} Over $\F_p$ the $(t-1)$-primary part of $H_\ast(X^\nu;\F_p)$ contains the
reductions of all the $(t-\zeta)$-primary parts with $\zeta$ of $p$-power order, since
$\Phi_{p^s}\equiv(t-1)^{\varphi(p^s)}\bmod p$, together with the contribution of the $p$-torsion
of $A_\nu$. Degeneration of the untwisted sequence over $\F_p$ is therefore a condition at the closed
point $(p,t-1)$, where all of these meet. The truncations of this sequence compute the mod-$p$ 
homology of the cyclic covers $X^{(p^r)}$ of \S\ref{subsec:setting}, and
degenerate at the latest at $E_{p^r}$ \cite[\S4]{Liu-Maxim-Wang-ejm26}.

The meeting at $(p,t-1)$ has consequences: the monodromy may be semisimple at a rational 
eigenvalue of order $p^s>2$, but, as the next lemma shows, its reduction modulo $p$ is not.

\begin{lemma}
\label{lem:collision}
Let $L$ be a nonzero $\LA$-module that is finitely generated and torsion-free as an abelian group,
and suppose that $\Phi_q(t)L=0$ for a prime power $q=p^a>2$. Then $t-1$ acts on $L/pL$ as a nonzero
nilpotent endomorphism. Consequently, if $H_i(X^\nu;\Z)$ is finitely generated as an abelian group
and a primitive $q$-th root of unity is an eigenvalue of $h_\ast$ on $H_i(X^\nu;\Q)$, then the 
$(t-1)$-primary part of $H_i(X^\nu;\F_p)$ is not semisimple.
\end{lemma}

\begin{proof}
Since $\Phi_q\equiv(t-1)^{\varphi(q)}\bmod p$, the endomorphism $t-1$ of $L/pL$ is nilpotent. The
group $L$ is a module over $\Z[\zeta_q]=\LA/(\Phi_q)$, torsion-free because it is torsion-free over
$\Z$, and in $\Z[\zeta_q]$ one has $(1-\zeta_q)^{\varphi(q)}=pu$ with $u$ a unit; hence
$(t-1)^{\varphi(q)}L=pL$. If $t-1$ vanished on $L/pL$, then $(t-1)L\subseteq pL$, whence
$pL=(t-1)^{\varphi(q)}L\subseteq p^{\varphi(q)}L\subseteq p^2L$, since $\varphi(q)\ge2$; as $L$ is
torsion-free this gives $L=pL$, which is impossible for a nonzero finitely generated free abelian
group.

For the second assertion, let $H'$ be the quotient of $H_i(X^\nu;\Z)$ by its torsion subgroup, a
free abelian group of finite rank, and let $L=\ker\big(\Phi_q(t)\colon H'\to H'\big)$. Then
$L\ne0$, and $L$ is saturated in $H'$, because $H'/L\cong\Phi_q(t)H'$ is torsion-free; hence
$L/pL\subseteq H'/pH'$, and by the first part $t-1$ is nonzero on the $(t-1)$-primary part of
$H'/pH'$. This part is a quotient of the $(t-1)$-primary part of $H_i(X^\nu;\Z)\otimes\F_p$, since
$H_i(X^\nu;\Z)\otimes\F_p\to H'/pH'$ is surjective and taking primary parts of finite-dimensional
$\F_p[t]$-modules is exact. Finally, $H_i(X^\nu;\Z)\otimes\F_p$ is a $\LA$-submodule of
$H_i(X^\nu;\F_p)$ by the universal coefficient theorem. So $t-1$ is nonzero on the $(t-1)$-primary
part of $H_i(X^\nu;\F_p)$, on which $1$ is the only eigenvalue, and this part is not semisimple.
\end{proof}

\begin{remark}
\label{rem:collision}
The prime power $q=2$ must be excluded: $\Phi_2=t+1\equiv t-1\bmod 2$, and $\LA/(t+1)$
reduces to $\F_2$ with trivial action. If $p\nmid d$, then $\Phi_d$ is squarefree modulo $p$, and
its reduction has the root $1$ only when $d=1$. If $d=p^am$ with $a\ge1$, $p\nmid m$, and $m>1$,
then $\Phi_d\equiv\Phi_m^{\varphi(p^a)}\bmod p$, which has a repeated factor when $p^a>2$, but away
from $t-1$, at closed points $(p,g)$ with $g\ne t-1$ that Definition~\ref{def:Aformal} below does
not examine (Remark~\ref{rem:scope}). So, one cyclotomic factor at a time, the eigenvalues
whose reduction modulo $p$ produces a Jordan block of size greater than one at $t=1$ are exactly
those of order $p^a>2$.
Across several factors the lattice matters as well: $\LA/(t^2-1)$ has semisimple rational
monodromy, with eigenvalues $\pm1$, but reduces modulo $2$ to the single Jordan block
$\F_2[t]/(t-1)^2$.
\end{remark}

\section{Alexander formality}
\label{sec:def}

We now define the notion. Definition~\ref{def:etak} isolates the condition attached to a single field
and a single root of unity, and Definition~\ref{def:Aformal} imposes it at every root of unity in
characteristic zero and at every prime. The remarks that follow record what the definition 
depends on and how it reads as a formality condition (Remark~\ref{rem:intrinsic}), which points
of $\Spec\LA$ it examines (Remark~\ref{rem:scope}), and in which degrees (Remark~\ref{rem:degree}).

\begin{definition}
\label{def:etak}
Let $\k$ be a field and $\zeta\in\k^\ast$ a root of unity. The space $X$ is \emph{\etaform{} over
$\k$ at $\zeta$} if the equivariant spectral sequence of $(X,\eta)$ with coefficients twisted by the
character $\rho_\zeta\colon g\mapsto\zeta^{\nu(g)}$, that is, the spectral sequence of the filtered complex
$C^\bullet(X;\k_{\rho_\zeta})\llbracket s\rrbracket$ of \S\ref{subsec:ss}, degenerates at $E_2$ in every
degree; for
$\zeta=1$ we say simply that $X$ is \emph{\etaform{} over $\k$}. We use this only for $\zeta=1$
and for $\k=\Q(\zeta)\subseteq\C$, where the following equivalences are available: by
\cite[Thm.~1.1]{Liu-Maxim-Wang-ejm26} for $\zeta=1$ and by \cite[Thm.~5.1]{Pajitnov-ejm17} 
for $\k\subseteq\C$, the condition holds if and only if all Massey products
$\langle\eta,\dots,\eta,a\rangle$ of length at least three with $a\in H^\bullet(X;\k_{\rho_\zeta})$
vanish, and if and only if the $(t-\zeta)$-primary part of the torsion of
$H_\ast(X^\nu;\k)$ is semisimple. Restricting the last two conditions to $a\in H^q(X;\k_{\rho_\zeta})$ and
to $H_q(X^\nu;\k)$, which keeps them equivalent, we say that $X$ is \emph{\etaform{} over $\k$ at
$\zeta$ in degree~$q$}. (Pajitnov works over $\C$, with twisted de Rham complexes, 
whose spectral sequences agree with these
from $E_1$ on by the proof of Proposition~\ref{prop:circle}; for $\k=\Q(\zeta)$ all three conditions are
invariant under extension of scalars to $\C$.) The Massey products are those
of \S\ref{subsec:ss}, defined by chains from a fixed representative of $\eta$; in positive characteristic
they are computed in $C^\bullet(X;\k)$, and they need not vanish when $X$ is formal over $\k$
(Remark~\ref{rem:restricted}).
\end{definition}

\begin{definition}
\label{def:Aformal}
A space $X$ is \emph{Alexander formal} if it is \etaform{} over $\Q(\zeta)$ at $\zeta$ for every root of unity
$\zeta$, and \etaform{} over $\F_p$ for every prime $p$. Equivalently: for every $d\ge1$ the
$\Phi_d$-primary part of the torsion of $A_\nu\otimes\Q$ is semisimple, and for every prime $p$ the
$(t-1)$-primary part of the torsion of $H_\ast(X^\nu;\F_p)$ is semisimple. In the language of
\S\ref{sec:setup}, the equivariant spectral sequence degenerates at $E_2$ at the horizontal primes
$(\Phi_d)$ of $\Spec\LA$ and at the closed points $(p,t-1)$.
\end{definition}

\begin{remark}
\label{rem:intrinsic}
The characteristic-zero clauses of Definition~\ref{def:Aformal} depend only on $A_\nu\otimes\Q$, and
the modular ones on the $\F_p[t^{\pm1}]$-modules $H_\ast(X^\nu;\F_p)$; these are determined by the
$\LA$-module $A_\nu$ when $A_\nu$ has no $p$-torsion, $H_q(X^\nu;\F_p)$ being then
$H_q(X^\nu;\Z)\otimes\F_p$, and in general by the $\LA$-chain complex $C_\ast(X^\nu)$. 

The clauses are nonetheless formality conditions in the literal sense. By
Proposition~\ref{prop:circle}, $X$ is \etaform{} over $\F_p$, or over $\Q(\zeta)$ at $\zeta$, if and only if
$C^\bullet(X;\F_p)$, or $C^\bullet(X;\Q(\zeta)_{\rho_\zeta})$, is formal as a DG module over the exterior
algebra $H^\bullet(S^1)$, acting by cup product with the cocycle $\eta$ of \S\ref{subsec:ss}. The twisted
coefficients can be replaced by finite covers. The infinite cyclic cover of the pair
$(X^{(d)},\nu^{(d)})$ of \S\ref{subsec:setting} is $X^\nu$, with deck transformation $t^d$, and in
characteristic zero $t^d$ has a Jordan block of size greater than one at the eigenvalue $1$ if and only
if $t$ has one at some $d$-th root of unity. 
So the clauses of Definition~\ref{def:Aformal} at the $d$-th roots of unity say together that $(X^{(d)},\nu^{(d)})$ 
is \etaform{} over $\Q$, and $X$ is Alexander formal if and only if $C^\bullet(X^{(d)};\Q)$ is formal as a
DG module over $H^\bullet(S^1;\Q)$ for every $d\ge1$ and $C^\bullet(X;\F_p)$ is formal as a DG module over
$H^\bullet(S^1;\F_p)$ for every prime $p$: one formality property, imposed on the rational tower 
of cyclic covers along $\nu$ and on $X$ modulo every prime.
\end{remark}

\begin{remark}
\label{rem:scope}
Definition~\ref{def:Aformal} imposes degeneration at the horizontal primes $(\Phi_d)$ and at the
closed points $(p,t-1)$. It omits two kinds of points of $\Spec\LA$ that may support the torsion of $A_\nu$. First,
the eigenvalue clause ranges over roots of unity, so it accounts for
the whole support of the torsion of $A_\nu\otimes\Q$ exactly when $h_\ast$ is quasi-unipotent. This holds
whenever $X$ is a smooth complex quasi-projective variety or a compact K\"ahler
manifold \cite[Prop.~1.4]{Budur-Liu-Wang-mathann18}, by the structure theorems for their cohomology
jump loci \cite{Budur-Wang-ens15,Wang-mrl16}, and for most of the examples below, but not for a general
fibered $X$. Second, at a prime $p$ the definition imposes degeneration only at $(p,t-1)$; degeneration 
at the remaining closed points $(p,g)$, with $g\in\F_p[t^{\pm1}]$ irreducible and $g\ne t-1$, is a
further requirement, not imposed here. Both omissions are deliberate, and both are realized by 
examples. The mapping torus of Proposition~\ref{prop:converse} is Alexander formal, yet its monodromy 
has a Jordan block of size two at the non-cyclotomic horizontal prime $(t^2-3t+1)$. A bielliptic surface of 
type $\Z/6$ is Alexander formal, yet $\Phi_6\equiv(t+1)^2\bmod3$ gives a Jordan block at the closed point 
$(3,t+1)$ with no rational counterpart (Example~\ref{ex:bielliptic}). We keep the definition as stated because 
the arithmetic content of Alexander formality, the collision of Lemma~\ref{lem:collision}, sits at $(p,t-1)$, 
and because the spaces it is designed for have quasi-unipotent monodromy; the natural extension would impose 
degeneration at every point of $\Spec\LA$ that supports the torsion of $A_\nu$.
\end{remark}

\begin{remark}
\label{rem:degree}
The classical invariants $\beta^p_1(X,\eta)$ and $b_1(X,\zps)$ are degree-one invariants, whereas the
clauses of Definition~\ref{def:Aformal} concern all degrees. A comparison of the two can
therefore detect a failure of Alexander formality, but never certify Alexander formality. When $X^\nu$ has
finite-dimensional homology, Proposition~\ref{prop:finite-type} reduces the modular clauses to the
acyclicity of the Aomoto complex in all degrees, and the characteristic-zero clauses to that of the
twisted Aomoto complexes of \S\ref{subsec:ss}; both are computable (Corollary~\ref{cor:finite-type}).
For arrangement complements the mod-$p$ Aomoto complex is often exact in degree one
(Example~\ref{ex:betti}). 
In characteristic zero, the torsion of $A_\nu\otimes\Q$ is semisimple in every degree when $h_\ast$ 
has finite order and when $X$ is strongly formal (Proposition~\ref{prop:sf-char0}), in particular 
when $X$ is compact K\"ahler \cite[Thm.~1.6]{Budur-Liu-Wang-mathann18}; for a formal $X$ 
its eigenvalue-$1$ part is semisimple in every degree \cite[Prop.~1.7]{Budur-Liu-Wang-mathann18}.
Finally, for $X$ smooth quasi-projective of dimension $n$, the Jordan blocks of $h_\ast$ on $H_q$, at
every eigenvalue, have size at most $\min\{q+1,2n-q\}$ \cite[Thm.~1.8]{Budur-Liu-Wang-mathann18}; so in
degree one, and in characteristic zero, all products $\langle\eta,\dots,\eta,a\rangle$ of length at
least four vanish.
\end{remark}

\section{The classical constraints}
\label{sec:first}

This section recalls the classical constraints on the algebraic monodromy that Alexander formality
assembles (\S\ref{subsec:constraints}), and answers the first clause of Question~\ref{q:main}
(\S\ref{subsec:assembly}). Proposition~\ref{thm:assembly} relates the clauses at $\lambda=1$ 
and at a prime to these constraints, Proposition~\ref{prop:finite-type} reduces each clause to 
the acyclicity of an Aomoto complex when $X^\nu$ has finite-dimensional homology, 
and Corollary~\ref{cor:assembly} collects the answer.

\subsection{Four constraints on the monodromy}
\label{subsec:constraints}
Formality properties of $X$ constrain $h_\ast$, and several such constraints are known. We list the
four that motivate this paper. Each is an instance of the degeneration described in \S\ref{subsec:ss},
taken over a different field or after a different twist.
\begin{enumerate}[label=\textup{(\arabic*)},leftmargin=2.2em,itemsep=4pt]
\item \label{f1}
\emph{Rationally, at the eigenvalue $\lambda=1$.} Recall from \S\ref{subsec:setting} that formality
of $X$ implies $1$-formality of $\pi_1(X)$. If $\pi_1(X)$ is $1$-formal and $H_1(X^\nu;\Q)$ is
finite-dimensional, then $h_\ast$ has no $(t-1)$-primary Jordan block of size greater than one on
$H_1(X^\nu;\Q)$, that is to say, $X$ is $\eta$-formal over $\Q$ in degree one. This is a theorem of
Papadima and Suciu \cite[Prop.~3.2]{Papadima-Suciu-forum10}, proved by means of the tangent cone
theorem for cohomology jump loci \cite{Dimca-Papadima-Suciu-duke09}.

\item \label{f2}
\emph{Modularly, at the prime-power characters.} Let $p$ be a prime, and let $\beta^p_q(X,\eta)$ and
$b_q(X,\zeta)$ be the Aomoto--Betti and twisted Betti numbers of \S\ref{subsec:ss}. Comparing the
$E_2$-page of the mod-$p$ spectral sequence with its limit bounds these twisted Betti numbers by the
mod-$p$ Aomoto--Betti numbers: if $H_\ast(X;\Z)$ is torsion-free, then
\[
b_1(X,\zps)\le\beta^p_1(X,\eta) \qquad (s\ge1) ,
\]
and more generally $\dim_\C H^q(X;\C_\rho)\le\beta^p_q(X,\eta)$ in every degree $q$ and for every
character $\rho$ of order $p^s$ obtained from $\nu$ \cite[Thm.~C]{Papadima-Suciu-tams10}. The
torsion-freeness hypothesis cannot be dropped \emph{loc.\ cit.}; for the horizontal arrangements of
\S\ref{subsec:diagonal} the bound reads $b_1(X,\zps)\le\beta^p_1(X,\eta)\le v_p(\tau)$ whenever the
signed spanning-tree number $\tau$ of the linking graph is nonzero \cite[Thm.~E]{Suciu-planes-r4}.
Sharpness of all these bounds, at every character of $p$-power order and in every degree, is equivalent
to $\eta$-formality
of $X$ over $\F_p$, provided the $p$-power cyclic covers of $X$ along $\nu$ have no $p$-torsion in
their homology; the defects are accounted for by the higher differentials of the mod-$p$ spectral
sequence (Proposition~\ref{thm:assembly}\ref{a2}). Sharpness at the characters of order $p$ in degree
one alone is a much weaker condition (Example~\ref{ex:betti}).

\item \label{f3}
\emph{At all characters, over $\C$.} By \cite[Thm.~6.3]{Kohno-Pajitnov-cejm14}, strong
formality (\S\ref{subsec:KP}) implies that for every $\rho\in\CC=H^1(X;\C^\ast)$ and every
$\alpha\in H^1(X;\C)$ the twisted spectral sequence $\E^\ast_r(\rho,\alpha)$ of that paper, whose
differentials are the twisted Massey products $\langle\alpha,\dots,\alpha,-\rangle$ on
$H^\bullet(X;\C_\rho)$, degenerates at its second page. In view of Pajitnov's theorem quoted in
\S\ref{subsec:ss}, a strongly formal space has no Jordan block of size greater than one at any
eigenvalue (Proposition~\ref{prop:sf-char0}).

\item \label{f4}
\emph{Integrally.} The bound in \ref{f2} can be strict, and then $A_\nu$ carries formality-theoretic
information that is lost over $\Q$. The complement $X$ of Yoshinaga's icosidodecahedral arrangement
of $16$ planes in $\C^3$ has $\beta^2_1(X,\eta)=1$ while $b_1(X,-1)=0$
\cite[Prop.~4.2, Thm.~4.3]{Yoshinaga-ejm20}. It is formal, being an arrangement complement, and the
algebraic monodromy of its Milnor fibration is trivial on $H_1(F;\Q)$; yet $X$ is not $\eta$-formal
over $\F_2$ (Theorem~\ref{thm:arr}),
and the strict inequality is accounted for by $2$-torsion in $A_\nu$, the first homology of the
Milnor fiber being $\Z^{15}\oplus\Z/2$ \cite{Yoshinaga-ejm20}, \cite[\S11.3]{Suciu-revroum24}. 
Further examples of this kind are in \cite{Ishibashi-Sugawara-Yoshinaga-advappl25}.
\end{enumerate}

Two asymmetries in the list are worth noting. First, formality of $X$ yields
\ref{f1}, indeed in every degree \cite[Thm.~6.1]{Kohno-Pajitnov-cejm14}, but nothing modular: the
space in \ref{f4} is formal, and the bound in \ref{f2} is strict for it at $p=2$; nor does formality
over $\F_p$, or even over $\Z$ (Example~\ref{ex:three-lines}). Strong formality yields the whole
characteristic-zero picture (Proposition~\ref{prop:sf-char0}), but again nothing modular. Formality of
DG modules over the cohomology of the circle, by contrast, makes sense over every field, and over
$\F_p$ it is $\eta$-formality (Proposition~\ref{prop:circle}). Second, conditions \ref{f1}, \ref{f2},
and \ref{f4} concern the $\LA$-module $A_\nu$ and the Massey products along the single class $\eta$,
whereas \ref{f3} involves the whole family of twisted models indexed by $\CC$. 

\subsection{Assembling the constraints}
\label{subsec:assembly}
We first split Alexander formality into its clauses, and relate the clause at $\lambda=1$ to
condition~\ref{f1} and the clause at a prime to condition~\ref{f2}.

\begin{proposition}
\label{thm:assembly}
\begin{enumerate}[label=\textup{(\roman*)},leftmargin=2.2em,itemsep=2pt]
\item \label{a1}
If $H_1(X^\nu;\Q)$ is finite-dimensional, $1$-formality of $\pi_1(X)$ implies $\eta$-formality over
$\Q$ in degree one; formality of $X$ implies it in every degree. Neither implication can be reversed.
\item \label{a2}
Let $p$ be a prime, let $X^{(p^r)}\to X$ be the $p^r$-fold cyclic cover along $\nu$, and let $b_q(X,\zeta)$
be the twisted
Betti numbers of \S\ref{subsec:ss}. If $H_\ast(X^{(p^r)};\Z)$ has no $p$-torsion for every $r\ge0$, then
$X$ is \etaform{} over
$\F_p$ if and only if $b_q(X,\zps)=\beta^p_q(X,\eta)$ for all $q\ge0$ and all $s\ge1$.
\end{enumerate}
\end{proposition}

\begin{proof}
\ref{a1} The first implication is \cite[Prop.~3.2]{Papadima-Suciu-forum10}, whose hypothesis that
$b_1(\ker\nu)$ be finite is the finite-dimensionality of $H_1(X^\nu;\Q)$. The second implication is
\cite[Thm.~6.1]{Kohno-Pajitnov-cejm14}.
For the failure of the converses, let $H$ be the Heisenberg nilmanifold, $X=S^1\times H$, and $\nu$
the projection onto $\pi_1(S^1)$. Then $X^\nu\simeq H$ and $t$ acts trivially on
$H_\ast(X^\nu;\Z)$, so $X$ is $\eta$-formal over every field, while $\pi_1(X)=\Z\times\pi_1(H)$ is
not $1$-formal.

\ref{a2} For $r\ge1$, the characters of $\pi_1(X)$ that factor through $\nu$ and have order dividing
$p^r$ are the trivial one and, for $1\le s\le r$, the $\varphi(p^s)$ characters of order $p^s$, which
are Galois conjugate to one another. Hence $b_q(X^{(p^r)};\C)=b_q(X)+\sum_{s=1}^r\varphi(p^s)\,
b_q(X,\zps)$, and in the absence of $p$-torsion $b_q(X^{(p^r)};\F_p)=b_q(X^{(p^r)};\C)$ and
$b_q(X;\F_p)=b_q(X)$. On the other hand, by \cite[Props.~1.3, 1.4]{Liu-Maxim-Wang-ejm26},
\[
b_q(X^{(p^r)};\F_p)\le b_q(X;\F_p)+(p^r-1)\,\beta^p_q(X,\eta)\qquad\text{for all $q$},
\]
with equality for all $q$ if and only if the truncation of the mod-$p$ equivariant spectral sequence
computing $H_\ast(X^{(p^r)};\F_p)$ degenerates at $E_2$ (which is automatic for $p^r=2$); the differentials of
that truncation are given by the same Massey products \cite[\S4]{Liu-Maxim-Wang-ejm26}, and the
inequality is strict in the degree of $\omega$ whenever some product $\langle\eta,\dots,\eta,\omega\rangle$
of length $k$ with $3\le k\le p^r$ does not vanish. Since $p^r-1=\sum_{s=1}^r\varphi(p^s)$, the
inequality reads
\[
\sum_{s=1}^r\varphi(p^s)\big(\beta^p_q(X,\eta)-b_q(X,\zps)\big)\ \ge\ 0 .
\]
If $X$ is \etaform{} over $\F_p$, all products $\langle\eta,\dots,\eta,\omega\rangle$ of length at
least three vanish, every truncation degenerates, equality holds for every $r$, and subtracting the
equalities for $r-1$ and $r$ gives $b_q(X,\zeta_{p^r})=\beta^p_q(X,\eta)$ for every $r\ge1$.
Conversely, if
$b_q(X,\zps)=\beta^p_q(X,\eta)$ for all $q$ and $s$, equality holds for every $r$; a nonvanishing product
of length $k\ge3$ would make it strict for any $r$ with $p^r\ge k$, so there is none, and $X$ is
\etaform{} over $\F_p$.
\end{proof}

When the infinite cyclic cover has finite-dimensional homology, $\eta$-formality becomes a condition
on an Aomoto complex alone: over an arbitrary field $\k$ at $\zeta=1$, and over $\C$, twisted by
$\rho_\zeta$, at a root of unity $\zeta$.

\begin{proposition}
\label{prop:finite-type}
\begin{enumerate}[label=\textup{(\roman*)},leftmargin=2.2em,itemsep=2pt]
\item \label{ft1}
Let $\k$ be a field such that $H_q(X^\nu;\k)$ is finite-dimensional for every $q$. Then $X$ is
\etaform{} over $\k$ if and only if the Aomoto complex $\big(H^\bullet(X;\k),\eta\cup-\big)$ is
acyclic, that is, $\beta^p_q(X,\eta)=0$ for all $q$, where $p$ is the characteristic of $\k$.
\item \label{ft2}
Suppose that $H_q(X^\nu;\Q)$ is finite-dimensional for every $q$. Then, for every root of unity
$\zeta$, the space $X$ is \etaform{} over $\Q(\zeta)$ at $\zeta$ if and only if the twisted Aomoto complex
$\big(H^\bullet(X;\C_{\rho_\zeta}),\eta\cup-\big)$ is acyclic.
\end{enumerate}
\end{proposition}

\begin{proof}
\ref{ft1} This is essentially \cite[Prop.~9.4]{Papadima-Suciu-tams10}, which shows moreover that
acyclicity of the Aomoto complex forces $H_\ast(X^\nu;\k)$ to be finite-dimensional; we include the
short argument. The equivariant spectral sequence of \S\ref{sec:setup} converges to $H^\bullet(X;L_\nu\otimes_R\hat
R)$, filtered by the powers of $s$. The complex computing this cohomology is the $\hat R$-dual of
$C_\ast(X^\nu;\k)\otimes_R\hat R$, a complex of free $\hat R$-modules whose homology is
$H_\ast(X^\nu;\k)\otimes_R\hat R$, since $\hat R$ is flat over $R$
\cite[\S2]{Liu-Maxim-Wang-ejm26}. By hypothesis this homology has finite length, so, $\hat R$ being
a principal ideal domain, the universal coefficient theorem shows that each
$H^m\coloneqq H^m(X;L_\nu\otimes_R\hat R)$ has finite length too. Multiplication by $s^i$ identifies
the complex with its $i$-th filtered piece, so the filtration induced on $H^m$ is the $s$-adic one,
and $E_\infty^{i,m-i}=s^iH^m/s^{i+1}H^m$ vanishes for $i\gg0$. If the sequence degenerates at $E_2$, then
$\beta^p_m(X,\eta)=\dim E_2^{i,m-i}=\dim E_\infty^{i,m-i}=0$ for $i\gg0$. Conversely, if the Aomoto
complex is acyclic, then $E_2^{i,j}=0$ for $i\ge1$, and each differential $d_r$ with $r\ge2$, which
raises the column index by $r$, has zero target; so the sequence degenerates at $E_2$.

\ref{ft2} By Definition~\ref{def:etak}, the condition is that the spectral sequence of the filtered
complex $C^\bullet(X;\C_{\rho_\zeta})\llbracket s\rrbracket$, with differential $\delta+s\eta$, degenerate
at $E_2$. As in \cite[\S2]{Liu-Maxim-Wang-ejm26}, this complex is the $\hat R$-dual of
$C_\ast(X^\nu;\C)\otimes_R\hat R$, where now $R=\C[t^{\pm1}]$ acts on $\hat R=\C\llbracket s\rrbracket$
through $t\mapsto\mu(1+s)$, with $\mu=\zeta$ or $\mu=\zeta^{-1}$ according to convention. Thus $\hat R$
is the completion of $R$ at $t=\mu$, and the homology of $C_\ast(X^\nu;\C)\otimes_R\hat R$ is the
completion of $H_\ast(X^\nu;\C)=H_\ast(X^\nu;\Q)\otimes_\Q\C$ at $t=\mu$, which has finite length. The
argument of part~\ref{ft1} now applies verbatim, the $E_2$-term in positive filtration being the
cohomology of the twisted Aomoto complex (\S\ref{subsec:ss}).
\end{proof}

\begin{corollary}
\label{cor:finite-type}
Suppose that $H_q(X^\nu;\Z)$ is finitely generated as an abelian group for every $q$. Then $X$ is
Alexander formal if and only if
the twisted Aomoto complex $\big(H^\bullet(X;\C_{\rho_\zeta}),\eta\cup-\big)$ is acyclic for every
root of unity $\zeta$, and the Aomoto complex $\big(H^\bullet(X;\F_p),\eta\cup-\big)$ is acyclic for
every prime $p$. Moreover, $H^\bullet(X;\C_{\rho_\zeta})=0$ unless $\zeta$ is an eigenvalue of
$h_\ast$ on $H_\ast(X^\nu;\C)$, so the twisted complexes need only be examined at these eigenvalues.
\end{corollary}

\begin{proof}
By the universal coefficient theorem, $H_q(X^\nu;\Q)$ and $H_q(X^\nu;\F_p)$ are finite-dimensional
for every $q$ and every prime $p$, so the first assertion follows from
Proposition~\ref{prop:finite-type}, applied to each clause of Definition~\ref{def:Aformal}. By
Shapiro's lemma, $H^\bullet(X;\C_{\rho_\zeta})$ is the cohomology of
$\Hom_{\C[t^{\pm1}]}\big(C_\ast(X^\nu;\C),\C_\mu\big)$, where $t$ acts on $\C_\mu=\C$ by $\mu=\zeta$ or
$\mu=\zeta^{-1}$, according to convention. By the universal coefficient theorem over the principal
ideal domain $\C[t^{\pm1}]$, its part of degree $m$ is an extension of
$\Hom\big(H_m(X^\nu;\C),\C_\mu\big)$ by $\Ext^1\big(H_{m-1}(X^\nu;\C),\C_\mu\big)$, and both vanish
unless $\mu$ is an eigenvalue of $h_\ast$, these modules being finite-dimensional. Finally, the roots
of unity that
are eigenvalues of $h_\ast$ form a set stable under $\zeta\mapsto\zeta^{-1}=\bar\zeta$, because
$h_\ast$ is defined over $\Q$.
\end{proof}

\begin{corollary}
\label{cor:collision}
Suppose that $H_\ast(X^\nu;\Z)$ is finitely generated as an abelian group and that, for some prime
power $q=p^a>2$, a
primitive $q$-th root of unity is an eigenvalue of $h_\ast$ on $H_\ast(X^\nu;\Q)$. Then $X$ is not
\etaform{} over $\F_p$, and so is not Alexander formal.
\end{corollary}

\begin{proof}
By Lemma~\ref{lem:collision} the $(t-1)$-primary part of $H_\ast(X^\nu;\F_p)$ is not semisimple;
apply Definition~\ref{def:etak}.
\end{proof}

\begin{remark}
\label{rem:allzeta}
The restriction to prime-power characters in Question~\ref{q:main} comes from the modular 
constraint~\ref{f2}, in which $\beta^p_1(X,\eta)$ concerns only the characters of
order~$p^s$. Over $\Q$, however, 
\emph{every} eigenvalue matters: a fibered $X$ may carry a rational size-$2$ block 
at a non-prime-power root of unity (say, a primitive $6$th root), and such a block is 
named neither by the $\lambda=1$ clause nor by the prime-power
clauses. Indeed a $\LA$-lattice in $\Q[t^{\pm1}]/(\Phi_6^2)$ reduces modulo $2$ to a module supported at
the irreducible $t^2+t+1$, and modulo $3$ to one supported at $t+1$; in neither case, nor
modulo any other prime, since $\Phi_6(1)=1$, does it meet the point $(p,t-1)$ to which
$\beta^p_1(X,\eta)$ refers. Definition~\ref{def:Aformal} demands $\eta$-formality 
at all $\Phi_d$, which closes the gap. The gap is genuine, and a closed $3$-manifold 
already realizes it.
\end{remark}

\begin{example}
\label{ex:phi6}
Let $\varphi$ be an orientation-preserving homeomorphism of the closed surface 
$\Sigma_2$ acting on $H_1(\Sigma_2;\Z)=\Z^4$, in a symplectic basis $a_1,a_2,b_1,b_2$, 
so that the intersection form has matrix $\left(\begin{smallmatrix}0&I_2\\-I_2&0\end{smallmatrix}\right)$, by
\[
M=\begin{pmatrix}1&-1&1&-1\\1&0&1&0\\0&0&0&-1\\0&0&1&1\end{pmatrix};
\]
$M$ is symplectic, and such a $\varphi$ exists because $\Mod(\Sigma_2)\to\Sp_4(\Z)$
is onto. The characteristic polynomial of $M$ is $\Phi_6(t)^2$, and $\Phi_6(M)\ne0=\Phi_6(M)^2$, so
over $\C$ the matrix is $J_2(\zeta_6)\oplus J_2(\bar\zeta_6)$. Let $X$ be the mapping torus and
$\nu$ the class of the fibration. The Alexander module is $\LA/(t-1)$ in degrees $0$ and $2$, and
$H_1(\Sigma_2;\Z)$ with $t$ acting by $M$ in degree one. Three things follow.
\begin{enumerate}[label=\textup{(\roman*)},leftmargin=2.2em,itemsep=2pt]
\item $M$ has no eigenvalue $1$, so $b_1(X)=1$; hence $\pi_1(X)$ is $1$-formal
\cite[\S7.1]{Suciu-Wang-forum19}, and $X$ is even formal, by Lemma~\ref{lem:b1one} below.
\item $M$ has no eigenvalue at a prime-power root of unity, so the rational prime-power clauses hold
vacuously; and since $\Phi_6(1)=1$ is a unit modulo every $p$ and $H_\ast(\Sigma_2;\Z)$ is
torsion-free, the $(t-1)$-primary part of $H_\ast(X^\nu;\F_p)$ lives in degrees $0$ and $2$, where
$t$ acts trivially. So $X$ is $\eta$-formal over every $\F_p$.
\item The block $J_2(\zeta_6)$ means that $X$ is not $\eta$-formal over $\Q(\zeta_6)$
(\S\ref{sec:setup}), so $X$ is not Alexander formal.
\end{enumerate}
Thus $X$ satisfies both conditions of Question~\ref{q:main} --- with $1$-formality at 
$\lambda=1$, not merely $\eta$-formality --- and fails Alexander formality, so the clauses of
Definition~\ref{def:Aformal} at roots of unity whose order is not a prime power cannot be omitted. 
Since a Jordan block of size $>1$ at any eigenvalue rules out strong formality 
\cite[Thm.~3.1]{Pajitnov-ejm17}, $X$ lies, like Example~\ref{ex:paj}, in the ``neither'' 
class of \S\ref{subsec:corners}.
\end{example}

We can now state the answer to the first clause of Question~\ref{q:main}. We read its second
condition, $\eta$-formality at every prime-power character, both in characteristic zero and modulo
$p$; it is condition~\ref{as2} of the following list, in which the clauses of
Definition~\ref{def:Aformal} at $\zeta\ne1$ are sorted according to the order of $\zeta$:
\begin{enumerate}[label=\textup{(\alph*)},leftmargin=2.2em,itemsep=2pt]
\item \label{as1} 
$X$ is \etaform{} over $\Q$.
\item \label{as2} 
$X$ is \etaform{} over $\Q(\zeta)$ at every root of unity $\zeta\ne1$ of prime-power
order, and \etaform{} over $\F_p$ for every prime $p$.
\item \label{as3} 
$X$ is \etaform{} over $\Q(\zeta)$ at every root of unity $\zeta\ne1$ whose order is
not a prime power.
\end{enumerate}
By Definition~\ref{def:Aformal}, $X$ is Alexander formal if and only if \ref{as1}, \ref{as2}, and
\ref{as3} hold.

\begin{corollary}
\label{cor:assembly}
The conjunction of $1$-formality of $\pi_1(X)$ with condition~\ref{as2}, which is what
Question~\ref{q:main} asks about, neither implies nor is implied by Alexander formality.
\end{corollary}

\begin{proof}
The mapping torus of Example~\ref{ex:phi6} has $1$-formal
fundamental group and satisfies \ref{as2}, but not \ref{as3}, so it is not Alexander formal.
Conversely, the space $S^1\times H$ of the proof of Proposition~\ref{thm:assembly}\ref{a1} 
has trivial monodromy on $H_\ast(X^\nu;\Z)$, so the $(t-\zeta)$-primary part of the torsion of
$H_\ast(X^\nu;\Q(\zeta))$ vanishes for $\zeta\ne1$ and $t$ acts trivially on $H_\ast(X^\nu;\k)$ for
every field $\k$; hence it is Alexander formal, while its fundamental group is not $1$-formal.
\end{proof}

\section{\texorpdfstring{$\mathcal{C}$}{C}-graded formality}
\label{sec:covers}

The Kohno--Pajitnov algebra $\bar A^\bullet(X)$ is graded by the group $\CC$ of rank-one characters,
and the CDGA of a finite abelian cover of $X$ is its part over a finite subgroup of $\CC$. This
section develops the consequences of that observation. Section~\ref{subsec:KP} recalls the algebra
and records its relation to finite abelian covers (Lemma~\ref{lem:graded}). Section~\ref{subsec:G}
proves Theorem~\ref{thm:G}, which passes formality of a CDGA graded by an abelian group to its
sub-CDGAs over all subgroups, and \S\ref{subsec:Fformal} places strong formality in Pajitnov's scale of
$\FF$-formality, as its rank-one case. Its consequences for strong formality are drawn in \S\ref{sec:strong}.

\subsection{The Kohno--Pajitnov algebra}
\label{subsec:KP}
We recall the construction of \cite[\S6]{Kohno-Pajitnov-cejm14}.
For a smooth manifold $M$, possibly noncompact, and a character $\rho\in\CC=H^1(M;\C^\ast)$, let
$E_\rho=\widetilde M\times_\rho\C$ be the flat line bundle with monodromy $\rho$, and let
$A^\bullet(M,E_\rho)$ be the complex of smooth $E_\rho$-valued forms, with the differential induced
by the flat connection; its cohomology is $H^\bullet(M;\C_\rho)$. The canonical isomorphisms
$E_\rho\otimes E_\sigma\cong E_{\rho\sigma}$ make 
\[
\bar A^\bullet(M)=\bigoplus_{\rho\in\CC} A^\bullet(M,E_\rho)
\]
a CDGA graded by $\CC$, in the sense of \S\ref{subsec:graded}. A smooth homotopy equivalence 
$f\colon M'\to M$ induces an isomorphism
$f^\ast\colon\CC_M\to\CC_{M'}$ and, through $f^\ast E_\rho\cong E_{f^\ast\rho}$, a quasi-isomorphism
$\bar A^\bullet(M)\to\bar A^\bullet(M')$ of CDGAs covering it. A finite CW-complex $X$ is homotopy
equivalent to a smooth manifold, for instance an open regular neighborhood of a homotopy equivalent
finite simplicial complex in some Euclidean space, and we put $\bar A^\bullet(X)=\bar A^\bullet(M)$
and $A^\bullet(X,E_\rho)=A^\bullet(M,E_\rho)$ for such an $M$. By the above, $\bar A^\bullet(X)$ is
well defined up to quasi-isomorphism of CDGAs graded by $\CC=H^1(X;\C^\ast)$, so formality of
$\bar A^\bullet(X)$, and of each of its sub-CDGAs $\bigoplus_{\rho\in\DD}A^\bullet(X,E_\rho)$ over a
subgroup $\DD\subseteq\CC$, depends only on the homotopy type of $X$. Following Kohno and Pajitnov
\cite[Def.~6.2]{Kohno-Pajitnov-cejm14}, $X$ is \emph{strongly formal} if $\bar A^\bullet(X)$ is formal
as a CDGA.

\begin{remark}
\label{rem:clash}
This is not the notion of the same name of Milivojevi\'c and Stelzig \cite{Milivojevic-Stelzig-cag24},
which we call \emph{bigraded} strong formality: a compact complex manifold has it if its algebra of
forms, with its bigrading and both differentials $\partial,\bar\partial$, is connected to one with
$\partial=\bar\partial=0$ by maps inducing isomorphisms on Bott--Chern and Aeppli cohomology
\cite[Def.~4.1]{Milivojevic-Stelzig-cag24}. Its obstructions are the Aeppli--Bott--Chern--Massey
products, which refine the ordinary ones \cite[\S3]{Milivojevic-Stelzig-cag24},
\cite[\S2.3]{Placini-Stelzig-Zoller-jtop26}. Both notions imply formality (strong formality 
by Theorem~\ref{thm:G} below, applied with $\DD=\{\one\}$), and compact K\"ahler 
manifolds separate them: they are strongly formal in the sense of Kohno and Pajitnov
\cite[Thm.~6.5]{Kohno-Pajitnov-cejm14}, but those admitting a nonconstant holomorphic map to 
a curve of genus at least two, for instance such a curve itself, carry nonzero ABC-Massey products
\cite[Thm.~A and Cor.~6]{Placini-Stelzig-Zoller-jtop26}.
\end{remark}

\begin{lemma}
\label{lem:graded}
Let $Y\to X$ be a connected finite abelian cover with deck group $G$, classified by an epimorphism
$H_1(X;\Z)\surj G$, and let $\widehat G=\Hom(G,\C^\ast)\subseteq\CC$ be the corresponding finite
subgroup of torsion characters. Then $A^\bullet(Y)\cong\bigoplus_{\rho\in\widehat G}A^\bullet(X,E_\rho)$
as CDGAs. Thus $A^\bullet(Y)$ is the sub-CDGA of $\bar A^\bullet(X)$ over the subgroup $\widehat G$,
that is, the identity-degree part of the induced grading by $\CC/\widehat G$, and not an arbitrary
subalgebra; in particular, it is a direct summand of $\bar A^\bullet(X)$ as a cochain complex, and
$H^\bullet(Y;\C)=\bigoplus_{\rho\in\widehat G}H^\bullet(X;\C_\rho)$.
\end{lemma}
\begin{proof}
For a finite abelian cover in characteristic zero, $p_\ast\C_Y\cong\bigoplus_{\rho\in\widehat
G}E_\rho$ as algebra bundles, whence the displayed isomorphism; it is an isomorphism of CDGAs because
on both sides the product is induced by $E_\rho\otimes E_\sigma\cong E_{\rho\sigma}$. Since $\widehat
G$ is a subgroup of $\CC$, the sum over $\widehat G$ is closed under multiplication, and the
differential preserves the $\CC$-grading.
\end{proof}

\subsection{Formality of group-graded CDGAs}
\label{subsec:G}
By Lemma~\ref{lem:graded}, what is needed is that formality of a CDGA graded by an abelian group $\CC$
passes to its $\CC$-graded structure: then the sub-CDGA $\bigoplus_{\rho\in\DD}\bar A_\rho$ is formal for
every subgroup $\DD\subseteq\CC$. This is Theorem~\ref{thm:G} below. For $\CC$ finite, the sub-CDGA is the
invariant subalgebra of a finite group action, and the statement reduces to the classical descent of
formality along finite regular covers \cite[\S2.4]{Suciu-revroum24}; the case needed here, a finite
subgroup of a divisible group, cannot be reduced to it (\S\ref{subsec:graded}).

Fix a field $\k$ of characteristic zero and an abelian group $\CC$, and put
$\Gamma=\Hom(\CC,\k^\ast)$. We assume throughout that $\Gamma$ \emph{separates the points}
of $\CC$: for every $\sigma\ne\one$ there is $\chi\in\Gamma$ with $\chi(\sigma)\ne1$. This
holds whenever $\k$ is algebraically closed, since $\k^\ast$ is then divisible, hence injective
as an abelian group, so that a nontrivial character of $\langle\sigma\rangle$ extends to $\CC$.
In the application $\k=\C$.

Under this hypothesis a $\CC$-grading $V=\bigoplus_\sigma V_\sigma$ is the same thing as an
action of $\Gamma$ in which $\chi$ acts on $V_\sigma$ by the scalar $\chi(\sigma)$, and
$V^\Gamma=V_{\one}$. We need a completed version of this, because spaces of multilinear maps
between graded spaces are graded only in the following sense. Call a $\Gamma$-module $V$, with
$\chi\in\Gamma$ acting by $\theta_\chi$, \emph{$\CC$-diagonal} if it carries $\k$-linear endomorphisms
$\pi_\sigma$, for $\sigma\in\CC$, with $\pi_\sigma\pi_\tau=\delta_{\sigma\tau}\pi_\sigma$ such that for 
each $v\in V$ the family
$(\pi_\sigma v)_\sigma$ is summable with sum $v$, and $\theta_\chi v=\sum_\sigma\chi(\sigma)\pi_\sigma
v$ for all $\chi\in\Gamma$. Write $v_\sigma=\pi_\sigma v$ and $V_{\one}=\pi_{\one}V=V^\Gamma$. In
every instance below $V$ is a space of multilinear maps between $\CC$-graded vector spaces and
$v_\sigma$ is the component raising $\CC$-degree by $\sigma$; summability means that on each
homogeneous argument all but finitely many $v_\sigma$ vanish, so that every identity asserted below
may be verified after evaluation on a single homogeneous argument, where all sums are finite.

\begin{lemma}
\label{lem:charcocycle}
Let $V$ be $\CC$-diagonal and let $c\colon\Gamma\to V$ satisfy
\begin{enumerate}[label=\textup{(\roman*)},leftmargin=2.2em,itemsep=1pt]
\item \label{cc1} $c(\chi\chi')=c(\chi)+\theta_\chi c(\chi')$ for all $\chi,\chi'\in\Gamma$.
\item \label{cc2} $c(\chi)=\sum_{\tau\in\CC}\chi(\tau)\,w_\tau$ for some summable family
$(w_\tau)_{\tau\in\CC}$ in $V$.
\end{enumerate}
Then $c(\chi)=u-\theta_\chi u$ for all $\chi$, where $u\coloneqq w_{\one}$, and $u_{\one}=0$. 
In particular $c$ is a coboundary.
\end{lemma}
\begin{proof}
Insert \ref{cc2} into \ref{cc1}. The left-hand side becomes
$\sum_\tau\chi(\tau)\chi'(\tau)w_\tau$, and the right-hand side becomes
$\sum_\tau\chi(\tau)w_\tau+\sum_{\sigma,\tau}\chi(\sigma)\chi'(\tau)(w_\tau)_\sigma$. Both are
combinations of the functions $(\chi,\chi')\mapsto\chi(\sigma)\chi'(\tau)$, which for distinct pairs
$(\sigma,\tau)\in\CC\times\CC$ are distinct characters of the group
$\Gamma\times\Gamma$ --- this is where separation is used --- and hence are linearly independent over
$\k$ by Artin's independence of characters. Evaluating on a fixed homogeneous argument, so that
all sums are finite, and comparing the coefficients at $(\sigma,\tau)$ gives
\[
[\sigma=\tau]\;w_\tau \;=\; (w_\tau)_\sigma \;+\; [\tau=\one]\;w_\sigma .
\]
Taking $\sigma=\tau=\one$ gives $(w_{\one})_{\one}=0$; taking $\tau=\one\ne\sigma$ gives
$w_\sigma=-(w_{\one})_\sigma$. Hence, with $u=w_{\one}$ and $u_{\one}=0$,
\[
c(\chi)=w_{\one}+\sum_{\tau\ne\one}\chi(\tau)w_\tau
=u-\sum_{\tau\ne\one}\chi(\tau)u_\tau=u-\sum_{\tau}\chi(\tau)u_\tau=u-\theta_\chi u . \qedhere
\]
\end{proof}

\noindent Hypothesis \ref{cc2} cannot be dropped: $\Gamma$ acts trivially on $V_{\one}$, so every
homomorphism $\Gamma\to V_{\one}$ is a $1$-cocycle, while every coboundary has zero component in $V_{\one}$. In
our situation \ref{cc2} holds automatically, because the cocycle is assembled by the grading from one
fixed group element.

\begin{lemma}
\label{lem:descent}
Let $\g$ be a DGLA over $\k$, complete with respect to a decreasing filtration
$\g=F^1\g\supseteq F^2\g\supseteq\cdots$ with $dF^i\subseteq F^i$ and
$[F^i,F^j]\subseteq F^{i+j}$, and let $\Gamma$ act on $\g$ by filtration-preserving
automorphisms of DGLAs in such a way that $\g$ is $\CC$-diagonal, the endomorphisms
$\pi_\sigma$ commuting with $d$, preserving each $F^i$, and satisfying
$\pi_\mu[\phi,\psi]=\sum_{\sigma\tau=\mu}[\pi_\sigma\phi,\pi_\tau\psi]$. Let $x_0,x_1$ be
$\Gamma$-invariant Maurer--Cartan elements of $\g$, and suppose that the filtration 
is split compatibly with $x_0$, that is, $\g=\prod_{k\ge1}\g_{(k)}$ with
$F^n\g=\prod_{k\ge n}\g_{(k)}$, each $\g_{(k)}$ is stable under the
$\pi_\sigma$, and $d_{x_0}=d+[x_0,-]$ maps $\g_{(k)}$ into $\g_{(k+e)}$ for some
$e\ge0$ independent of $k$. If $x_0$ and $x_1$ are gauge equivalent in $\g$, then they are
gauge equivalent by a $\Gamma$-invariant gauge transformation: there is
$\alpha\in(\g^0)^\Gamma=\g^0_{\one}$ with $x_1=\exp(\alpha)\ast x_0$.
\end{lemma}
\begin{proof}
Since $F^1\g=\g$ and the filtration is complete and multiplicative, $\g^0$
is pro-nilpotent and $\GG\coloneqq \exp(\g^0)$ is a pro-unipotent group; $\exp$ and $\log$ are
mutually inverse and $\Gamma$-equivariant, so $\GG^\Gamma=\exp(\g^0_{\one})$. Put
$d_{x_0}=d+[x_0,-]$ and $\z=\ker(d_{x_0})\cap\g^0$. Twisting by $x_0$ and using
\[
\exp(\zeta)\ast0=-\frac{e^{\ad\zeta}-1}{\ad\zeta}\big(d_{x_0}\zeta\big),
\]
one gets $\ZZ\coloneqq \mathrm{Stab}_{\GG}(x_0)=\exp(\z)$. Indeed, as $\zeta\in\g^0=F^1\g^0$, the
operator $\ad\zeta$ raises filtration degree, $\ad\zeta(F^i\g)\subseteq F^{i+1}\g$, so the series
$\frac{e^{\ad\zeta}-1}{\ad\zeta}=\sum_{n\ge0}\frac{(\ad\zeta)^n}{(n+1)!}$ converges in the complete
space $\g$ to the identity plus an operator raising filtration degree, and is therefore invertible.
As $x_0$ is $\Gamma$-invariant,
$d_{x_0}$ commutes with the $\pi_\sigma$, so $\z$ is a $\Gamma$-stable, $\CC$-diagonal,
complete filtered Lie subalgebra of $\g^0$. Write $\ZZ_n=\exp(\z\cap F^n)$, so
that $\ZZ_1=\ZZ$ and $\bigcap_n\ZZ_n=\{1\}$. Since $d_{x_0}$ is continuous and
homogeneous for the splitting, $\z=\prod_k\z_{(k)}$ with
$\z_{(k)}=\z\cap\g_{(k)}$, so that $\z\cap F^n=\prod_{k\ge n}\z_{(k)}$; 
and since $d_{x_0}$ commutes with the $\pi_\sigma$, each $\z_{(k)}$ contains 
the sum of every summable family of homogeneous elements of $\z_{(k)}$, 
all components of $d_{x_0}$ of that sum being zero.

Pick $g\in\GG$ with $g\ast x_0=x_1$. Each $\theta_\chi$ fixes $x_0$ and $x_1$, so
$\theta_\chi(g)\ast x_0=x_1$ as well, and
\[
c_g(\chi)\coloneqq g^{-1}\theta_\chi(g)\in\ZZ,
\qquad
c_g(\chi\chi')=c_g(\chi)\cdot\theta_\chi\big(c_g(\chi')\big).
\]
We record that $\log c_g$ has the shape required by Lemma~\ref{lem:charcocycle}\ref{cc2} modulo any
$F^{n+1}$. Write $g=\exp(\gamma)$ and $\gamma=\sum_\sigma\gamma_\sigma$. Then
$\log c_g(\chi)=\mathrm{BCH}(-\gamma,\theta_\chi\gamma)$ is a convergent sum of iterated brackets in
$-\gamma$ and $\theta_\chi\gamma$; expanding each entry into its components produces terms
\[
\chi(\sigma_{i_1})\cdots\chi(\sigma_{i_k})\,\big[\gamma_{\sigma_1},[\gamma_{\sigma_2},\dots]\big]
=\chi\big(\sigma_{i_1}\cdots\sigma_{i_k}\big)\,\big[\gamma_{\sigma_1},[\gamma_{\sigma_2},\dots]\big],
\]
because $\chi$ is a homomorphism. A bracket of length $k$ lies in $F^k$, so modulo $F^{n+1}$ only the
finitely many BCH terms of length $\le n$ survive, and $\log c_g(\chi)\equiv\sum_\tau\chi(\tau)w_\tau$
for a summable family $(w_\tau)$ in $\g^0$. The $w_\tau$ are brackets of components of $\gamma$,
and need not lie in $\z$. Nothing here is special to $g$: the same holds for
every fixed element of $\GG$, in particular for the successive modifications of $g$ made below.

Now approximate. Suppose $c_g(\chi)\in\ZZ_n$ for all $\chi$, which holds for $n=1$. Then
$\log c_g(\chi)\in\z\cap F^n$; let $\bar c(\chi)\in\z_{(n)}$ be its component in
$\g_{(n)}$, that is, its class modulo $F^{n+1}$. We apply Lemma~\ref{lem:charcocycle} in the
$\CC$-diagonal space $V=\g^0/F^{n+1}\g^0=\prod_{k\le n}\g^0_{(k)}$, which
contains $\z_{(n)}$. The Baker--Campbell--Hausdorff correction in
$\log\big(c_g(\chi)\theta_\chi c_g(\chi')\big)-\log c_g(\chi)-\theta_\chi\log c_g(\chi')$ is a sum of
brackets of elements of $F^n$, so it lies in $F^{2n}\subseteq F^{n+1}$, and $\bar c$ satisfies
\ref{cc1}; by the previous paragraph it satisfies \ref{cc2} in $V$, with the family
$(w_\tau\bmod F^{n+1})$. So $\bar c(\chi)=u-\theta_\chi u$ with $u\in V$ and $u_{\one}=0$. For
$\tau\ne\one$, taking $\tau$-components gives $(1-\chi(\tau))\,u_\tau=\pi_\tau\bar c(\chi)\in\mathfrak
z_{(n)}$ for every $\chi$, and choosing $\chi$ with $\chi(\tau)\ne1$, as separation allows, gives
$u_\tau\in\z_{(n)}$. Hence $u=\sum_\tau u_\tau$ lies in $\z_{(n)}\subseteq\mathfrak
z\cap F^n$, by the closure property recorded above. Put $h=\exp(u)\in\ZZ_n$ and $g'=gh$. Then
$g'\ast x_0=x_1$ and
$c_{g'}(\chi)=h^{-1}c_g(\chi)\theta_\chi(h)$ has leading term $-u+\bar c(\chi)+\theta_\chi u=0$, i.e., 
$c_{g'}(\chi)\in\ZZ_{n+1}$ for all $\chi$.

Iterating gives $h_n\in\ZZ_n$ such that $g_n\coloneqq g\,h_1\cdots h_n$ satisfies $g_n\ast x_0=x_1$ 
and $c_{g_n}(\chi)\in\ZZ_{n+1}$. Since $h_n\in\exp(F^n)$, the sequence $(g_n)$ is Cauchy and
converges to some $g_\infty\in\GG$ with $g_\infty\ast x_0=x_1$ and
$c_{g_\infty}(\chi)\in\bigcap_n\ZZ_n=\{1\}$; that is, $\theta_\chi(g_\infty)=g_\infty$ for all
$\chi$. Thus $g_\infty=\exp(\alpha)$ with $\alpha\in\g^0_{\one}$.
\end{proof}

\begin{theorem}
\label{thm:G}
Let $\CC$ be an abelian group and $\k$ a field of characteristic zero such that
$\Gamma=\Hom(\CC,\k^\ast)$ separates the points of $\CC$ (for instance, $\k$ algebraically
closed, and in particular $\k=\C$). Let $\bar A$ be a
$\CC$-graded CDGA over $\k$. If $\bar A$ is formal, then $\bar A$ is formal as a
$\CC$-graded CDGA. Consequently
$\bar A_{\DD}\coloneqq \bigoplus_{\rho\in\DD}\bar A_\rho$ is formal for every subgroup
$\DD\subseteq\CC$.
\end{theorem}
\begin{proof}
Write $H=H^\bullet(\bar A)=\bigoplus_\sigma H^\bullet(\bar A_\sigma)$, a $\CC$-graded
graded-commutative algebra, and let $\Gamma$ act on $H$ by the graded algebra automorphisms
$\theta_\chi|_{H_\sigma}=\chi(\sigma)$.

\vspace*{2pt}
\emph{Step 1: a $\CC$-graded minimal $C_\infty$-model.} In each $\bar A_\sigma$ choose a
splitting $\bar A^i_\sigma=B^i_\sigma\oplus\mathcal H^i_\sigma\oplus L^i_\sigma$ with
$Z^i_\sigma=B^i_\sigma\oplus\mathcal H^i_\sigma$ and
$d\colon L^i_\sigma\isom B^{i+1}_\sigma$, arranging $1\in\mathcal H^0_{\one}$;
assembling over $\sigma$ gives a contraction of $\bar A$ onto $H$ all of whose structure maps 
preserve the $\CC$-grading. Homotopy transfer along a contraction \cite[\S10.3]{Loday-Vallette-12}
equips $H$ with a strictly unital $C_\infty$-structure $m=(m_2,m_3,\dots)$, with $m_1=0$ and $m_2$ 
the cup product, together with a $C_\infty$-quasi-isomorphism $(H,m)\rightsquigarrow\bar A$. 
The operations $m_n$, and the components of arity $n\ge2$ of the quasi-isomorphism, 
whose linear part is the inclusion $i\colon H\to\bar A$, are signed sums over planar binary 
trees of composites of $i$, the projection $p\colon\bar A\to H$, the homotopy $h\colon\bar A\to\bar A$ 
of the contraction, and the product of $\bar A$ \cite[\S9.4 and Thm.~10.3.3]{Loday-Vallette-12}. 
Since $i$, $p$, and $h$ preserve the $\CC$-grading and the product of $\bar A$ is $\CC$-graded, 
so do $m$ and the quasi-isomorphism; equivalently, both are $\Gamma$-invariant. For an additional 
$\Z$-grading and $A_\infty$-structures this is \cite[Cors.~2.19 and~2.20]{Budur-Rubio-adv19}, 
where it is applied to the weight grading of a mixed Hodge model; the argument does not depend 
on the grading group.

\vspace*{2pt}
\emph{Step 2: formality as a gauge condition.} Let $\g$ be the Lie algebra of coderivations
of the cofree Lie coalgebra $L^c(sH)$ that vanish in arity one, with zero differential, bracket the
commutator, and filtration $F^k$ by arity $\ge k+1$; it is complete, and
$[F^i,F^j]\subseteq F^{i+j}$. Its Maurer--Cartan elements are exactly the $C_\infty$-structures on $H$
with vanishing linear part, and gauge equivalence is $C_\infty$-isomorphism with identity linear part.
Conjugation by $\theta_\chi$ makes $\Gamma$ act on $\g$ by filtration-preserving DGLA
automorphisms, and $\g$ is $\CC$-diagonal, $\pi_\sigma\phi$ being the component of
$\phi$ that raises $\CC$-degree by $\sigma$; both $m$ and $m_2$ are $\Gamma$-invariant
Maurer--Cartan elements.

Since $\bar A$ is formal, $\bar A$ and $H$ are joined by a zigzag of CDGA quasi-isomorphisms, so their
minimal $C_\infty$-models are $C_\infty$-isomorphic \cite{Kadeishvili-soobshch82},
\cite[\S\S10.3, 11.4]{Loday-Vallette-12}: there is a $C_\infty$-isomorphism $f\colon(H,m)\to(H,m_2)$.
Because $m_1=0$ on both sides, its linear part $f_1$ is an automorphism of the graded algebra
$(H,m_2)$; replacing $f$ by $f_1^{-1}\circ f$ we may assume $f_1=\id$. So $m$ and $m_2$ are gauge
equivalent in $\g$. The filtration is split by arity, $\g_{(k)}$ consisting of the
coderivations determined by their component of arity $k+1$, and since the differential of $\g$
is zero, $d_{m_2}=[m_2,-]$ raises arity by exactly one. So Lemma~\ref{lem:descent} applies with
$x_0=m_2$, $x_1=m$, and $e=1$, and makes them gauge equivalent by a $\Gamma$-invariant gauge, that is,
by a $\CC$-graded $C_\infty$-isomorphism $(H,m_2)\cong(H,m)$.

\vspace*{2pt}
\emph{Step 3: conclusion.} Composing that isomorphism with the quasi-isomorphism of Step 1 yields a
$\CC$-graded $C_\infty$-quasi-isomorphism $F\colon(H,m_2)\rightsquigarrow\bar A$, whose linear
part $F_1$ is a quasi-isomorphism; rectifying it gives a zigzag of $\CC$-graded CDGA
quasi-isomorphisms between $H$ and $\bar A$, which is $\CC$-graded formality. (For instance
$H\leftarrow\Omega BH\to\bar A$, with $\Omega B$ the bar--cobar resolution, which is functorial for
$C_\infty$-morphisms and built from their structure maps, and so preserves $\CC$-degrees.) Finally let
$\DD\subseteq\CC$ be a subgroup. Each structure map $F_k$ preserves $\CC$-degrees
and $\DD$ is closed under multiplication, so $F_k$ restricts to
$(H_{\DD})^{\otimes k}\to\bar A_{\DD}$; since $m_2$ restricts as well, we obtain a
$C_\infty$-quasi-isomorphism $(H_{\DD},m_2)\rightsquigarrow\bar A_{\DD}$, its linear part
being the restriction of $F_1$ and hence a quasi-isomorphism onto
$H_{\DD}=H^\bullet(\bar A_{\DD})$. Therefore $\bar A_{\DD}$ is formal.
\end{proof}

\begin{corollary}
\label{cor:massey}
Let $\bar A$ be a CDGA graded by an abelian group $\CC$, over a field $\k$ as in Theorem~\ref{thm:G}.
If $\bar A$ is formal, then for every subgroup $\DD\subseteq\CC$ all Massey products in
$H^\bullet(\bar A_\DD)$, of every length, vanish. In particular, if $X$ is strongly formal, then all
Massey products in the cohomology of every connected finite abelian cover of $X$ vanish.
\end{corollary}

\noindent The first assertion follows from Theorem~\ref{thm:G}, since Massey products vanish in a
formal CDGA, and the second from the first and Lemma~\ref{lem:graded}.

\begin{remark}
\label{rem:Gproof}
No finiteness is used in the proof of Theorem~\ref{thm:G} --- not of $\CC$, not of $\DD$, not of
$\bar A$ --- and $\DD$ is not required to have finite index, so the divisibility obstruction recorded
above is bypassed rather than removed. Characteristic zero enters in the $C_\infty$-transfer, in the
identification of the quasi-isomorphism type of a CDGA with the $C_\infty$-isomorphism type of its
minimal model, and in the exponential and Baker--Campbell--Hausdorff series of
Lemma~\ref{lem:descent}; so the argument says nothing over $\F_p$, and we do not know whether a modular
version holds. Finally, the grading group must be abelian; what this means for nonabelian covers is
discussed in \S\ref{subsec:limits}.
\end{remark}

\subsection{\texorpdfstring{$\FF$}{F}-formality}
\label{subsec:Fformal}
Formality and strong formality are the two ends of the rank-one part of a scale due to Pajitnov
\cite[Def.~3.9]{Pajitnov-topappl19}. Call a family $\FF$ of finite-dimensional complex representations
of $\pi_1(X)$ \emph{tensor-closed} if it contains the trivial representation and is closed under tensor
products. (Pajitnov asks only for the second condition; the first, which makes $A^\bullet(X)$ a summand
of the algebra below and the products with classes in $H^1(X;\C)$ products in that algebra, is the
convention of \cite[Prop.~6.8]{Kohno-Pajitnov-cejm14}.) For such a family, the space $X$ is \emph{$\FF$-formal} 
if the DGA $A^\bullet(X,\FF)=\bigoplus_{\rho\in\FF}A^\bullet(X,E_\rho)$ is formal, $E_\rho$ being the flat bundle
with monodromy $\rho$. For $\FF=\{\one\}$ this is Sullivan formality; for
$\FF$ the set of all rank-one characters it is strong formality; and for $X$ compact K\"ahler and $\rho$
semisimple, the tensor powers $\rho^{\otimes n}$ with $n\ge0$ form a family for which $X$ is $\FF$-formal,
by the work of Corlette and Simpson \cite[Prop.~6.8]{Kohno-Pajitnov-cejm14},
\cite[\S3.3]{Pajitnov-topappl19}. The algebra $A^\bullet(X,\FF)$ is graded by
$\FF$ as an abelian monoid, which can be a group only when every member of $\FF$ 
has dimension one, since an inverse $\sigma$ of $\rho$ forces $\dim\rho\cdot\dim\sigma=1$. 
In higher rank, $A^\bullet(X,\FF)$ is moreover only a DGA, since the symmetry of 
$\rho\otimes\sigma\cong\sigma\otimes\rho$ is not the identity. So Theorem~\ref{thm:G} 
reaches the rank-one part of the scale, and its proof reaches a little beyond.

\begin{proposition}
\label{prop:Fformal}
\begin{enumerate}[label=\textup{(\roman*)},leftmargin=2.2em,itemsep=2pt]
\item \label{ff1} If $X$ is $\FF$-formal for a tensor-closed family $\FF$ of rank-one characters,
then $X$ is $\FF'$-formal for every tensor-closed subfamily $\FF'\subseteq\FF$. In particular, if $X$ is
strongly formal, then $X$ is $\FF$-formal for every such $\FF$.
\item \label{ff2} Let $\rho$ be a representation of $\pi_1(X)$ of dimension at least two. If $X$ is
$\FF$-formal for $\FF=\{\rho^{\otimes n}\}_{n\ge0}$, then $X$ is $\FF_S$-formal for
$\FF_S=\{\rho^{\otimes n}\}_{n\in S}$, for every submonoid $S$ of $\N$.
\end{enumerate}
\end{proposition}

\begin{proof}
\ref{ff1} Such a family is a submonoid $\FF$ of $\CC$, and $A^\bullet(X,\FF)=\bar A_\FF$ is a CDGA
graded by $\CC$. If it is formal, then it is formal as a $\CC$-graded CDGA, by Theorem~\ref{thm:G};
Step~3 of the proof of that theorem uses only that $\DD$ is closed under multiplication, so it applies
with $\DD=\FF'$ and shows that $\bar A_{\FF'}$ is formal. The last assertion is the case $\FF=\CC$. 
\ref{ff2} The powers of $\rho$ have distinct dimensions, so $\FF\cong\N$, and $A^\bullet(X,\FF)$ 
is a DGA graded by $\N\subset\Z$.
The proof of Theorem~\ref{thm:G} goes through for DGAs graded by an abelian group, with
$A_\infty$-structures and coderivations of the tensor coalgebra in place of $C_\infty$-structures and the
cofree Lie coalgebra. So $A^\bullet(X,\FF)$ is formal as a $\Z$-graded DGA, and, by Step~3, its part over
$S$, which is $A^\bullet(X,\FF_S)$, is formal.
\end{proof}

\begin{question}
\label{q:Fformal}
Let $\FF$ be a tensor-closed family of finite-dimensional complex representations of $\pi_1(X)$, that
is, a family that contains the trivial representation and is closed under tensor products.
\begin{enumerate}[label=\textup{(\alph*)},leftmargin=2.2em,itemsep=2pt]
\item \label{qf1} If $X$ is $\FF$-formal, is $X$ also $\FF'$-formal for every tensor-closed subfamily
$\FF'\subseteq\FF$?
\item \label{qf2} Suppose $\FF$ is generated, under tensor products, by the irreducible representations
of $\pi_1(X)$ that factor through a finite quotient $G$, and let $Y\to X$ be the corresponding regular
cover. If $X$ is $\FF$-formal, is $Y$ formal?
\end{enumerate}
\end{question}

\noindent Proposition~\ref{prop:Fformal} answers~\ref{qf1} for families of rank one, and for the
tensor-closed subfamilies of the powers of a single representation. For $G$ abelian, the family
in~\ref{qf2} consists of the characters in $\widehat G$, and the answer is yes, by
Lemma~\ref{lem:graded}. For $G$ nonabelian, a
positive answer would show, with Example~\ref{ex:nonab-cover}, that strong formality does not imply
$\FF$-formality in higher rank. Along a single direction, Remark~\ref{rem:nonab-eta} gives a
necessary condition for the formality of $Y$ that can be checked on $X$.

\section{Strong formality: covers and jump loci}
\label{sec:strong}

Applied to the Kohno--Pajitnov algebra, Theorem~\ref{thm:G} has two kinds of consequences for strong
formality, according to the part of $\CC$ to which the $\CC$-graded formality is restricted: a subgroup
gives the cover principle and the stability properties of strong formality (\S\ref{subsec:cover}), and a
pair of degrees $\{\one,\rho\}$ gives the twisted tangent cone theorem (\S\ref{subsec:jump}).
Section~\ref{subsec:limits} shows that the cover principle neither extends to nonabelian covers nor
admits a converse.

\subsection{The cover principle}
\label{subsec:cover}

\begin{proposition}
\label{prop:cover}
If $X$ is strongly formal, then every connected finite abelian cover of $X$ is formal. 
\end{proposition}

\begin{proof}
By Lemma~\ref{lem:graded}, $A^\bullet(Y)=\bigoplus_{\rho\in\widehat G}A^\bullet(X,E_\rho)$ is
the sub-CDGA of $\bar A^\bullet(X)$ over the subgroup $\widehat G\subseteq\CC$. Strong
formality says $\bar A^\bullet(X)$ is formal, so Theorem~\ref{thm:G}, applied over $\k=\C$ with
$\DD=\widehat G$, gives formality of $A^\bullet(Y)$, and hence of $Y$, since formality 
over $\C$ implies formality over $\Q$ \cite{Sullivan-ihes77}.
\end{proof}

\noindent Theorem~\ref{thm:G} says more than Proposition~\ref{prop:cover} uses, in two ways: the
subgroup $\DD$ is arbitrary, and formality of $\bar A_{\DD}$ kills its Massey products of all lengths
(Corollary~\ref{cor:massey}).

\begin{corollary}
\label{cor:allD}
Let $X$ be strongly formal, and let $\DD$ be any subgroup of the character group
$\CC=H^1(X;\C^\ast)$. Then the CDGA $\bigoplus_{\rho\in\DD}A^\bullet(X,E_\rho)$ is formal, and every
Massey product in its
cohomology $\bigoplus_{\rho\in\DD}H^\bullet(X;\C_\rho)$, of every length, vanishes. For $\DD$ finite
this recovers Proposition~\ref{prop:cover}. For $\DD=\DD_\nu\coloneqq\{\rho:\rho|_{\ker\nu}=\one\}\cong\C^\ast$, 
attached to an epimorphism $\nu\colon\pi_1(X)\surj\Z$, it says that the CDGA assembling the twisted
models along the entire line of characters that factor through $\nu$ is formal.
\end{corollary}

\noindent The same machinery gives the stability properties of strong formality that we shall need, and a
criterion for it.

\begin{lemma}
\label{lem:products}
Let $X$ and $Y$ be connected finite CW-complexes.
\begin{enumerate}[label=\textup{(\roman*)},leftmargin=2.2em,itemsep=2pt]
\item \label{pr1} If $Y$ is a homotopy retract of $X$ and $X$ is strongly formal, then $Y$ is
strongly formal.
\item \label{pr2} $\bar A^\bullet(X\times Y)$ is quasi-isomorphic to $\bar A^\bullet(X)\otimes\bar
A^\bullet(Y)$ as a CDGA graded by $\CC_{X\times Y}=\CC_X\times\CC_Y$; in particular, a product of
strongly formal spaces is strongly formal.
\item \label{pr3} If $X$ and $Y$ are strongly formal, then so is $X\vee Y$.
\item \label{pr4} Suppose that there is $d\ge1$ such that $H^q(Y;\C_\rho)=0$ for all $\rho\in\CC_Y$ when
$q>d$, and for all $\rho\ne\one$ when $q<d$. If $Y$ is formal, then $Y$ is strongly formal. For $d=1$,
for instance if $Y$ is homotopy equivalent to a graph, $Y$ is automatically formal.
\end{enumerate}
\end{lemma}

\begin{proof}
We use the strictly unital minimal $C_\infty$-models of Step~1 of the proof of
Theorem~\ref{thm:G}: a CDGA is formal if and only if its minimal model $(H,m)$ is
$C_\infty$-isomorphic to $(H,m_2)$, and the isomorphism may be taken to have identity linear part.

\ref{pr1} Let $i\colon Y\to X$ and $r\colon X\to Y$ satisfy $r\circ i\simeq\id_Y$. Pulling back local
systems and forms, they induce CDGA maps $i^\ast\colon\bar A^\bullet(X)\to\bar A^\bullet(Y)$ and
$r^\ast\colon\bar A^\bullet(Y)\to\bar A^\bullet(X)$, covering the homomorphisms of character groups
induced by $i$ and $r$. The composite $\theta=H(i^\ast)\circ H(r^\ast)$ is an automorphism of the
$\CC_Y$-graded algebra $H_Y$: a homotopy $r\circ i\simeq\id_Y$ identifies the two maps induced in
twisted cohomology up to transport along the loop $\gamma$ traced by the base point, so that $\theta$
acts on $H^\bullet(Y;\C_\rho)$ as the scalar $\rho(\gamma)$; it is the identity when the homotopy
fixes the base point, but not in general. Transfer gives a $C_\infty$-morphism
$i_\sharp\colon(H_X,m^X)\to(H_Y,m^Y)$ with linear part $H(i^\ast)$. If
$g\colon(H_X,m_2)\to(H_X,m^X)$ is a $C_\infty$-isomorphism with identity linear part, then
$i_\sharp\circ g\circ H(r^\ast)\colon(H_Y,m_2)\to(H_Y,m^Y)$, in which $H(r^\ast)$ is a strict
morphism because it is multiplicative, is a $C_\infty$-morphism with invertible linear part $\theta$,
hence an isomorphism.

\ref{pr2} The characters of $\pi_1(X\times Y)=\pi_1(X)\times\pi_1(Y)$ are the pairs $(\rho,\rho')$,
with $E_{(\rho,\rho')}=E_\rho\boxtimes E_{\rho'}$. The exterior product of forms gives maps
$A^\bullet(X,E_\rho)\otimes A^\bullet(Y,E_{\rho'})\to A^\bullet(X\times Y,E_{(\rho,\rho')})$,
which are quasi-isomorphisms by the K\"unneth theorem for local systems; they are multiplicative and
respect the gradings, so they assemble into a quasi-isomorphism $\bar A^\bullet(X)\otimes\bar
A^\bullet(Y)\to\bar A^\bullet(X\times Y)$ of $\CC_X\times\CC_Y$-graded CDGAs. A tensor product of
formal CDGAs is formal.

\ref{pr3} The characters of $\pi_1(X\vee Y)=\pi_1(X)\ast\pi_1(Y)$ are the pairs $\rho=(\rho_X,\rho_Y)$,
so $\CC=\CC_X\times\CC_Y$. Realize $X\vee Y$, as in \S\ref{subsec:KP}, by a manifold $U\cup V$ with $U$
and $V$ open, $U\simeq X$, $V\simeq Y$, and $U\cap V$ contractible. Over $U\cap V$ the bundles $E_\rho$
are trivialized compatibly with products, so evaluation at a point is a quasi-isomorphism from
$\bigoplus_\rho A^\bullet(U\cap V,E_\rho)$ onto the group algebra $\C[\CC]$, in degree $0$; and
$\bigoplus_\rho A^\bullet(U,E_\rho)=\bar A^\bullet(X)\otimes\C[\CC_Y]$, and similarly for $V$. By
Mayer--Vietoris, $\bar A^\bullet(X\vee Y)$ is quasi-isomorphic, as a $\CC$-graded CDGA, 
to the homotopy pullback of
\[
\begin{tikzcd}[column sep=32pt]
\bar A^\bullet(X)\otimes\C[\CC_Y]  \ar[r, "\varepsilon_X\otimes\id"] & \C[\CC]
& \ar[l, "\ \:\id\otimes\varepsilon_Y"']  \C[\CC_X]\otimes\bar A^\bullet(Y),
\end{tikzcd}
\]
the maps evaluating $0$-forms at the base point. A $\CC_X$-graded minimal Sullivan model of
$\bar A^\bullet(X)$ exists, since $H^0=\C$, and admits exactly one CDGA map to $\C[\CC_X]$, which is
concentrated in degree $0$. As $X$ is strongly formal, Theorem~\ref{thm:G} and the lifting property
therefore allow us to replace $\bar A^\bullet(X)$ by $H_X$ and $\varepsilon_X$ by the augmentation
$\varepsilon\colon H_X\to\C\cdot\one\subseteq\C[\CC_X]$; similarly for $Y$. So
$\bar A^\bullet(X\vee Y)$ is quasi-isomorphic to the homotopy pullback $P$ of
$H_X\otimes\C[\CC_Y]\to\C[\CC]\leftarrow\C[\CC_X]\otimes H_Y$. Computed with the path object
$\C[\CC]\otimes\C[s]\otimes\bwedge(ds)$, the elements of $P$ are triples $(b,\gamma,b')$ with
$\gamma(0)=\varepsilon(b)$ and $\gamma(1)=\varepsilon(b')$. The unit, the triples $(h,0,0)$ and
$(0,0,h')$ with $h,h'$ homogeneous of positive degree, and the triples $e_\rho=(0,ds\cdot[\rho],0)$ with
$\rho_X\ne\one\ne\rho_Y$, are cocycles spanning a sub-CDGA $S$ with zero differential, since every
product involving some $e_\rho$ vanishes. To see that $S\to P$ is a quasi-isomorphism, 
fix $\rho=(\rho_X,\rho_Y)$. In degree $\rho$, evaluation at $s=0$ and $s=1$ maps 
$\C[s]\otimes\bwedge(ds)$ onto $\C\oplus\C$, with kernel $K$ satisfying $H^0(K)=0$ and 
$H^1(K)=\C\cdot[ds]$; so projection onto the outer entries gives a short exact sequence
\[
\begin{tikzcd}[column sep=18pt]
0\ar[r]& K \ar[r]& P_\rho\ar[r]& H^\bullet(X;\C_{\rho_X})\oplus H^\bullet(Y;\C_{\rho_Y}) \ar[r]& 0 .
\end{tikzcd}
\]
Its connecting homomorphism vanishes except in degree $0$, where it sends $(b,b')$ to
$(\varepsilon(b')-\varepsilon(b))\,[ds]$; as $H^0(X;\C_{\rho_X})$ is $\C$ if $\rho_X=\one$ and zero
otherwise, and likewise for $Y$, it is onto $H^1(K)$ unless $\rho_X\ne\one\ne\rho_Y$, and zero in that
case. Hence $H^0(P_\rho)$ is $\C\cdot\one$ for $\rho=\one$ and zero otherwise, and the projection
$H^q(P_\rho)\to H^q(X;\C_{\rho_X})\oplus H^q(Y;\C_{\rho_Y})$ is an isomorphism for $q\ge2$, and for
$q=1$ unless $\rho_X\ne\one\ne\rho_Y$, when its kernel is spanned by the class of $e_\rho$. These are
exactly the classes of $S_\rho$, so $S\to P$ is a quasi-isomorphism; this is the Mayer--Vietoris
count for $H^\bullet(X\vee Y;\C_\rho)$, the classes $e_\rho$ accounting for those that come from
neither summand. So $P$, and hence $\bar A^\bullet(X\vee Y)$, is formal.

\ref{pr4} Since $H^0(Y;\C_\rho)=0$ for $\rho\ne\one$, the cohomology
$H_Y=\bigoplus_\rho H^\bullet(Y;\C_\rho)$ is $H_{\one}\oplus N$, where $H_{\one}=H^\bullet(Y;\C)$ lies
in degrees $0,\dots,d$, with $H_{\one}^0=\C\cdot1$, and $N=\bigoplus_{\rho\ne\one}H^d(Y;\C_\rho)$ lies
in degree $d$. Let $(H_Y,m)$ be the strictly unital minimal model, so that $m_k$ has degree $2-k$. The
argument has three steps.

\emph{Operations involving $N$ are trivial.} For $k\ge2$ the operation $m_k$ vanishes on arguments one 
of which lies in $N$, except for $m_2(1,x)=x=m_2(x,1)$. Indeed, by strict unitality we may assume that no
argument is $1$, so that all have positive degree; if $j\ge1$ of the $k$ arguments lie in $N$, the
value has degree at least $jd+(k-j)+2-k=j(d-1)+2>d$, while $H_Y$ vanishes above degree $d$. The same
holds for the cup product $m_2$ alone.

\emph{A strictly unital isomorphism on $H_{\one}$.} Restricting the contraction of Step~1 to the
identity degree shows that $(H_{\one},m)$ is a minimal model of $A^\bullet(Y)$. As $Y$ is formal, there
is a $C_\infty$-isomorphism $g\colon(H_{\one},m_2)\to(H_{\one},m)$ with identity linear part. We may
take $g$ strictly unital, that is, with $g_k$ vanishing on arguments involving $1$ for $k\ge2$. To see
this, note that $g_k$, of degree $1-k$, takes arguments of positive degree to values of positive degree,
and so do $m_k$ and $m_2$; hence the equations for a $C_\infty$-morphism, evaluated on arguments of
positive degree, involve only components of $g$ on such arguments, and they still hold after the
components of $g$ on arguments involving $1$ are replaced by~$0$. The equations evaluated on arguments
involving $1$ then hold because both structures are strictly unital: the only terms that can be nonzero
are those in which a $1$ is multiplied by an adjacent argument under $m_2$, and for a $1$ inside the
string of arguments the two such terms on the source side cancel, there being none on the target side,
while for a $1$ at either end there is one such term on each side, and they agree.

\emph{Extension by the identity on $N$.} Extend $g$ to $H_Y$ by $g_1=\id$ on $N$, and by $g_k=0$,
for $k\ge2$, on arguments involving $N$. Evaluated on arguments involving $N$, the equations for a
$C_\infty$-morphism $(H_Y,m_2)\to(H_Y,m)$ then reduce, by the first step, to the terms in which an
element of $N$ is multiplied by $1$: the components $g_k$ with $k\ge2$ are strictly unital, so their
values have positive degree. These terms occur only for two arguments, and are the same on both sides.
So the extension is a $\CC_Y$-graded $C_\infty$-isomorphism $(H_Y,m_2)\to(H_Y,m)$, and
$\bar A^\bullet(Y)$ is formal.

Finally, for $d=1$ the count in the first step shows that $m_k$ takes $(H_{\one}^1)^{\otimes k}$ to
degree $2$, where $H_Y$ is zero; so $(H_{\one},m)$ is $(H_{\one},m_2)$, and $Y$ is formal.
\end{proof}

\begin{corollary}
\label{cor:milnor-strong}
Let $\A$ be a central arrangement such that some connected finite abelian cover of $U(\A)$ or
of $M(\A)$, for instance the Milnor fiber $F(\A)\to U(\A)$, is not formal. Then the complement $M(\A)$
is formal but not strongly formal.
\end{corollary}

\begin{proof}
Brieskorn's theorem \cite{Brieskorn-71} gives formality of $M(\A)$. By Proposition~\ref{prop:cover},
the space covered, $U(\A)$ or $M(\A)$, is not strongly formal; in the first case neither is
$M(\A)\cong U(\A)\times\C^\ast$, which retracts onto $U(\A)$, by Lemma~\ref{lem:products}\ref{pr1}.
\end{proof}

\noindent The hypothesis holds for the Ceva arrangement $\A(3,3,3)$, whose Milnor fiber is
not even $1$-formal \cite{Zuber-blms10}, and for the arrangements $\A(3k,3k,3)$ 
\cite[Thm.~1.2]{Suciu-arx26a}. When $b_1(X)=0$, the cover principle has an 
elementary proof, independent of Theorem~\ref{thm:G}.

\begin{corollary}
\label{cor:b1zero}
Suppose $b_1(X)=0$, so that $\CC$ is finite and equal to $\Hom(H_1(X;\Z),\C^\ast)$. Then
$\bar A^\bullet(X)\cong A^\bullet(X^{\ab})$, where $X^{\ab}\to X$ is the maximal abelian cover, 
finite with deck group $H_1(X;\Z)$; so $X$ is strongly formal if and only if $X^{\ab}$ is formal.
\end{corollary}

\begin{proof}
Apply Lemma~\ref{lem:graded} with $G=H_1(X;\Z)$, for which $\widehat G=\CC$.
\end{proof}

\noindent Every connected finite abelian cover $Y$ of $X$ is covered by $X^{\ab}$, regularly, and
formality descends along finite regular covers \cite[\S2.4]{Suciu-revroum24}; so
if $X$ is strongly formal, $Y$ is formal. Corollary~\ref{cor:b1zero} also lets one certify 
by hand that a space is not strongly formal, on the standard example of a formal space 
with a non-formal finite cover.

\begin{example}
\label{ex:dp}
Let $K=F_2/\gamma_3(F_2)$ be the free group of rank two modulo triple commutators, let $\varphi$ be
the order-two automorphism $x\mapsto x^{-1}$, $y\mapsto y^{-1}$, and let $X$ be a finite $2$-complex
with $\pi_1(X)=G=K\rtimes_\varphi\Z/2$; this is \cite[Ex.~2.9]{Dimca-Papadima-pisa11}. Since
$\varphi$ acts as $-1$ on $K^{\ab}=\Z^2$ we get
$H_1(X;\Z)=\Z^2/2\Z^2\oplus\Z/2=(\Z/2)^3$, so $b_1(X)=0$. Hence $H^\bullet(X;\Q)$ is concentrated in
degrees $0$ and $2$, and $X$ is formal: the operations $m_k$, of degree $2-k$, of a strictly unital
minimal model take $k\ge2$ classes of degree two to degree $k+2$, where the cohomology vanishes.
Let $Y\to X$ be the double cover with $\pi_1(Y)=K$; the deck group acts as $-1$ on $H^1(Y;\C)$, so
\cite[Cor.~2.5]{Dimca-Papadima-pisa11} does not apply, and
indeed $K$ is not $1$-formal (Morgan; see \cite[Ex.~2.9]{Dimca-Papadima-pisa11}), so $Y$ is not
formal. The maximal abelian cover $X^{\ab}$, of degree $8$, covers $Y$ regularly, so it is not formal
either, and $X$ is \emph{not} strongly formal, by Corollary~\ref{cor:b1zero}.
\end{example}

\noindent So Proposition~\ref{prop:cover} needs strong formality, not merely formality;
Corollary~\ref{cor:milnor-strong} and Examples~\ref{ex:link623} and~\ref{ex:paj} below also separate
the two.

\subsection{The twisted tangent cone theorem}
\label{subsec:jump}
Restricting the $\CC$-graded formality of Theorem~\ref{thm:G} to the pair of $\CC$-degrees
$\{\one,\rho\}$, rather than to a subgroup, gives information about the cohomology jump loci near an
arbitrary character $\rho$, including characters of infinite order, which no finite cover reaches. Let
$\VV^i_k(X)$ and $\RR^i_k(X,\rho)$ be the characteristic varieties of $X$ and the resonance varieties of
the twisted Aomoto complexes, defined in \S\ref{subsec:graded}; for $\rho=\one$ the latter are the usual
resonance varieties $\RR^i_k(X)$. We identify the tangent space to $\CC$ at
$\rho$ with $H^1(X;\C)$ through $\omega\mapsto\rho\exp(\omega)$, and write $\mathrm{TC}_\rho$ for the
tangent cone at $\rho$. For any finite complex, Libgober proved the inclusion
$\mathrm{TC}_\rho\big(\VV^i_k(X)\big)\subseteq\RR^i_k(X,\rho)$
\cite[Thm.~3.1]{Libgober-topappl02}. At $\rho=\one$ it is an equality, and the two germs agree, when
$X$ is formal \cite[Thm.~A]{Dimca-Papadima-Suciu-duke09}, \cite[Thm.~B and \S1.4]{Dimca-Papadima-ccm14};
the irreducible components of $\VV^i_k(X)$ through $\one$ are then subtori
\cite[Thm.~C(2)]{Dimca-Papadima-ccm14}, \cite[Thm.~1.3(1)]{Budur-Wang-pamq20}. 
Strong formality gives the same at every character.

\begin{theorem}
\label{thm:jump}
Let $X$ be a strongly formal space with the homotopy type of a finite CW-complex. Then, for every
character $\rho\in\CC$ and all $i\ge0$ and $k\ge1$, the following hold.
\begin{enumerate}[label=\textup{(\roman*)},leftmargin=2.2em,itemsep=2pt]
\item \label{j1}
The map $\omega\mapsto\rho\exp(\omega)$ restricts to an isomorphism from the germ of
$\RR^i_k(X,\rho)$ at $0$ onto the germ of $\VV^i_k(X)$ at $\rho$; in particular,
$\mathrm{TC}_\rho\big(\VV^i_k(X)\big)=\RR^i_k(X,\rho)$.
\item \label{j2}
Every irreducible component of $\RR^i_k(X,\rho)$ is a linear subspace of $H^1(X;\C)$ defined
over $\Q$.
\item \label{j3}
Every irreducible component of $\VV^i_k(X)$ is a translated subtorus of $\CC$.
\end{enumerate}
\end{theorem}

\begin{proof}
By Theorem~\ref{thm:G}, the CDGA $\bar A^\bullet(X)$ is joined to its cohomology
$\bigoplus_{\sigma\in\CC}H^\bullet(X;\C_\sigma)$, with zero differential, by a zigzag of
$\CC$-graded quasi-isomorphisms. In a $\CC$-graded CDGA $B$, the summand $B_{\one}$ is a sub-CDGA, 
and $B_\rho$ is a differential graded module over it, as $B_{\one}B_\rho\subseteq B_\rho$ and the differential
preserves the grading. A $\CC$-graded map $f\colon B\to B'$ of CDGAs restricts to a map
$f_{\one}\colon B_{\one}\to B'_{\one}$ of CDGAs and a map $f_\rho\colon B_\rho\to B'_\rho$ with
$f_\rho(ab)=f_{\one}(a)f_\rho(b)$ for $a\in B_{\one}$ and $b\in B_\rho$, that is, to a map of pairs,
consisting of a map of algebras and a map of modules over it. Since $f$ preserves the grading, $H(f)$
is the direct sum of the maps $H(f_\sigma)$; so if $f$ is a quasi-isomorphism, then so are $f_{\one}$
and $f_\rho$. Hence the pair $\big(A^\bullet(X),A^\bullet(X,E_\rho)\big)$ is joined to
$\big(H^\bullet(X;\C),H^\bullet(X;\C_\rho)\big)$, with zero differentials and the module structure
given by cup product, by a zigzag of maps of pairs whose two components are quasi-isomorphisms. 
We now check the hypotheses of Budur and Wang. Regard $A^\bullet(X)$ as a DGLA with zero bracket and
$A^\bullet(X,E_\rho)$ as a module over it in the sense of \cite[Def.~3.2]{Budur-Wang-comp15a}, the
identity required there reducing to graded commutativity; the maps of pairs above are then morphisms
of DGLA pairs, and those of the zigzag are $\infty$-equivalences \cite[Def.~3.4]{Budur-Wang-comp15a}.
This module is bounded and has finite-dimensional
cohomology, as their standing assumption requires \cite[\S3]{Budur-Wang-comp15a}, and the zigzag
above shows that the pair $\big(A^\bullet(X),A^\bullet(X,E_\rho)\big)$ is formal
\cite[Def.~3.6]{Budur-Wang-comp15a} (the modules of the intermediate pairs of the zigzag need not be
bounded, but their cohomology, isomorphic to $H^\bullet(X;\C_\rho)$, is bounded and
finite-dimensional, which is all that the arguments of \cite{Budur-Wang-comp15a} use). 
The manifold of \S\ref{subsec:KP} on which $A^\bullet(X)$ is
defined has the homotopy type of a finite CW-complex, the local system $E_\rho$ is irreducible,
having rank one, and $\mathcal{E}nd(E_\rho)$ is the trivial line bundle; so
\cite[Thm.~7.4]{Budur-Wang-comp15a}, applied with the trivial local system as coefficients,
identifies the formal germ of $\VV^i_k(X)$ at $\rho$ with the cohomology jump functor of this pair.
That functor is unchanged along the zigzag, whose maps are quasi-isomorphisms on both components
\cite[Thm.~1.2]{Budur-Wang-comp15a}. As the pair is formal and its bracket is zero, so that
$[H^0,H^1]=0$, the functor is prorepresented, for every $i$, by the formal germ at $0$ of the
resonance variety of the cohomology pair \cite[Cor.~4.7]{Budur-Wang-comp15a}. That variety consists
of the classes $\omega$ in the quadratic cone $\{\omega\in H^1(X;\C):[\omega,\omega]=0\}$, which is all
of $H^1(X;\C)$, such that the complex $\big(H^\bullet(X;\C_\rho),\omega\cup-\big)$ has $i$-th
cohomology of dimension at least $k$ \cite[Lem.~4.5]{Budur-Wang-comp15a}; that is, it is
$\RR^i_k(X,\rho)$. Here the formal germ of $\CC$ at $\rho$ is identified with that of $H^1(X;\C)$ at
$0$ through $\omega\mapsto\rho\exp(\omega)$. This identification is induced by an analytic isomorphism
of the ambient germs, and an ideal of convergent power series is determined by its completion, the
completion being faithfully flat; so it is an isomorphism of analytic germs. Since
$\RR^i_k(X,\rho)$ is a homogeneous cone, this proves \ref{j1}.

Part~\ref{j3} follows from \cite[Thm.~1.3(2)]{Budur-Wang-pamq20}: for $X$ with the homotopy type of a
finite CW-complex, that theorem shows that every irreducible component of $\VV^i_k(X)$ through $\rho$ is
a translated subtorus as soon as the pair $\big(A^\bullet(X),A^\bullet(X,E_\rho)\big)$ is homotopy
equivalent to a pair whose module is finite-dimensional in each degree, as the cohomology pair above
is, for every $\rho$. Its proof rests on the exponential Ax--Lindemann theorem in the form
\cite[Prop.~1.4 and Cor.~2.2]{Budur-Wang-pamq20}. By \ref{j1} and \ref{j3}, near $0$ the cone
$\RR^i_k(X,\rho)$ is a finite union of the Lie algebras of subtori of $\CC$, which are defined
over $\Q$; being a cone, it is that union, which proves \ref{j2}.
\end{proof}

\noindent In the language of the $L_\infty$ pairs of Budur and Rubi\'o \cite{Budur-Rubio-adv19}, the pair
$\big(H^\bullet(X;\C),H^\bullet(X;\C_\rho)\big)$ carries a canonical $L_\infty$-module structure, with $m_2$
the cup product, weakly equivalent to $\big(A^\bullet(X),A^\bullet(X,E_\rho)\big)$
\cite[Thm.~2.25 and \S4.1]{Budur-Rubio-adv19}, and along a class $\omega\in H^1(X;\C)$ the operations
$m_{n+1}(\omega,\dots,\omega,-)$ determine the twisted sequence $\E_r^\ast(\rho,\omega)$ from $E_1$ on,
through the twisted differential of \cite[Def.~3.6]{Budur-Rubio-adv19}. The components of the jump loci
through $\rho$ are translated subtori as soon as $m_n(\omega,\dots,\omega,-)=0$ for all large $n$, uniformly
in $\omega$ \cite[Thm.~1.5(2)]{Budur-Rubio-adv19}: bounded length along all classes, rather than vanishing
along a single class as in Proposition~\ref{prop:circle}. The proof above shows that, for strongly formal
$X$, the canonical structure may be taken with $m_n=0$ for all $n\ge3$, by the uniqueness in
\cite[Thm.~1.6]{Budur-Rubio-adv19}.

\begin{example}
\label{ex:link623}
Let $X$ be the complement of the $2$-component link $6^2_3$ in Rolfsen's table. Its linking number is
$2$, and its Alexander polynomial is $\Delta_L=t_1t_2-2(t_1+t_2)+1$. Since the linking number is
nonzero, $X$ is formal \cite[Thm.~10.3]{Suciu-mm22}. On the other hand, the curve $\{\Delta_L=0\}$ is
an irreducible component of $\VV^1_1(X)\subseteq(\C^\ast)^2$ which does not pass through
$\one$ and is not a translated subtorus \cite[Ex.~10.5]{Suciu-mm22}. By
Theorem~\ref{thm:jump}\ref{j3}, $X$ is not strongly formal.
\end{example}

\begin{remark}
\label{rem:jump}
Theorem~\ref{thm:jump} has three precursors. At $\rho=\one$ it is the tangent cone theorem for formal
spaces. For compact K\"ahler manifolds, Budur and Wang proved the analogue of part~\ref{j1} at an
irreducible local system of any rank, through the formality of the corresponding DGLA pairs
\cite[Thm.~7.5 and Cor.~7.6]{Budur-Wang-comp15a}. And Kohno and
Pajitnov showed, under the degeneration of the twisted spectral sequences that strong formality
guarantees, that the directions along which the twisted Betti numbers jump form a finite union of
rationally defined subspaces \cite[Thm.~6.3 and Prop.~10.2]{Kohno-Pajitnov-cejm14}. The theorem
upgrades directions to germs; this needs a finite-dimensional model of each pair
$\big(A^\bullet(X),A^\bullet(X,E_\rho)\big)$ \cite[Thm.~1.3(2)]{Budur-Wang-pamq20}, and
Theorem~\ref{thm:G} provides one at every character at once, namely the cohomology pair itself.

For smooth quasi-projective varieties part~\ref{j3} holds unconditionally \cite{Budur-Wang-ens15}, 
but part~\ref{j1} does not, contrary to \cite[p.~161]{Libgober-topappl02}: it fails already at $\rho=\one$ 
for the configuration spaces of $n\ge3$ points on an elliptic curve, whose first resonance variety is 
irreducible and not linear \cite[Ex.~10.1 and Rem.~10.3]{Dimca-Papadima-Suciu-duke09}, 
\cite[Prop.~22.4]{Suciu-arx26b}. A character at which the inclusion is strict would thus show 
that an arrangement complement is not strongly formal.
\end{remark}

\subsection{The limits of the cover principle}
\label{subsec:limits}
The cover principle raises two questions: does it extend to nonabelian covers, and does it have a
converse? Both answers are negative (Example~\ref{ex:nonab-cover} and Proposition~\ref{prop:converse}).
For nonabelian covers, a necessary condition can nevertheless be checked on $X$ itself, along a single
direction (Remark~\ref{rem:nonab-eta}).

For a finite Galois cover $Y\to X$ with nonabelian deck group $G$, the proof of
Proposition~\ref{prop:cover} breaks down: $p_\ast\C_Y$ contains local systems of rank greater than one,
so $A^\bullet(Y)$ is not a sub-CDGA of $\bar A^\bullet(X)$; a finite cover is dominated by a finite
abelian cover only if it is itself abelian; and Theorem~\ref{thm:G}, applied with $\DD=\widehat G$,
yields only the formality of the maximal abelian intermediate cover, which is $X$ itself when $G$ is
perfect. Compact K\"ahler manifolds give no counterexample: their finite covers are again compact
K\"ahler, hence formal \cite{DGMS-inv75}. The natural place to look is the opposite extreme, $H_1(X;\Z)=0$: 
then $\CC=\{\one\}$, strong formality is plain formality, and every nontrivial finite Galois cover of $X$ has
perfect deck group.

\begin{example}
\label{ex:nonab-cover}
Let $M=\Sigma(2,3,7)$ be the Brieskorn sphere, let $\Gamma$ be the kernel of the Hurwitz
epimorphism $\Delta(2,3,7)\surj\PSL(2,7)$ from the triangle group, and let $N\to M$ be the
degree-$168$ cover obtained by pulling back $\Gamma$ along the Seifert fibration of $M$. We check that
$M$ is strongly formal, that the cover is Galois with simple deck group $\PSL(2,7)$, and that $N$ is
not $1$-formal, hence not formal. Consequently the hypothesis that the cover be abelian cannot be
removed from Proposition~\ref{prop:cover}, and strong formality is not inherited by finite covers.

\vspace*{1pt}
\emph{The base.} $M$ is an integral homology sphere, so $\CC=\Hom(H_1(M;\Z),\C^\ast)$ is
trivial and $\bar A^\bullet(M)=A^\bullet(M)$: strong formality of $M$ is plain formality. And
$H^\bullet(M;\Q)\cong H^\bullet(S^3;\Q)$, which is intrinsically formal, so $M$ is formal. Since
$\pi_1(M)$ is perfect, so is every quotient of it; $\PSL(2,7)$ is one, and it is nonabelian.

\vspace*{1pt}
\emph{The cover.} Write $1\to\Z\langle h\rangle\to\pi_1(M)\to\Delta(2,3,7)\to1$ for the central
extension given by the Seifert fibration, whose Euler number is $e(M)=\pm1/42$. The subgroup $\Gamma$
is torsion-free of index $168$, with $\mathbb{H}^2/\Gamma$ the Klein quartic $\Sigma_3$
\textup{(}the Hurwitz bound $84(g-1)=168$ is attained at $g=3$\textup{)}. Since $\Gamma$ is
torsion-free, the Seifert fibration of $N$ induced from that of $M$ has no exceptional fibers, so $N$ is
the circle bundle over $\Sigma_3$ with Euler number
$e(N)=168\cdot e(M)=\pm4$, and $\pi_1(N)=\langle a_1,b_1,\dots,a_3,b_3,h\mid
[a_i,h]=[b_i,h]=1,\ \prod_i[a_i,b_i]=h^{\pm4}\rangle$. In particular $H_1(N;\Z)=\Z^6\oplus\Z/4$.

\vspace*{1pt}
\emph{Non-$1$-formality of the cover.} The space $N$ is an orientable Seifert fibered space whose base
has genus $3>0$ and whose Euler number is nonzero, so $\pi_1(N)$ is not graded-formal, by
\cite[Prop.~10.7]{Suciu-Wang-forum19}. Since $1$-formality is the conjunction of graded-formality and
filtered-formality \cite[Prop.~7.6]{Suciu-Wang-forum19}, the group $\pi_1(N)$ is not $1$-formal, and
so $N$ is not formal.
\end{example}

We turn to the converse of Proposition~\ref{prop:cover}. Its failure is witnessed by a closed
$3$-manifold fibered over the circle, and the formality of its covers rests on the following standard
fact.

\begin{lemma}
\label{lem:b1one}
A closed orientable $3$-manifold $N$ with $b_1(N)=1$ is formal.
\end{lemma}

\noindent This is a special case of a theorem of Fern\'andez and Mu\~noz \cite{Fernandez-Munoz-proc05};
see also \cite[Thm.~4.2 and Prop.~10.2]{Papadima-Suciu-bullroum09}, where it is applied to mapping 
tori of surfaces whose monodromy does not have $1$ as an eigenvalue.

\begin{proposition}
\label{prop:converse}
The converse of Proposition~\ref{prop:cover} holds when $b_1(X)=0$ and fails in general: there is a
closed orientable $3$-manifold $X$, fibered over the circle with fiber of genus two, such that $X$ and
all of its finite abelian covers are formal, but $X$ is not strongly formal. This $X$ is moreover
Alexander formal.
\end{proposition}
\begin{proof}
If $b_1(X)=0$, the maximal abelian cover $X^{\ab}$ is finite, and it is formal if and only if $X$ is
strongly formal (Corollary~\ref{cor:b1zero}). For the general case, let $\varphi$ be an
orientation-preserving homeomorphism of $\Sigma_2$ acting on $H_1(\Sigma_2;\Z)$, 
in the symplectic basis of Example~\ref{ex:phi6}, by the symplectic matrix
\[
M=\begin{pmatrix}2&1&-2&0\\1&1&-1&0\\0&0&1&-1\\0&0&-1&2\end{pmatrix},
\]
and let $X$ be its mapping torus, with $\nu$ the class of the fibration. The characteristic polynomial of $M$ is $f(t)^2$ with
$f(t)=t^2-3t+1$, and $f(M)\ne0=f(M)^2$, so over $\C$ the matrix is $J_2(\lambda)\oplus J_2(\lambda^{-1})$
with $\lambda=(3+\sqrt5)/2$, which is not a root of unity. Since $\det(M-I)=f(1)^2=1$, we have
$H_1(X;\Z)=\Z$, so the finite abelian covers of $X$ are exactly the cyclic covers $X^{(k)}$ along the
base, the mapping tori of $\varphi^k$. As $\lambda$ is not a root of unity, $M^k$ has no eigenvalue
$1$, so $b_1(X^{(k)})=1$ and $X^{(k)}$ is formal by Lemma~\ref{lem:b1one}. On the other hand,
$\VV^1(X)$ is finite, while $H^1(X;\C_{\rho_\lambda})\cong\C$ is spanned by a class $a$ with
$\eta\cup a=0$: in the Wang sequence, $\eta\cup-$ is the natural map from the invariants to the
coinvariants of the monodromy on $H^1(\Sigma_2;\C)\otimes\C_{\rho_\lambda}$, which has a Jordan block of
size two at~$1$. So $\RR^1_1(X,\rho_\lambda)=H^1(X;\C)\ne\mathrm{TC}_{\rho_\lambda}(\VV^1_1(X))$, 
and $X$ is not strongly formal, by Theorem~\ref{thm:jump}; numerically, $\dim H^1(X;\C_{\rho_\lambda})=1$ 
is less than the multiplicity of $\lambda$ in $f(t)^2$ \cite[Cor.~6.15]{Dimca-Papadima-Suciu-imrn08}, and
$\langle\eta,\eta,a\rangle$ does not vanish \cite[Thm.~5.1]{Pajitnov-ejm17}. Finally, $M$ has no eigenvalue 
at a root of unity, $f(1)^2=1$ is a unit modulo every $p$, and $H_\ast(\Sigma_2;\Z)$ is torsion-free, 
so, exactly as in Example~\ref{ex:phi6}, every clause of Definition~\ref{def:Aformal} holds.
\end{proof}

\noindent Whether the converse of Proposition~\ref{prop:cover} holds for smooth
quasi-projective $X$, where the cohomology jump loci in all degrees are finite unions of
torsion-translated subtori \cite{Budur-Wang-ens15}, we do not know.

\begin{remark}
\label{rem:cover-open}
Besides Question~\ref{q:Fformal}, the solvable range remains open: a regular cover $Y\to X$ with
solvable deck group is a tower of finite abelian covers, so the cover principle would extend to it if
strong formality were inherited by finite abelian covers. We do not expect this, since the restriction
map $\CC_X\to\CC_Y$ is typically far from onto.
\end{remark}

\begin{remark}
\label{rem:nonab-eta}
Along a single direction, the formality of a nonabelian cover can be tested on $X$ itself. Let
$Y\to X$ be a connected finite regular cover with deck group $G$, let $\nu\colon\pi_1(X)\surj\Z$, and let
$\nu_Y$ be the restriction of $\nu$ to $\pi_1(Y)$, divided by its index $m$. Since
$A^\bullet(Y)\cong\bigoplus_\sigma A^\bullet(X,E_\sigma)^{\oplus\dim\sigma}$ as DG modules over
$A^\bullet(X)$, the sum running over the irreducible representations $\sigma$ of $G$, and since $\eta$
pulls back to $m$ times the class of $\nu_Y$, the pair $(Y,\nu_Y)$ is \etaform{} over $\C$ if and only if,
for every $\sigma$, the monodromy has no Jordan block of size greater than one for the eigenvalue $1$ on
the torsion of the twisted Alexander module $H_\ast(X^\nu;\sigma)$ \cite[Thm.~C]{Pajitnov-topappl19};
this holds if $Y$ is formal \cite[Thm.~6.1]{Kohno-Pajitnov-cejm14}. For $G$ abelian, the condition
follows from strong formality of $X$, through Proposition~\ref{prop:cover}; for $G$ nonabelian we do not
know whether it does, as $b_1(X)=0$ in Example~\ref{ex:nonab-cover}. It holds for every $G$ when $X$ has a
circle action whose orbit class $\gamma$ satisfies $\nu(\gamma)\ne0$, as an arrangement complement does
with the class of its Milnor fibration: the action lifts to $Y$ through a finite cover of the circle, so
the monodromy on the homology of the infinite cyclic cover of $Y$ has finite order, hence is semisimple
over $\C$.
\end{remark}

\section{Incomparability}
\label{sec:second}

This section shows that neither condition implies the other (Theorem~\ref{thm:incomp}). An Alexander
formal space can fail to be strongly formal through Massey products not of the form
$\langle\eta,\dots,\eta,a\rangle$ (\S\ref{subsec:dir}); a strongly formal space can fail to be Alexander
formal only at a prime, and compact K\"ahler manifolds do (\S\ref{subsec:reverse}).
Section~\ref{subsec:corners} gives witnesses in all four classes cut out by the two conditions, and
\S\ref{subsec:toric} finds all four among toric complexes and characterizes the right-angled Artin
groups whose classifying spaces are strongly formal.

\begin{theorem}
\label{thm:incomp}
\leavevmode
\begin{enumerate}[label=\textup{(\roman*)},leftmargin=2.2em,itemsep=2pt]
\item \label{i1} A strongly formal space is \etaform{} over $\Q(\zeta)$ at every root of unity
$\zeta$, but need not be Alexander formal: bielliptic surfaces of types $\Z/3$ and $\Z/4$, with the class
$\nu$ of Example~\ref{ex:bielliptic}, are strongly formal and not Alexander formal.
\item \label{i2} An Alexander formal space need not be strongly formal: if $Y$ is not strongly formal, then
$S^1\times Y$, with $\nu$ the projection onto the first factor, is Alexander formal and not strongly
formal; and the mapping torus of Proposition~\ref{prop:converse}, whose first Betti number is one, is
Alexander formal and not strongly formal.
\item \label{i3} All four classes cut out by the two conditions are nonempty.
\end{enumerate}

\end{theorem}

\begin{proof}
Part~\ref{i1} is Propositions~\ref{prop:sf-char0} and \ref{prop:kahler} together with
Corollary~\ref{cor:bielliptic}, and part~\ref{i2} is Propositions~\ref{prop:dir}
and~\ref{prop:converse}. For~\ref{i3}, the two mixed classes are nonempty by \ref{i1}
and~\ref{i2}, complex tori are in both classes (\S\ref{subsec:corners}), and Examples~\ref{ex:paj}
and~\ref{ex:phi6} are in neither. 
\end{proof}

\subsection{The directional gap: Alexander formal but not strongly formal}
\label{subsec:dir}
The simplest way to be Alexander formal is to have trivial monodromy.

\begin{proposition}
\label{prop:dir}
Let $Y$ be a connected finite CW-complex, let $X=S^1\times Y$, and let $\nu$ be the
projection onto $\pi_1(S^1)=\Z$. Then $X$ is Alexander formal, and $X$ is strongly formal 
if and only if $Y$ is.
\end{proposition}

\begin{proof}
Here $X^\nu\simeq Y$, and the deck transformation is homotopic to the identity, so $t$ acts
trivially on $H_\ast(X^\nu;\k)$ for every field $\k$. Every primary part is then semisimple, and $X$
is Alexander formal by Definition~\ref{def:Aformal}. If $X$ is strongly formal, so is its retract $Y$, by
Lemma~\ref{lem:products}\ref{pr1}; conversely, $S^1$ is strongly formal by
Lemma~\ref{lem:products}\ref{pr4}, so if $Y$ is strongly formal, so is $X$, by
Lemma~\ref{lem:products}\ref{pr2}.
\end{proof}

\noindent For arrangements, the product $\A_3\times\A_1$ gives an instance of the gap with the class of
the Milnor fibration (Corollary~\ref{cor:arr-corners}); the obstruction to strong formality there is
Narkawicz's triple Massey product in a rank-two piece of $\bar A^\bullet$, which we now describe.

\begin{example}[Narkawicz]
\label{ex:nark}
Let $\A$ be the decone of the braid arrangement $\A_3$,
\[
\A\colon\quad \ell_1=x-1,\quad \ell_2=x,\quad \ell_3=y,\quad \ell_4=y-1,\quad \ell_5=x-y,
\]
five lines in $\C^2$ with two triple points, $\{2,3,5\}$ at the origin and $\{1,4,5\}$ at $(1,1)$,
and two pairs of parallels, $\{1,2\}$ and $\{3,4\}$. Put $X=\C^2\smallsetminus\bigcup_i\ell_i$, so
that $X\cong U(\A_3)$ and $M(\A_3)\cong\C^\ast\times X$. Write $A^\bullet=A^\bullet(\A)$ 
for the Orlik--Solomon algebra, generated by the logarithmic forms 
$\omega_i=\frac{1}{2\pi\sqrt{-1}}\,d\log\ell_i$; here $\dim A^1=5$ and $\dim A^2=6$.

\vspace*{1pt}
\emph{Two directions.} Fix $r\in\C$ and set
\[
\alpha_1=r(\omega_2+\omega_3-2\omega_5),\qquad
\alpha_2=r(-\omega_1-\omega_4+2\omega_5),
\]
the two rows of the rank-two matrix
$a=r\!\left(\begin{smallmatrix}0&1&1&0&-2\\-1&0&0&-1&2\end{smallmatrix}\right)\in M_{2\times5}(\C)$
paired against $(\omega_1,\dots,\omega_5)$. These are not arbitrary: $\alpha_j$ is a \emph{local
resonance class} at one of the two triple points, supported on the three lines through it, with
coefficients summing to zero. In particular neither is a multiple of the diagonal class
$\omega_1+\dots+\omega_5$, and, as we record at the end, the two cannot be detected 
together along any single epimorphism
$\nu\colon\pi_1(X)\surj\Z$.

\vspace*{1pt}
\emph{The rank-two piece of the Kohno--Pajitnov algebra.} Exponentiating gives rank-one characters
$\rho_j=\exp(\alpha_j)$, explicitly $\rho_1=(1,e^{r},e^{r},1,e^{-2r})$ and
$\rho_2=(e^{-r},1,1,e^{-r},e^{2r})$ in the coordinates $\CC=H^1(X;\C^\ast)=(\C^\ast)^5$. For all
but countably many $r$, namely whenever $e^r$ is not a root of unity, the character
$\rho=(\rho_1,\rho_2)\colon\pi_1(X)\to(\C^\ast)^2$ has Zariski dense image and
$k\mapsto\rho^k\coloneqq \rho_1^{k_1}\rho_2^{k_2}$ is injective on $\Z^2$. 
Fix such an $r$ and put
\[
\DD=\langle\rho_1,\rho_2\rangle\cong\Z^2\subseteq\CC ,
\qquad
\bar A_{\DD}=\bigoplus_{k\in\Z^2}A^\bullet(X,E_{\rho^k}) .
\]
Narkawicz's algebra $A^\bullet_a$ is the Aomoto model of $\bar A_{\DD}$:
\[
A^\bullet_a=\bigoplus_{k\in\Z^2}\big(A^\bullet,\ \lambda_k\wedge-\big)\cdot q^{k},
\qquad \lambda_k=k_1\alpha_1+k_2\alpha_2 ,
\]
and the inclusion $A^\bullet_a\inj \bar A_{\DD}$ is a quasi-isomorphism of CDGAs: summand by
summand, this is the comparison theorem for Aomoto complexes
\cite{Esnault-Schechtman-Viehweg-inv92,Schechtman-Terao-Varchenko-jpaa95}, whose
non-resonance condition holds because each sum of residues of $\lambda_k$ that it involves has
the form $rm/(2\pi\sqrt{-1})$ with $m\in\Z$, and so is a positive integer only if $e^r$ is a
root of unity. Here $q^k=q_1^{k_1}q_2^{k_2}$ 
only records the summand of the $\Z^2$-grading, which is the restriction to $\DD$ of the 
$\CC$-grading of \S\ref{subsec:KP}: $\varphi\cdot q_1^2q_2$ is the Orlik--Solomon form
$\varphi$ in the summand $k=(2,1)$, with differential 
$\lambda_{(2,1)}\wedge-=(2\alpha_1+\alpha_2)\wedge-$.

\vspace*{1pt}
\emph{The triple product.} Since $\alpha_1\wedge(\omega_2-\omega_3)=0$ and
$\alpha_2\wedge(\omega_1-\omega_4)=0$ --- each being a product of two 
elements of one local component --- the classes
\[
u=\big[(\omega_1-\omega_4)\,q_2\big],\qquad v=\big[(\omega_2-\omega_3)\,q_1\big]
\]
are defined, in $\Z^2$-degrees $(0,1)$ and $(1,0)$. One has $v^2=0$, since $v$ has degree one, and
$uv=d\big[\tfrac1r(\omega_2+\omega_3)\,q_1q_2\big]$, so the triple Massey product
$\langle u,v,v\rangle$ is defined, and, up to the sign fixed by the convention for Massey products,
\[
\big\langle\,u,\,v,\,v\,\big\rangle
=\Big[\tfrac1r(\omega_2+\omega_3)(\omega_2-\omega_3)\cdot q_1^2q_2\Big]
=-\tfrac{2}{r}\,\big[\omega_2\omega_3\cdot q_1^2q_2\big]\ \ne\ 0
\]
in $H^2(A^\bullet_a)/H^1(A^\bullet_a)^2$ \cite[Thm.~6.6.1]{Narkawicz-thesis07}.
We verify the nonvanishing, a finite computation inside $A^\bullet$. Modulo
$\omega_1\omega_2=\omega_3\omega_4=0$, together with $\omega_2\omega_3=\omega_2\omega_5-\omega_3\omega_5$
and $\omega_1\omega_4=\omega_1\omega_5-\omega_4\omega_5$, the products
$\omega_1\omega_3,\omega_1\omega_5,\omega_2\omega_4,\omega_2\omega_5,\omega_3\omega_5,\omega_4\omega_5$
form a basis of $A^2$. For $\lambda=x\alpha_1+y\alpha_2$, the coordinates of
$\lambda\omega_1,\lambda\omega_2,\lambda\omega_4,\lambda\omega_5$ on
$\omega_1\omega_3,\omega_1\omega_5,\omega_2\omega_4,\omega_4\omega_5$ form a matrix of determinant
$2r^4x^2y^2$; so if $xy\ne0$, the kernel of $\lambda\wedge-$ on $A^1$ is $\C\lambda$, and
$H^1(A^\bullet,\lambda\wedge-)=0$. Hence $H^1(A^\bullet_a)_k=0$ when $k_1k_2\ne0$; as the only degrees
on the coordinate axes adding up to $(2,1)$ are $(2,0)$ and $(0,1)$, the decomposables in degree
$(2,1)$ are the products of $H^1(A^\bullet_a)_{(2,0)}$ with $H^1(A^\bullet_a)_{(0,1)}$. These are
spanned by $[(\omega_2-\omega_3)\,q_1^2]$ and by $u$, since the kernel of $\alpha_1\wedge-$ on $A^1$ is
the local component spanned by $\omega_2-\omega_5$ and $\omega_3-\omega_5$, and similarly for
$\alpha_2$. Finally, the linear form on $A^2$ with values $2,1,2,1,-1,-1$ on the basis above vanishes
on $(2\alpha_1+\alpha_2)\wedge A^1$ and on
$(\omega_1-\omega_4)(\omega_2-\omega_3)=\omega_2\omega_4-\omega_1\omega_3$, but takes the value $2$ on
$\omega_2\omega_3$. So $[\omega_2\omega_3\cdot q_1^2q_2]$ is not decomposable, and
$\langle u,v,v\rangle$, whose indeterminacy $u\cdot H^1(A^\bullet_a)_{(2,0)}$ consists of
decomposables, does not contain $0$.
\end{example}

Two conclusions follow. For an epimorphism $\nu\colon\pi_1(X)\surj\Z$, let $\DD_\nu\subseteq\CC$ be
the group of characters that factor through $\nu$, as in Corollary~\ref{cor:allD}.

\begin{proposition}
\label{prop:nark}
Let $X=U(\A_3)$ and let $u,v$ be the classes of Example~\ref{ex:nark}.
\begin{enumerate}[label=\textup{(\alph*)},leftmargin=2.2em,itemsep=2pt]
\item \label{na1} $X$ and $M(\A_3)\cong\C^\ast\times X$ are not strongly formal.
\item \label{na2} For no epimorphism $\nu\colon\pi_1(X)\surj\Z$ is the product
$\langle u,v,v\rangle$ defined in the sub-CDGA $\bar A_{\DD_\nu}$ of $\bar A^\bullet(X)$.
\end{enumerate}
\end{proposition}

\begin{proof}
\ref{na1} A nonvanishing triple Massey product makes $\bar A_{\DD}$ non-formal, so by
Corollary~\ref{cor:allD}, applied to the subgroup $\DD$ of Example~\ref{ex:nark}, the space $X$ is not
strongly formal; nor is $M(\A_3)$, of which $X$ is a retract, by Lemma~\ref{lem:products}\ref{pr1}.

\ref{na2} The classes $v$ and $u$ lie in the $\CC$-degrees $\rho_1$ and $\rho_2$. A character of
infinite order in $\DD_\nu$ is trivial exactly on the meridians $\gamma_i$ with $\nu(\gamma_i)=0$; so
$\rho_1\in\DD_\nu$ forces $\nu(\gamma_1)=\nu(\gamma_4)=0$ and $\nu(\gamma_2)\ne0$, while
$\rho_2\in\DD_\nu$ forces $\nu(\gamma_2)=\nu(\gamma_3)=0$ and $\nu(\gamma_1)\ne0$. No $\nu$ does
both.
\end{proof}

\noindent Part~\ref{na2} is the precise sense in
which the obstruction is of rank at least two: the products $\langle\eta,\dots,\eta,a\rangle$ that govern
the characteristic-zero clauses of Alexander formality of $(X,\nu)$ live in $\bar A_{\DD_\nu}$
(Remark~\ref{rem:line}), so no invariant of a single Alexander module can carry Narkawicz's product.

Narkawicz's product is also invisible to the twisted tangent cone theorem, which concerns one
character at a time. A \emph{flat} of a central arrangement $\A$ is an intersection $Y$ of some of
its hyperplanes; its rank is $\codim Y$, and the \emph{localization} of $\A$ at $Y$ is the
central arrangement $\A_Y=\{H\in\A : H\supseteq Y\}$, of the same rank. For a flat $Y$ of rank two, let
$j_Y\colon U(\A)\inj U(\A_Y)$ be the inclusion; when $\lvert\A_Y\rvert\ge3$, the torus
$T_Y=j_Y^\ast H^1(U(\A_Y);\C^\ast)$ is the local component of $\VV^1(U(\A))$ attached to $Y$, and
$L_Y=j_Y^\ast H^1(U(\A_Y);\C)$ is its tangent space at $\one$.

\begin{proposition}
\label{prop:A3-tc}
Let $\A$ be an essential central arrangement in $\C^3$, let $U=U(\A)$, let $Y$ be a flat of rank two
with $m=\lvert\A_Y\rvert\ge3$, and let $\rho\in T_Y\smallsetminus\{\one\}$.
\begin{enumerate}[label=\textup{(\roman*)},leftmargin=2.2em,itemsep=2pt]
\item \label{tc1} The space $H^1(U;\C_\rho)$ has dimension $m-2$, and $\RR^1_k(U,\rho)$ is $L_Y$ for
$1\le k\le m-2$ and is empty for $k>m-2$.
\item \label{tc2} The torus $T_Y$ is the only irreducible component of $\VV^1(U)$ through $\rho$.
\item \label{tc3} For all $i\ge0$ and $k\ge1$, the map $\omega\mapsto\rho\exp(\omega)$ restricts to an
isomorphism from the germ of $\RR^i_k(U,\rho)$ at $0$ onto the germ of $\VV^i_k(U)$ at $\rho$.
\end{enumerate}
Hence, if all components of $\VV^1(U)$ are local, $U$ satisfies the conclusions of
Theorem~\ref{thm:jump} at every character.
\end{proposition}

\begin{proof}
Write $\rho=j_Y^\ast\sigma$, with $\sigma\ne\one$. Projection from $Y$ makes $U(\A_Y)$ a $\C$-bundle over
$\Sigma=\mathbb{P}^1\smallsetminus\{m\text{ points}\}$, so $H^1(U(\A_Y);\C_\sigma)\cong H^1(\Sigma;\C_\sigma)$
has dimension $m-2$ and $H^2(U(\A_Y);\C_\sigma)=0$. We obtain $U$ from $U(\A_Y)$ by removing the lines
of $\A$ not through $Y$ one at a time; projection from $Y$ maps each removed curve $C$ isomorphically
onto an open subset of $\Sigma$, so $H^0(C;\C_\sigma)=0$, and the Gysin sequences
$0\to H^1(U';\C_\sigma)\to H^1(U'\smallsetminus C;\C_\sigma)\to H^0(C;\C_\sigma)$ of the successive steps show that
$j_Y^\ast\colon H^1(U(\A_Y);\C_\sigma)\to H^1(U;\C_\rho)$ is an isomorphism.

Let $\omega\in H^1(U;\C)$ and $0\ne y\in H^1(U(\A_Y);\C_\sigma)$. If $\omega=j_Y^\ast\omega'\in L_Y$, then
$\omega\cup j_Y^\ast y=j_Y^\ast(\omega'\cup y)=0$. Otherwise $\omega$ has a nonzero residue $c$ along some
line $L\in\A$ not through $Y$. Let $L^\circ$ be $L$ minus its intersections with the other lines, and
write $\omega$ as $c$ times a class with a simple pole along $L$ plus a class defined on $U\cup L^\circ$.
The Gysin boundary $\partial\colon H^2(U;\C_\rho)\to H^1(L^\circ;\C_\sigma)$ is linear over classes defined
on $U\cup L^\circ$, so $\partial(\omega\cup j_Y^\ast y)=\pm c\,(j_Y^\ast y)|_{L^\circ}$, which is nonzero
since $L^\circ$ maps isomorphically onto an open subset $\Sigma'\subseteq\Sigma$ and
$H^1(\Sigma;\C_\sigma)\to H^1(\Sigma';\C_\sigma)$ is injective. As $H^0(U;\C_\rho)=0$, the twisted Aomoto
complex at $\omega$ thus has first cohomology of dimension $m-2$ if $\omega\in L_Y$, and zero
otherwise, which proves~\ref{tc1}.

By \cite[Thm.~3.1]{Libgober-topappl02}, $\mathrm{TC}_\rho(\VV^1_1(U))\subseteq\RR^1_1(U,\rho)=L_Y$; as
$\VV^1_1(U)$ is a finite union of translated subtori \cite{Budur-Wang-ens15}, every component through
$\rho$ lies in $\rho T_Y=T_Y$, which is~\ref{tc2}. With~\ref{tc1} at every point of
$T_Y\smallsetminus\{\one\}$, the germs of $\VV^1_k(U)$ and $\rho\exp(\RR^1_k(U,\rho))$ at $\rho$ are both
$T_Y$ for $k\le m-2$ and both empty for $k>m-2$. Since $U$, the complement of a curve in $\mathbb{P}^2$ 
and so a smooth affine surface, has the homotopy
type of a $2$-dimensional CW-complex \cite{Dimca-92}, both germs are empty in degree zero and above
degree two, and $\dim H^2(U;\C_\tau)=\dim H^1(U;\C_\tau)+\chi(U)$ for $\tau\ne\one$, as for the twisted
Aomoto complexes; so $\RR^2_k(U,\rho)=\RR^1_{k-\chi(U)}(U,\rho)$ and $\VV^2_k(U)=\VV^1_{k-\chi(U)}(U)$ 
near $\rho$, proving~\ref{tc3}. Finally, the germs agree at $\one$ as $U$ is formal, and at 
$\tau\ne\one$ off $\VV^1(U)$ as above.
\end{proof}

\noindent For $\A_3$, the characteristic variety
$\VV^1(U(\A_3))$ consists of the four local components and an essential component
$W=f^\ast H^1(\Sigma;\C^\ast)$, and $\VV^1_2(U(\A_3))=\{\one\}$ \cite[Ex.~3.1]{Suciu-topappl02},
\cite[Prop.~6.9]{Cohen-Suciu-mpcps99}; here $f$ is the pencil of conics through the four triple points
defined by the $3$-net. For this $(3,2)$-net (base locus the four triple points, and the three
double points of $\A_3$ off it) both sides of the inequality in \cite[Thm.~4.2]{Falk-Yuzvinsky-comp07}
equal $7$, so the three pairs of lines of the net are the only singular members of the pencil, and
$f\colon U(\A_3)\to\Sigma=\mathbb{P}^1\smallsetminus\{3\text{ points}\}$ is a smooth bundle whose fiber $C$
is a conic minus the four base points \cite[Cor.~4.4]{Falk-Yuzvinsky-comp07}; these punctures are
distinguished, so the monodromy acts trivially on $H^1(C;\C)$. For $\rho=f^\ast\sigma\ne\one$, the Leray
spectral sequence of $f$ with coefficients in $\C_\rho$ is concentrated in the column $p=1$, as
$H^0(\Sigma;\C_\sigma)=0=H^2(\Sigma;\C_\sigma)$; so $H^1(U(\A_3);\C_\rho)=f^\ast H^1(\Sigma;\C_\sigma)\cong\C$
and $\omega\cup f^\ast y=y\otimes\omega|_C$ in $H^2(U(\A_3);\C_\rho)\cong H^1(\Sigma;\C_\sigma)\otimes
H^1(C;\C)$, which for $y\ne0$ vanishes exactly when $\omega\in f^\ast H^1(\Sigma;\C)=T_\rho W$. The
argument above, with $W$ in place of $T_Y$, shows that the germs agree at $\rho$. So the conclusions of
Theorem~\ref{thm:jump} hold for $U(\A_3)$ at every character, although $U(\A_3)$ is not strongly formal
(Proposition~\ref{prop:nark}). By \cite[Thm.~11.4 and Ex.~12.4]{Suciu-conm26}, the last assertion of 
Proposition~\ref{prop:A3-tc} applies when $\A$ is decomposable over $\Q$ and the rational Alexander 
invariant of $\pi_1(M(\A))$ is $I$-adically separated, and to the non-Pappus arrangement.

\begin{remark}
\label{rem:higher-rank}
Raising the rank of the coefficients does not by itself create obstructions beyond the reach of twisted
Alexander modules. For a finite-dimensional complex representation $V$ of $\pi_1(X)$ and
$\lambda\in\C^\ast$, Proposition~\ref{prop:circle}, together with \cite[Thm.~C]{Pajitnov-topappl19},
shows that the following are equivalent: the cochains $C^\bullet(X;V_\lambda)$ are formal as a DG module
over $\bwedge(\eta)$; all products $\langle\eta,\dots,\eta,x\rangle$ of length at least three with
$x\in H^\bullet(X;V_\lambda)$ vanish; and the monodromy on the torsion of the twisted Alexander 
module $H_\ast(X^\nu;V^\vee)$ has no Jordan block of size greater than one at $\lambda$;
Remark~\ref{rem:nonab-eta} uses this to test nonabelian covers. Strong formality 
is a condition of a different kind. Along the line $\DD_\nu$ of
characters that factor through $\nu$, it makes $\bar A_{\DD_\nu}$ formal as a CDGA
(Corollary~\ref{cor:allD}), whereas the characteristic-zero clauses of Alexander formality, even when
extended to every $\lambda\in\C^\ast$ as in Remark~\ref{rem:scope}, ask only that each summand
$A^\bullet(X,E_{\rho_\lambda})$ be formal as a DG module over $\bwedge(\eta)$. 
\end{remark}

\subsection{The reverse direction: strongly formal but not Alexander formal}
\label{subsec:reverse}
Strong formality already supplies the characteristic-zero half of Alexander formality, for every $X$ and every $\nu$.

\begin{proposition}
\label{prop:sf-char0}
Let $X$ be strongly formal. Then for every $\nu$ and every $\lambda\in\C^\ast$ the monodromy on
the torsion of $H_\ast(X^\nu;\C)$ has no Jordan block of size $>1$ at $\lambda$; in particular,
$X$ is \etaform{} over $\Q(\zeta)$ at every root of unity $\zeta$. Hence $X$ is Alexander formal if and
only if it is \etaform{} over $\F_p$ for every prime $p$.
\end{proposition}

\begin{proof}
By \cite[Thm.~6.3]{Kohno-Pajitnov-cejm14} the twisted sequence $\E_r^\ast(\rho_\lambda,\eta)$
degenerates at $E_2$ for every $\lambda\in\C^\ast$, so every product $\langle\eta,\dots,\eta,a\rangle$
of length at least three with $a\in H^\bullet(X;\C_{\rho_\lambda})$ vanishes. By
\cite[Thm.~5.1]{Pajitnov-ejm17} the maximal size of a $\lambda$-Jordan block is one less than the
maximal length of a nonvanishing such product, and is therefore at most $1$. Only formality of
$\bar A^\bullet(X)$ as a CDGA enters here: \cite[Thm.~6.3]{Kohno-Pajitnov-cejm14} uses only that each
$A^\bullet(X,E_\rho)$ is a direct summand of $\bar A^\bullet(X)$ as a DG $A^\bullet(X)$-module
\cite[Def.~3.15 and Prop.~3.16]{Kohno-Pajitnov-cejm14}, so Theorem~\ref{thm:G} is not needed. The last
sentence follows from
Definition~\ref{def:Aformal}, degeneration over $\Q(\zeta)$ being equivalent to degeneration over
$\C$.
\end{proof}

\noindent In Examples~\ref{ex:phi6} and~\ref{ex:paj}, where the sixth power and the square of the
monodromy are unipotent with blocks of size two, strong formality also fails through a finite cyclic
cover: $X^{(6)}$ and $X^{(2)}$, respectively, carry nonvanishing triple Massey products
\cite[Thm.~13]{Bazzoni-Fernandez-Munoz-tams15}, and Proposition~\ref{prop:cover} applies.

\noindent Compact K\"ahler manifolds supply strongly formal spaces that are not Alexander formal.
Their strong formality \cite[Thm.~6.5]{Kohno-Pajitnov-cejm14} comes from Simpson's theory: every rank-one
character is semisimple, and the twisted sequences degenerate at semisimple characters.

\begin{proposition}
\label{prop:kahler}
Let $X$ be a compact K\"ahler manifold and $\nu\colon\pi_1(X)\surj\Z$. Then $X$ is strongly formal,
and it is Alexander formal if and only if it is \etaform{} over $\F_p$ for every prime $p$. If
$H_\ast(X^\nu;\Z)$ is finitely generated as an abelian group and, for some prime power 
$q=p^a>2$, a primitive $q$-th root of unity is an eigenvalue of $h_\ast$ on $H_\ast(X^\nu;\Q)$, 
then $X$ is not \etaform{} over $\F_p$; in particular, $X$ is strongly formal and not Alexander formal.
\end{proposition}

\begin{proof}
Strong formality is \cite[Thm.~6.5]{Kohno-Pajitnov-cejm14}, the first equivalence is
Proposition~\ref{prop:sf-char0}, and the last assertion is Corollary~\ref{cor:collision}.
\end{proof}

\begin{example}
\label{ex:bielliptic}
Let $X$ be a bielliptic surface. By the classification of bielliptic surfaces \cite[List~VI.20]{Beauville-96},
$X=(E_1\times E_2)/G$, with $G$ one of seven finite abelian groups. We 
consider the four types with $G=\Z/d$ cyclic, $d\in\{2,3,4,6\}$, acting on $E_1$ 
by a free translation of order $d$ and on $E_2$ by an automorphism $\zeta$ of 
order $d$ fixing $0$. Then $X$ is smooth projective, $b_1(X)=2$, and the Albanese 
map $\pi\colon X\to B\coloneqq E_1/G$ is a holomorphic $E_2$-bundle with monodromy 
$\mu\colon\pi_1(B)=\Z^2\surj\langle\zeta\rangle=\Z/d$. Pick a basis $e_1,e_2$ of 
$\pi_1(B)$ with $\mu(e_1)=\zeta$ and $\mu(e_2)=1$, and set $\nu=\alpha\circ\pi_\#$, 
where $\alpha(e_1)=1$ and $\alpha(e_2)=0$; this class detects the monodromy.

Let $X^\nu=X\times_B B^\alpha$. This is an $E_2$-bundle over the cylinder $B^\alpha\simeq S^1$ whose
monodromy around the circle is $\mu(e_2)=1$, so the bundle is trivial, $X^\nu\simeq E_2\times S^1$,
and the deck transformation acts by $\zeta\times\id$. For $d\in\{3,4,6\}$ the map $\zeta_\ast$ acts on
$H_1(E_2;\Z)=\Z^2$ with characteristic polynomial $\Phi_d$, and $H_1(E_2;\Z)\cong\LA/(\Phi_d)$, the
ring $\Z[\zeta_d]$ being a principal ideal domain; for $d=2$ we have $\zeta_\ast=-1$ and
$H_1(E_2;\Z)\cong(\LA/(\Phi_2))^2$. Since $\det\zeta_\ast=1$, the K\"unneth theorem gives, as
$\LA$-modules and for $d\in\{3,4,6\}$,
\[
H_0(X^\nu)=\tfrac{\LA}{(t-1)},\quad
H_1(X^\nu)=\tfrac{\LA}{(\Phi_d)}\oplus\tfrac{\LA}{(t-1)},\quad
H_2(X^\nu)=\tfrac{\LA}{(t-1)}\oplus\tfrac{\LA}{(\Phi_d)},\quad
H_3(X^\nu)=\tfrac{\LA}{(t-1)} ,
\]
with $\LA/(\Phi_d)$ replaced by $(\LA/(\Phi_2))^2$ for $d=2$.
\end{example}

\begin{corollary}
\label{cor:bielliptic}
Let $X$ be a bielliptic surface of type $\Z/d$, with the class $\nu$ of Example~\ref{ex:bielliptic}.
Then $X$ is strongly formal, and it is Alexander formal if and only if $d\in\{2,6\}$. For $d=3$ it is
not \etaform{} over $\F_3$, and for $d=4$ it is not \etaform{} over $\F_2$.
\end{corollary}

\begin{proof}
Strong formality, and the characteristic-zero clauses, hold by Proposition~\ref{prop:kahler}. Since
$X^\nu$ is a finite complex with torsion-free homology, $H_\ast(X^\nu;\F_p)=H_\ast(X^\nu;\Z)\otimes\F_p$.
For $d=3$ and $p=3$, and for $d=4$ and $p=2$, one has $\Phi_d\equiv(t-1)^2\bmod p$, and
$\F_p[t]/(\Phi_d)$ is a Jordan block of size two at the eigenvalue $1$. For $d=2$, the
$(t-1)$-primary part of $(\F_p[t]/(t+1))^2$ is zero for $p$ odd and is $\F_2^2$ with trivial action
for $p=2$. For $d=6$ we have $\Phi_6(1)=1$, so the $(t-1)$-primary part of $\F_p[t]/(\Phi_6)$ is zero for
every $p$.
\end{proof}

\noindent The failures for $d=3,4$ are those predicted by Proposition~\ref{prop:kahler}: the
eigenvalues of order $3$ and $4$ collide with $1$ modulo $3$ and $2$. In type $\Z/6$ there is a
collision as well, $\Phi_6\equiv(t+1)^2\bmod 3$, but it produces a Jordan block at the closed point
$(3,t+1)$, which Definition~\ref{def:Aformal} does not examine (Remark~\ref{rem:scope}). Finally, for
any $d$, choosing instead $\nu=\alpha'\circ\pi_\#$ with $\alpha'(e_1)=0$ and $\alpha'(e_2)=1$ gives an
infinite cyclic cover whose deck transformation is homotopic to the identity, hence an Alexander formal
pair: the failure is tied to a class that detects the monodromy.

\subsection{The four classes}
\label{subsec:corners}
The two mixed classes have been treated above; for the class of spaces that are neither Alexander
formal nor strongly formal, the simplest witness is a mapping torus due to Pajitnov.

\begin{example}
\label{ex:paj}
Fix $n\ge1$, let $\Sigma_n$ be the closed orientable surface of genus $n$, and let
$\varphi\colon\Sigma_n\to\Sigma_n$ be
a diffeomorphism inducing
\[
\varphi^\ast=\begin{pmatrix}-I_n&I_n\\0&-I_n\end{pmatrix}
\qquad\text{on } H^1(\Sigma_n;\Z)=\Z^{2n};
\]
this matrix is symplectic, and the mapping class group surjects onto
$\Sp_{2n}(\Z)$, so such a $\varphi$ exists. Let $X$ be its mapping torus and
$\nu\colon\pi_1(X)\surj\Z$ the fibration class, so that $X^\nu\simeq\Sigma_n$ and $t$ acts on
$H_q(X^\nu;\C)=H_q(\Sigma_n;\C)$ as $\varphi_\ast$. Since $(\varphi^\ast+I)^2=0\ne\varphi^\ast+I$, 
the monodromy on $H^1$ is $n$ Jordan blocks of size $2$ at the eigenvalue $-1$. Therefore:
\begin{enumerate}[label=\textup{(\roman*)},leftmargin=2.2em,itemsep=2pt]
\item $X$ is \emph{not} Alexander formal: its rational Jordan blocks of size $2$ at $\zeta=-1$ make it fail
$\eta$-formality over $\Q(\zeta_2)=\Q$ after the twist $t\mapsto-t$. This is
a failure at an eigenvalue $\zeta\ne1$ in characteristic zero, which neither a compact K\"ahler
manifold nor an arrangement complement with the class of its Milnor fibration can exhibit, their
monodromy being semisimple over $\C$.
\item $X$ is \emph{not} strongly formal, by Proposition~\ref{prop:sf-char0}, which excludes such
blocks at every eigenvalue.
\end{enumerate}
So $X$ is neither Alexander formal nor strongly formal. The example is Pajitnov's
\cite[\S4.2]{Pajitnov-ejm17}. He observes further that $\varphi^\ast-1$ is invertible on $H^1(\Sigma_n;\C)$,
its only eigenvalue being $-2$, so the Wang sequence gives $H^\ast(X;\C)\cong H^\ast(S^1\times
S^2;\C)$, and $X$, a closed orientable $3$-manifold with $b_1(X)=1$, is formal by
Lemma~\ref{lem:b1one}. The obstruction to strong formality is thus read
directly off the Alexander module, from its rational Jordan blocks at $-1$.
\end{example}

\begin{example}
\label{ex:kodaira}
In genus one, Example~\ref{ex:paj} is one of a family of torus bundles indexed by Kodaira's list of
singular fibers of elliptic fibrations \cite{Kodaira-annals63}. For $A\in\SL_2(\Z)$, let $T_A$ be the
mapping torus of the diffeomorphism of $T^2=\R^2/\Z^2$ induced by $A$, and let $\nu$ be the class of
the bundle projection $T_A\to S^1$. Then $T_A^\nu\simeq T^2$, the monodromy acts by $A$ on
$H_1(T^2;\Z)=\Z^2$ and trivially on $H_0$ and $H_2$, and the $d$-fold cyclic cover of $T_A$ along $\nu$
is $T_{A^d}$. By the Wang sequence, $b_1(T_A)=1+\dim\ker(A-1)$, so $T_A$ is formal by
Lemma~\ref{lem:b1one} unless $A$ is unipotent. Up to conjugacy and orientation, the monodromy of a
singular fiber of type $\mathrm{I}_b$ or $\mathrm{I}_b^\ast$, with $b\ge1$, is $U^b$ or $-U^b$, where
$U=\left(\begin{smallmatrix}1&1\\0&1\end{smallmatrix}\right)$; it has order six for the types
$\mathrm{II}$ and $\mathrm{II}^\ast$, order four for the types $\mathrm{III}$ and $\mathrm{III}^\ast$,
order three for the types $\mathrm{IV}$ and $\mathrm{IV}^\ast$, and order two for the type
$\mathrm{I}_0^\ast$. Accordingly:
\begin{enumerate}[label=\textup{(\roman*)},leftmargin=2.2em,itemsep=2pt]
\item For type $\mathrm{I}_b$, the monodromy has a Jordan block of size two at the eigenvalue $1$, so
$\pi_1(T_A)$ is not $1$-formal \cite[Prop.~4.1]{Papadima-Suciu-forum10}, and $T_A$ is not formal. For
$b=1$, it is the Heisenberg nilmanifold.
\item For type $\mathrm{I}_b^\ast$, the bundle $T_A$ is formal, but its Jordan block of size two at $-1$
makes it fail \ref{AF1} at $\zeta=-1$, and its double cover $T_{A^2}$, of type $\mathrm{I}_{2b}$, is not
formal. So, like Example~\ref{ex:paj}, which in genus one is the case $b=1$, it is neither Alexander
formal nor strongly formal (Proposition~\ref{prop:sf-char0}).
\item For the types of finite order $d$, the monodromy is semisimple over $\C$, so $T_A$ satisfies
\ref{AF1}, and $T_A$ and all of its cyclic covers along $\nu$ are formal. As $\LA$-modules,
$H_1(T^2;\Z)\cong\LA/(\Phi_d)$ for $d\in\{3,4,6\}$, and $(\LA/(\Phi_2))^2$ for $d=2$, as in
Example~\ref{ex:bielliptic}. By Lemma~\ref{lem:collision}, the types $\mathrm{III}$ and $\mathrm{III}^\ast$
fail \ref{AF2} at $p=2$, and the types $\mathrm{IV}$ and $\mathrm{IV}^\ast$ at $p=3$; at every other
prime, $\Phi_4(1)=2$ and $\Phi_3(1)=3$ are units, and there is no Jordan block at $1$. The types
$\mathrm{II}$ and $\mathrm{II}^\ast$, and also the type $\mathrm{I}_0^\ast$, are Alexander formal, since
$\Phi_6(1)=1$ and $-1$ acts trivially modulo $2$.
\end{enumerate}
\end{example}

\begin{example}
\label{ex:tubes}
Punctured versions of these bundles appear around the irregular fibers of polynomial maps
$f\colon\C^2\to\C$, that is, the fibers over the values near which $f$ fails to be locally trivial. If
$c$ is such a value and $\varepsilon>0$ is small, the tube $f^{-1}(D^\ast)$ over the punctured disk
$D^\ast$ of radius $\varepsilon$ about $c$ is a fiber bundle over $D^\ast$, with fiber a regular fiber
$F$ of $f$; up to homotopy, it is the mapping torus of the local monodromy of $f$ at $c$.

For $f=x^3+y^3+xy$, the curves $f=s$ are the affine parts of the cubics $x^3+y^3+xyz=sz^3$, which pass
through the three points of the line $z=0$ where $x^3+y^3=0$, and are smooth there; the only singular
members with $s\ne\infty$ are the nodal cubic $s=0$ and the triangle of lines $s=1/27$. So $F$ is an
elliptic curve with three points removed, the member $s=0$ is a fiber of type $\mathrm{I}_1$, and the
tube around $f^{-1}(0)$ is homotopy equivalent to the Heisenberg nilmanifold with three circles removed,
the sections of the bundle through these three points. This accounts for the Jordan block of size two at
the eigenvalue $1$ on $H_1(F;\C)$ found by Dimca \cite[Ex.~6.1]{Papadima-Suciu-forum10}. The tube is
not formal, and even its fundamental group is not $1$-formal \cite[Prop.~3.2]{Papadima-Suciu-forum10},
although the complement of the link of the node, a Hopf link, is formal.

The Brian\c{c}on polynomial $f_B$ has no critical points and exactly two irregular values, $0$ and
$-16/9$, and its regular fiber is a torus with three punctures; its local monodromy at $0$ is the product
of the Dehn twists about two disjoint curves that cut the fiber into two pieces of genus zero, one of
them an annulus with one puncture, and its local monodromy at $-16/9$ is isotopic to a map of order $6$
\cite[\S5]{Neumann-Norbury-mathann00}. In the closed torus both pieces become annuli, so the two curves
are parallel and essential, and the tube around $f_B^{-1}(0)$ is homotopy equivalent to the torus bundle
of type $\mathrm{I}_2$ with three circles removed; again its fundamental group is not $1$-formal.
\end{example}

A complex torus is compact K\"ahler, hence strongly formal, and for every $\nu$ its infinite cyclic
cover is the product of a line and a torus, on which the deck transformation acts by a translation, 
hence trivially in homology, so it is Alexander formal. The table in \S\ref{subsec:results} lists the 
witnesses in all four classes, among them the arrangement complements of Corollary~\ref{cor:arr-corners}.

\begin{remark}
\label{rem:corners}
Yoshinaga's icosidodecahedral arrangement is not Alexander formal (Theorem~\ref{thm:arr}), but
we do not know whether its complement $X$ is strongly formal: its algebraic monodromy is trivial on
$H_1(F;\Q)$, and an arrangement of rank $3$ with this property has formal Milnor fiber
(\cite[Cor.~6.2(2)]{Suciu-revroum24}, from \cite{Dimca-Papadima-pisa11}), so
Corollary~\ref{cor:milnor-strong} does not apply through the Milnor fiber. By Theorem~\ref{thm:jump},
a character $\rho$ with $\mathrm{TC}_\rho\big(\VV^i_k(X)\big)\ne\RR^i_k(X,\rho)$ would settle the question
negatively; a computer search finds none among the characters of order at most three, 
and by Proposition~\ref{prop:A3-tc} none lies on a local component. Characters at which 
two components of $\VV^1$ meet have finite order \cite[Thm.~4.2(3)]{Dimca-Papadima-Suciu-imrn08}, 
so a search over characters of finite order can reach them.
\end{remark}

\subsection{Toric complexes}
\label{subsec:toric}
Let $L$ be a finite simplicial complex on a vertex set $V$. The \emph{toric complex} $T_L$ is the
subcomplex of the torus $(S^1)^V$, with its product cell structure, consisting of the cells $c_\sigma$
indexed by the simplices $\sigma\in L$, the empty simplex included, with $\dim c_\sigma=\lvert\sigma\rvert$
\cite[\S\S2.1--2.2]{Papadima-Suciu-adv09}. Its fundamental group is the right-angled Artin group $G_L$
of the graph $L^{(1)}$, it is aspherical when $L$ is a flag complex, and its cohomology ring is the
exterior face ring of $L$ \cite[(11)--(13)]{Papadima-Suciu-adv09}. A join gives a product,
$T_{K\ast L}=T_K\times T_L$, and for $W\subseteq V$ the coordinate projection retracts $T_L$ onto
$T_{L_W}$, where $L_W$ is the full subcomplex on $W$. Every toric complex is formal
\cite[Thm.~5.5 and Rem.~5.7]{Notbohm-Ray-agt05}, so $T_L$ is \etaform{} over
$\Q$ for every $\nu$ (Proposition~\ref{thm:assembly}\ref{a1}).

An epimorphism $\nu\colon G_L\surj\Z$ is a vector of weights $m_v=\nu(v)$ for $v\in V$, with
$\gcd_v m_v=1$. The equivariant cellular chain complex of $T_L^\nu$ is $\LA\otimes C_\bullet(T_L)$, with
\[
\partial(c_\sigma)=\sum_{r=1}^{k}(-1)^{r-1}\big(t^{m_{v_r}}-1\big)\,c_{\sigma\smallsetminus\{v_r\}},
\qquad \sigma=\{v_1<\dots<v_k\}
\]
\cite[(14)]{Papadima-Suciu-adv09}. Papadima and Suciu compute its homology over a power series ring
\cite[Prop.~5.6, Thm.~5.9]{Papadima-Suciu-adv09}. Fix a field $\k$ and an irreducible $f\in\k[t]$ other
than $t$; for $m_v\ne0$ let $b_v$ be the multiplicity of $f$ in $t^{m_v}-1$ and $\xi_v=s^{b_v}$, and let
$\xi_v=0$ if $m_v=0$. Let $\partial_\xi$ be the boundary map over $\k\llbracket s\rrbracket$, with $s$
standing for $f$, obtained by replacing each $t^{m_v}-1$ with $\xi_v$. Then the $f$-primary part of the
torsion of $H_i(T_L^\nu;\k)$ is $\bigoplus_j\k[t^{\pm1}]/(f^{a_j})$, where the $s^{a_j}$ are the
elementary divisors of
$\partial_\xi\colon\k\llbracket s\rrbracket\otimes C_{i+1}\to\k\llbracket s\rrbracket\otimes C_i$ that are 
neither zero nor units. In particular the monodromy is quasi-unipotent \cite[Prop.~5.7]{Papadima-Suciu-adv09}.

\begin{proposition}
\label{prop:toric}
Let $\nu\colon G_L\surj\Z$ have weights $m_v$, and let $\k$ be a field of characteristic $p\ge0$.
\begin{enumerate}[label=\textup{(\roman*)},leftmargin=2.2em,itemsep=2pt]
\item \label{to1} If $p$ divides no nonzero weight, which is automatic for $p=0$, then the
$(t-1)$-primary part of the torsion of $H_\ast(T_L^\nu;\k)$ is semisimple. In particular, $T_L$ is
\etaform{} over $\F_p$ for every prime $p$ that divides no nonzero weight.
\item \label{to2} If every weight lies in $\{0,1,-1\}$, then for every field $\k$ the torsion of
$H_\ast(T_L^\nu;\k)$ is a direct sum of copies of $\k[t^{\pm1}]/(t-1)$; in particular, $(T_L,\nu)$ is
Alexander formal.
\end{enumerate}
\end{proposition}

\begin{proof}
\ref{to1} For $f=t-1$ and $m_v\ne0$, the polynomial $t^{m_v}-1$ has a simple root at $1$, since
$p\nmid m_v$; so $b_v=1$, and every $\xi_v$ is $0$ or $s$. Hence $\partial_\xi=s\,\partial'$ with
$\partial'$ a matrix over $\k$, and diagonalizing $\partial'$ over $\k$ shows that every nonzero
elementary divisor of $\partial_\xi$ equals $s$: every Jordan block at the eigenvalue $1$ has size one.
\ref{to2} No irreducible $f\ne t-1$ divides $t^{\pm1}-1$, so the $f$-primary parts vanish, and
part~\ref{to1} applies in every characteristic. Alexander formality follows from
Definition~\ref{def:Aformal}.
\end{proof}

\noindent Part~\ref{to2} covers the Bestvina--Brady class, $\nu(v)=1$ for all $v$
\cite[Cor.~5.10]{Papadima-Suciu-adv09}. Outside the range of Proposition~\ref{prop:toric}, both kinds
of failure occur.

\begin{example}
\label{ex:toric-pencil}
Let $L$ be the path $a$---$b$---$c$, so that $T_L=S^1\times(S^1\vee S^1)$, with weights 
$1$, $m$, $1$, where $m\ge2$. By the Künneth formula over $\LA_\k=\k[t^{\pm1}]$,
\[
H_1(T_L^\nu;\k)\cong\LA_\k/(t^m-1)\oplus\LA_\k/(t-1),\qquad H_2(T_L^\nu;\k)=0 .
\]
Over $\Q$ the torsion is semisimple. If $p^s$ is the exact power of a prime $p$ dividing $m$, with
$s\ge1$, then $t^m-1\equiv(t^{m/p^s}-1)^{p^s}\bmod p$, so $H_1(T_L^\nu;\F_p)$ has a Jordan 
block of size $p^s$ at $1$, and $(T_L,\nu)$ is not \etaform{} over $\F_p$. Yet $T_L$ is strongly 
formal, by Lemma~\ref{lem:products}\ref{pr2},\ref{pr4}, and formal over $\Z$ (Example~\ref{ex:three-lines}). 
For $m=3$, the identification $M\cong U\times\C^\ast$ of the proof of Theorem~\ref{thm:arr} makes 
$T_L$ homotopy equivalent to the complement of three lines through the origin of $\C^2$, 
with $\nu$ the class of its Milnor fibration (Corollary~\ref{cor:arr-corners}\ref{co2}).
\end{example}

\begin{example}
\label{ex:toric-P4}
Let $L$ be the path $a$---$b$---$c$---$d$. It is a flag complex, so $T_L$ is a $2$-dimensional
$K(G_L,1)$, and $T_L$ is formal. Give $a$ and $d$ the weight $1$ and $b$ and $c$ the weight $2$. In
the bases $(c_{ab},c_{bc},c_{cd})$ and $(c_a,c_b,c_c,c_d)$, the boundary map $\partial_2$ has matrix
\[
\begin{pmatrix}
-(t^2-1)&0&0\\
t-1&-(t^2-1)&0\\
0&t^2-1&-(t-1)\\
0&0&t^2-1
\end{pmatrix},
\]
and $\partial_1=(t-1,\ t^2-1,\ t^2-1,\ t-1)$ has rank one; so $H_2(T_L^\nu;\Q)=0$, and
$H_1(T_L^\nu;\Q)$ is the torsion of $\coker\partial_2$. Over $\Q[t^{\pm1}]$, the greatest common
divisors of the $1\times1$, $2\times2$, and $3\times3$ minors
of $\partial_2$ are $t-1$, $(t-1)^2$, and $(t-1)^3(t+1)^2$, so the invariant factors are $t-1$, $t-1$,
and $(t-1)(t+1)^2$, and
\[
H_1(T_L^\nu;\Q)\cong\big(\Q[t^{\pm1}]/(t-1)\big)^{3}\oplus\Q[t^{\pm1}]/(t+1)^2 .
\]
The monodromy has a Jordan block of size two at $\zeta=-1$. So $(T_L,\nu)$ is not \etaform{} over $\Q$
at $-1$, hence not Alexander formal, and $T_L$ is not strongly formal, by
Proposition~\ref{prop:sf-char0}. By the dictionary of \S\ref{sec:setup}, the block corresponds to a
nonvanishing twisted triple Massey product $\langle\eta,\eta,a\rangle$ with $a\in H^1(T_L;\C_\rho)$, where
$\rho$ takes the value $-1$ on the end vertices and $1$ on the middle ones. With the weights
$1,1,1,1$ instead, $(T_L,\nu)$ is Alexander formal, by Proposition~\ref{prop:toric}\ref{to2}.
\end{example}

\begin{corollary}
\label{cor:toric-P4}
\begin{enumerate}[label=\textup{(\roman*)},leftmargin=2.2em,itemsep=2pt]
\item \label{tp1} If $L$ has a full subcomplex isomorphic to the path on four vertices, then $T_L$ is
formal but not strongly formal.
\item \label{tp2} Let $\Gamma$ be a finite graph and $\Delta_\Gamma$ its flag complex, so that
$T_{\Delta_\Gamma}$ is a classifying space for the right-angled Artin group of $\Gamma$. Then
$T_{\Delta_\Gamma}$ is strongly formal if and only if $\Gamma$ has no induced path on four vertices.
\end{enumerate}
\end{corollary}

\begin{proof}
\ref{tp1} The complex $T_L$ retracts onto the toric complex of that full subcomplex, which is not
strongly formal by Example~\ref{ex:toric-P4}, so $T_L$ is not strongly formal, by
Lemma~\ref{lem:products}\ref{pr1}.

\ref{tp2} An induced path on four vertices of $\Gamma$ spans a full subcomplex of $\Delta_\Gamma$
isomorphic to that path, so one implication is part~\ref{tp1}. Conversely, a graph without an induced
path on four vertices is a cograph: it is obtained from single vertices by disjoint unions and joins
\cite{Corneil-Lerchs-Burlingham-dam81}. The flag complex of a disjoint union or a join of graphs 
is the disjoint union or the join of their flag complexes, and, directly from the definition,
$T_{K\sqcup L}=T_K\vee T_L$ and $T_{K\ast L}=T_K\times T_L$. So $T_{\Delta_\Gamma}$ 
is obtained from circles by wedges and products, and it is strongly formal by
Lemma~\ref{lem:products}\ref{pr2}--\ref{pr4}.
\end{proof}

\noindent For complexes $L$ that are not flag, part~\ref{tp1} still applies, but we do not know whether
its converse holds in general. It does hold for the skeleta of a simplex.

\begin{example}
\label{ex:skeleta}
Let $L$ be the $(d-1)$-skeleton of the simplex on a finite set $V$, with $1\le d\le\lvert V\rvert$, so
that $T_L$ is the $d$-skeleton of the torus $(S^1)^V$; for $2\le d<\lvert V\rvert$ the complex $L$ is not a flag
complex. Its exterior face ring is the exterior algebra on $V$, truncated above degree $d$. By
\cite[Lem.~3.11]{Papadima-Suciu-adv09}, for $\rho\in(\C^\ast)^V$ the dimension of $H_i(T_L;\C_\rho)$
is that of the $i$-th cohomology of this ring under multiplication by $z_\rho=\sum_v(\rho(v)-1)v$. If
$\rho\ne\one$, then $z_\rho\ne0$, and multiplication by $z_\rho$ is exact on the whole exterior
algebra; so the cohomology vanishes in degrees below $d$, where the truncation has no effect. Thus
$H_i(T_L;\C_\rho)=0$ for $i<d$, and so is $H^i(T_L;\C_\rho)$, which is dual to
$H_i(T_L;\C_{\rho^{-1}})$. As $T_L$ is formal and has dimension $d$, it is strongly formal, by
Lemma~\ref{lem:products}\ref{pr4}. Since $L$ has no full subcomplex isomorphic to the path on four
vertices, this confirms the converse of Corollary~\ref{cor:toric-P4}\ref{tp1} for these complexes.
\end{example}

\noindent Toric complexes thus occupy the four classes of Theorem~\ref{thm:incomp}: the torus $T^n$
lies in both for every $\nu$, and Examples~\ref{ex:toric-pencil} and~\ref{ex:toric-P4} supply the other
three.

\section{Formality properties of arrangement complements}
\label{sec:functorial}

Let $\A$ be a central arrangement of $n$ hyperplanes, defined by $Q=\prod_i\ell_i$, with
complement $M=M(\A)$, projectivized complement $U=U(\A)$, and Milnor fiber $F=F(\A)=Q^{-1}(1)$. The
cohomology ring $H^\bullet(M;\Z)$ is the Orlik--Solomon algebra $A^\bullet=A^\bullet(\A)$, generated
by the classes $\omega_1,\dots,\omega_n$ of the logarithmic forms
$\frac{1}{2\pi\sqrt{-1}}\,d\log\ell_i$ \cite{Brieskorn-71}; it is torsion-free, and determined by
the intersection lattice. The diagonal class $\nu=\omega_1+\dots+\omega_n$ is the class of the Milnor fibration
$Q\colon M\to\C^\ast$. Throughout this section $X=M$. Then $X^\nu\simeq F$, a finite
complex, and $t$ acts on $A_\nu=H_\ast(F;\Z)$ as the monodromy $h_\ast$, induced by
$x\mapsto e^{2\pi\sqrt{-1}/n}x$ and of order dividing $n$.

This section proves Theorem~\ref{thmx:arr}: Alexander formality of $M(\A)$ is decided by the sizes of
the irreducible factors of $\A$ (\S\ref{subsec:criterion}). It then realizes the four classes of
Theorem~\ref{thm:incomp} by arrangement complements, and relates the obstructions to strong formality
found there to the Massey products that govern Alexander formality (\S\ref{subsec:arr-corners}).
Finally, \S\ref{subsec:localization} propagates the failure of strong formality for $\A_3$ along
localizations, and \S\ref{subsec:graphic} shows that the complements of generic arrangements, and of
many graphic arrangements, are strongly formal, deciding strong formality for the graphic arrangements
of chordal graphs.

\subsection{The criterion for Alexander formality}
\label{subsec:criterion}
The characteristic-zero clauses of Definition~\ref{def:Aformal} are automatic here: since $h_\ast^n=\id$,
the monodromy is semisimple over $\C$, so $M$ is \etaform{} over $\Q(\zeta)$ at every root of unity
$\zeta$ (\S\ref{sec:setup}).
The modular clauses are decided by the Aomoto complex, by Proposition~\ref{prop:finite-type}\ref{ft1}, 
and away from the primes dividing $n$ that complex is acyclic for a simple reason.

\begin{lemma}
\label{lem:homotopy}
Let $\k$ be a field in which $n$ is invertible. Then the Aomoto complex 
$(A^\bullet\otimes\k,\nu\cdot-)$ is acyclic.
\end{lemma}

\begin{proof}
The exterior algebra on $e_1,\dots,e_n$ carries the derivation $\partial$ of degree $-1$ with
$\partial e_i=1$. The Orlik--Solomon ideal is generated by the elements $\partial e_S$, and it is
stable under $\partial$, because $\partial$ is a derivation and $\partial^2=0$; so $\partial$ descends
to $A^\bullet\otimes\k$. For every $a$ one has $\partial(\nu a)+\nu\,\partial a=\partial(\nu)\,a=na$,
so $n^{-1}\partial$ is a contracting homotopy.
\end{proof}

At the primes dividing $n$ the answer depends on the decomposition of $\A$ into irreducible factors.
Recall that $\A$ is \emph{reducible} if, after a linear change of coordinates, it is a product
$\A'\times\A''$ of arrangements in complementary subspaces, and that every central arrangement is,
up to a factor $\C^k$ that does not affect the complement up to homotopy, a product of irreducible
ones, uniquely up to order.

\begin{theorem}
\label{thm:arr}
Let $\A=\A_1\times\dots\times\A_k$ be the decomposition of $\A$ into irreducible arrangements, and
let $n_i=\lvert\A_i\rvert$, so that $n=n_1+\dots+n_k$. For a prime $p$ the following are equivalent:
\begin{enumerate}[label=\textup{(\alph*)},leftmargin=2.2em,itemsep=2pt]
\item \label{ar1} 
$M(\A)$ is \etaform{} over $\F_p$.
\item \label{ar2} 
The Aomoto complex $(A^\bullet\otimes\F_p,\nu\cdot-)$ is acyclic.
\item \label{ar3} 
$p\nmid n_i$ for some $i$.
\end{enumerate}
Consequently $M(\A)$ is Alexander formal if and only if, for every prime $p$ dividing $n$, some $n_i$ is
prime to $p$, that is, if and only if $\gcd(n_1,\dots,n_k)=1$. In particular, if $\A$ is irreducible and $n\ge2$, 
then $M(\A)$ is not \etaform{} over $\F_p$ for any prime $p\mid n$, and so is not Alexander formal.
\end{theorem}

\begin{proof}
The equivalence of \ref{ar1} and \ref{ar2} is Proposition~\ref{prop:finite-type}\ref{ft1}, since $X^\nu\simeq F$
is a finite complex.

For \ref{ar2}$\Leftrightarrow$\ref{ar3}, note first that $M(\A)=\prod_iM(\A_i)$, and that, the rings
$H^\bullet(M(\A_i);\Z)$ being torsion-free, the K\"unneth theorem gives
$H^\bullet(M;\F_p)=\bigotimes_iH^\bullet(M(\A_i);\F_p)$ as algebras, with $\nu=\sum_i\nu_i$
for the diagonal classes $\nu_i$ of the factors. Since the $\nu_i$ have odd degree,
multiplication by $\nu$ is the tensor product differential; so the Aomoto complex of $\A$ over
$\F_p$ is the tensor product of those of the $\A_i$, and by the K\"unneth formula over a field it is
acyclic if and only if one of the factors is. It therefore suffices to show that, for irreducible
$\A_i$, the Aomoto complex of $\A_i$ over $\F_p$ is acyclic if and only if $p\nmid n_i$. If
$p\nmid n_i$, this is Lemma~\ref{lem:homotopy}. Suppose $p\mid n_i$. Then $n_i\ge2$, so $\A_i$ has
rank $r_i\ge2$. Choose a hyperplane $\ker\ell_{j_0}$ of $\A_i$. The map $x\mapsto([x],\ell_{j_0}(x))$ identifies
$M(\A_i)$ with $U_i\times\C^\ast$, where $U_i=U(\A_i)$; under this identification $\omega_{j_0}$
corresponds to $1\otimes e$, with $e$ the generator of $H^1(\C^\ast;\Z)$, and every other $\omega_j$
to $\bar\omega_j\otimes1+1\otimes e$, where $\bar\omega_j\in H^1(U_i;\Z)$ is the class of
$d\log(\ell_j/\ell_{j_0})$. Hence $\nu_i=\bar\nu_i\otimes1+n_i(1\otimes e)$, with $\bar\nu_i$
the sum of the $\bar\omega_j$, and modulo $p$ the Aomoto complex of $\A_i$ becomes
\[
\big(H^\bullet(U_i;\F_p),\bar\nu_i\cdot-\big)\otimes\big(\bwedge(e),0\big) .
\]
The first factor has Euler characteristic $\chi(U_i)$, because $H^\bullet(U_i;\Z)$ is torsion-free.
Write $\pi(\A_i,t)$ for the Poincar\'e polynomial of $M(\A_i)$, so that
$\pi(\A_i,t)=(1+t)P(U_i,t)$; a direct comparison with the characteristic polynomial shows that
$\chi(U_i)=P(U_i,-1)$ equals $(-1)^{r_i-1}$ times Crapo's beta invariant of the matroid of $\A_i$.
Since $\A_i$ is irreducible, this matroid is connected, and its beta invariant is nonzero
\cite{Crapo-jct67}. So the first factor is not acyclic, and neither is the tensor product.

The consequences follow: the characteristic-zero clauses hold, the clause at $p$ is \ref{ar3}, and
\ref{ar3} holds automatically when $p\nmid n$.
\end{proof}

\begin{remark}
\label{rem:arr-collision}
The Euler characteristic that obstructs acyclicity can also be read off the monodromy. Since $\mu_n$
acts freely on $F$ with quotient $U$, for each $n$-th root of unity $\zeta$ one has
$\sum_q(-1)^q\dim_\C H_q(F;\C)_\zeta=\chi(U)$, the subscript denoting the $\zeta$-eigenspace. So for
irreducible $\A$ every $n$-th root of unity is an eigenvalue of the monodromy in some degree, and when
$n$ has a prime-power divisor $q=p^a>2$, Lemma~\ref{lem:collision} converts the eigenvalue of order $q$
into a Jordan block of size greater than one at $t=1$ modulo $p$, as for the bielliptic surfaces of
types $\Z/3$ and $\Z/4$ (Example~\ref{ex:bielliptic}). When $2$ divides $n$ exactly, the only
eigenvalues of $2$-power order are $\pm1$, and the failure at $p=2$ given by Theorem~\ref{thm:arr}
comes instead from the integral lattice, through the eigenvalues $\pm1$ or the $2$-torsion of
$H_\ast(F;\Z)$ (compare Remark~\ref{rem:collision}). Since $t^n=1$ on $H_\ast(F;\F_p)$, every Jordan
block at $t=1$ modulo $p$ has size at most $p^s$, where $p^s$ is the largest power of $p$ dividing $n$:
modulo $p$ we have $t^n-1=(t^m-1)^{p^s}$, with $m=n/p^s$, and $(t^m-1)/(t-1)$ acts invertibly on the
$(t-1)$-primary part. The bound is attained by the pencil of $n\ge3$ lines through the origin of
$\C^2$, for every prime
$p\mid n$: its Milnor fiber is the $n$-fold cyclic cover of a wedge of $n-1$ circles, each sent to a
generator of $\Z/n$, so $H_1(F;\Z)\cong(\LA/(t^n-1))^{n-2}\oplus\LA/(t-1)$.
\end{remark}

\begin{example}
\label{ex:betti}
The failure of $\eta$-formality over $\F_p$ need not be visible in degree one. For an arrangement of rank
$3$ and a prime $p\mid n$, the Aomoto complex $\big(H^\bullet(U;\F_p),\bar\nu\cdot-\big)$ of the proof of
Theorem~\ref{thm:arr} has no cohomology in degree zero, cohomology of dimension $\beta^p_1(M,\eta)$ 
in degree one, and cohomology of dimension $\chi(U)+\beta^p_1(M,\eta)$ in degree two. 
For the braid arrangement $\A_3$ and the
MacLane arrangement at $p=2$, and for the $(9_3)_2$ configuration and the reflection arrangement
$\mathrm{B}_3$ at $p=3$, one has $\beta^p_1(M,\eta)=0$: for $\A_3$ and MacLane at $p=2$ because 
all their flats of rank two and multiplicity greater than two are triple points \cite[Cor.~3.2]{Papadima-Suciu-plms17}; 
for $\mathrm{B}_3$ by \cite[Ex.~7.7]{Papadima-Suciu-plms17}; and for $(9_3)_2$ by
\cite[Thm.~1.6]{Papadima-Suciu-plms17}, as this configuration supports no $3$-net: a $3$-net on its
nine lines, which meet in nine triple points and nine double points, would partition them into three
triples meeting pairwise at double points, whereas the pairs of lines meeting at double points form a
$9$-cycle. So in these 
four cases the Aomoto complex is exact in degree one, and its cohomology lies in degree two,
where it has dimension $\chi(U)$, namely 2, 7, 12, and 8, respectively. For $\A_3$ at $p=3$, 
the other prime dividing $n=6$, one has $\beta^3_1(M,\eta)=1$ \cite[Ex.~7.8]{Papadima-Suciu-plms17}, 
and the failure is already visible in degree one.
\end{example}

\noindent By Theorem~\ref{thm:arr}, Alexander formality of $M(\A)$ depends only on the intersection
lattice, since the decomposition into irreducible factors is that of the matroid into connected
components; in particular, it cannot distinguish the members of a Zariski pair.

\subsection{The four classes}
\label{subsec:arr-corners}

\begin{corollary}
\label{cor:arr-corners}
Each of the four classes of Theorem~\ref{thm:incomp} contains the complement of a central
arrangement, taken with the class of its Milnor fibration:
\begin{enumerate}[label=\textup{(\alph*)},leftmargin=2.2em,itemsep=2pt]
\item \label{co1} 
The Boolean arrangement of $n\ge1$ coordinate hyperplanes is Alexander formal and
strongly formal.
\item \label{co2} 
A pencil of $n\ge3$ lines through the origin of $\C^2$ is strongly formal and not
Alexander formal.
\item \label{co3} 
The Ceva arrangement $\A(3,3,3)$ and the braid arrangement $\A_3$ are neither.
\item \label{co4} 
The product $\A_3\times\A_1$ of the braid arrangement with a single hyperplane is
Alexander formal and not strongly formal.
\end{enumerate}
\end{corollary}

\begin{proof}
Alexander formality is decided by Theorem~\ref{thm:arr}. The Boolean arrangement and $\A_3\times\A_1$ have
an irreducible factor with one hyperplane, so their complements are Alexander formal; a pencil of $n\ge3$
lines, the Ceva arrangement, and $\A_3$ are irreducible, so theirs are not.

For strong formality: $M$ of the Boolean arrangement is $(\C^\ast)^n$, and $M$ of the pencil is
$\C^\ast\times(\mathbb{P}^1\smallsetminus\{n\text{ points}\})$; both are products of spaces homotopy
equivalent to graphs, and so they are strongly formal, by parts~\ref{pr2} and~\ref{pr4} of
Lemma~\ref{lem:products}. The Milnor
fiber of the Ceva arrangement is not $1$-formal \cite{Zuber-blms10}, so its complement is not
strongly formal by Corollary~\ref{cor:milnor-strong}; the complement of $\A_3$ is not strongly formal
by Proposition~\ref{prop:nark}. Finally $M(\A_3\times\A_1)=M(\A_3)\times\C^\ast$ retracts onto $M(\A_3)$,
so it is not strongly formal, by Lemma~\ref{lem:products}\ref{pr1}.
\end{proof}

\noindent In~\ref{co4}, with $w$ the extra coordinate, the diffeomorphism $(x,w)\mapsto(x,Q(x)w)$ of
$M(\A_3)\times\C^\ast$ carries the Milnor fibration $Q\cdot w$ to the projection onto the second
factor; so~\ref{co4} is an instance of Proposition~\ref{prop:dir}. By Theorem~\ref{thm:arr}, an
arrangement complement is Alexander formal only if the arrangement is reducible or consists of a single
hyperplane, but reducibility is not enough: $\A_2\times\A_2$,
the product of two pencils of three lines, has $n=6$ and $3\mid n_1,n_2$, so its complement is not
Alexander formal.

\begin{remark}
\label{rem:line}
The obstructions to strong formality in Corollary~\ref{cor:arr-corners} are of a different nature from
the Massey products that govern Alexander formality, which depend only on the structure of
$\bar A_{\DD_\nu}$, for $\DD_\nu$ the group of characters that factor through $\nu$, as a DG module over
$\bwedge(\eta)$ (Remark~\ref{rem:higher-rank}). Narkawicz's product lies outside every $\bar A_{\DD_\nu}$
(Proposition~\ref{prop:nark}). The obstruction for the Ceva arrangement lies inside $\bar A_{\DD_\nu}$, in
its algebra structure. Indeed, the subgroup $\mu_n\subseteq\DD_\nu$ of characters
of order dividing $n$ corresponds, by Lemma~\ref{lem:graded}, to the $n$-fold cyclic cover $M^{(n)}$ of
$M$ along $\nu$, which is $\{(x,z):z^n=Q(x)\}\cong F\times\C^\ast$. For the Ceva arrangement, the
non-formality of the Milnor fiber makes $\bar A_{\mu_n}$, and hence $\bar A_{\DD_\nu}$ by
Theorem~\ref{thm:G}, non-formal. Yet every product $\langle\eta,\dots,\eta,a\rangle$ vanishes, and the
monodromy of $M^{(n)}$ along $\nu^{(n)}$ is trivial, so the cochains of $M^{(n)}$ are formal as a DG module
over the cohomology of the circle (Proposition~\ref{prop:circle}), although $M^{(n)}$ is not formal.
Neither obstruction bears on the modular clauses, which by Theorem~\ref{thm:arr} depend 
only on the sizes of the irreducible factors.
\end{remark}

\subsection{Localization}
\label{subsec:localization}
The arrangements in Corollary~\ref{cor:arr-corners}\ref{co3} are not isolated: failures of strong
formality propagate along localizations. Flats and localizations were defined in \S\ref{subsec:dir}; 
the flats of $\A$ form the intersection lattice $L(\A)$.

\begin{lemma}
\label{lem:localization}
For every flat $Y$ of $\A$, the inclusion $j\colon M(\A)\inj M(\A_Y)$ has a homotopy section. Hence,
if $M(\A)$ is strongly formal, then so is $M(\A_Y)$.
\end{lemma}

\begin{proof}
The first assertion is \cite[Lem.~4.1]{Denham-Suciu-Yuzvinsky-sel16}. The second assertion follows from
Lemma~\ref{lem:products}\ref{pr1}.
\end{proof}

For a graph $\Gamma$ on the vertices $1,\dots,m$, the \emph{graphic arrangement} $\A_\Gamma$ consists
of the hyperplanes $x_i=x_j$ in $\C^m$, one for each edge $ij$; for the complete graph $K_m$, it is the
braid arrangement $\A_{m-1}$. For $r\ge1$, the \emph{monomial arrangement} $\A(r,r,m)$ is defined by
$\prod_{i<j}(x_i^r-x_j^r)$, and, for $r\ge2$, the monomial arrangement $\A(r,1,m)$ by
$x_1\cdots x_m\prod_{i<j}(x_i^r-x_j^r)$; they are the reflection arrangements of the groups $G(r,p,m)$
with $p=r$ and with $p<r$, respectively.
The arrangements $\A(1,1,m)$, $\A(2,2,m)$, and $\A(2,1,m)$ are the Coxeter arrangements of types
$\mathrm{A}_{m-1}$, $\mathrm{D}_m$, and $\mathrm{B}_m$.

\begin{corollary}
\label{cor:A3-loc}
If $\A$ has a flat $Y$ with $L(\A_Y)\cong L(\A_3)$, then $M(\A)$ is not strongly formal. 
This applies to the following arrangements.
\begin{enumerate}[label=\textup{(\alph*)},leftmargin=2.2em,itemsep=2pt]
\item \label{al1}
The monomial arrangements $\A(r,r,m)$ and $\A(r,1,m)$ for $m\ge4$, among them the braid arrangements
$\A_m$ with $m\ge3$, and the Coxeter arrangements of types $\mathrm{B}_m$ and $\mathrm{D}_m$ with $m\ge4$.
\item \label{al2}
The graphic arrangement of every graph that contains $K_4$.
\end{enumerate}
\end{corollary}

\begin{proof}
By Lemma~\ref{lem:localization}, it suffices to show that $M(\A_Y)$ is not strongly formal. The six
hyperplanes of $\A_Y$ give six lines in the projective plane $\mathrm{P}(\C^\ell/Y)$, with the
intersection lattice of $\A_3$: four triple points, each line through exactly two of them, and any two
of them on one of the lines. No three of the points are collinear, since the six lines joining them are
distinct, and $\mathrm{PGL}_3(\C)$ is transitive on such quadruples; so the six lines are projectively
equivalent to those of $\A_3$, and $M(\A_Y)\simeq M(\A_3)$, which is not strongly formal by
Proposition~\ref{prop:nark}\ref{na1}. In~\ref{al1}, take $Y=\{x_1=x_2=x_3=x_4\}$: it lies on
$x_i=\zeta x_j$ only if $\zeta=1$ and $i,j\le4$, and on no coordinate hyperplane, so $\A_Y$ consists of
the hyperplanes $x_i=x_j$ with $1\le i<j\le4$. In~\ref{al2}, take $Y=\{x_a=x_b=x_c=x_d\}$, for four
pairwise adjacent vertices $a,b,c,d$.
\end{proof}

\noindent More generally, if $Y$ is a flat of the reflection arrangement $\A(W)$ of a finite reflection
group $W$, then, by Steinberg's theorem \cite{Steinberg-tams64}, $\A(W)_Y$ is the reflection
arrangement of the pointwise stabilizer $W_Y$, so Corollary~\ref{cor:A3-loc} applies as soon as $W$ has
a parabolic subgroup
of type $\mathrm{A}_3$, as do the Coxeter groups of types
$\mathrm{E}_6$, $\mathrm{E}_7$, $\mathrm{E}_8$, and $\mathrm{H}_4$. As the braid arrangements are
irreducible, the $\A_m$ with $m\ge3$ form an infinite family in the class of
Corollary~\ref{cor:arr-corners}\ref{co3}. Localization says nothing about
arrangements of rank $3$, such as $\mathrm{B}_3$ or Yoshinaga's arrangement.

\subsection{Graphic and generic arrangements}
\label{subsec:graphic}
We now turn to arrangements with strongly formal complements. For graphic arrangements, a gluing
construction produces many of them, and for chordal graphs the obstruction of
Corollary~\ref{cor:A3-loc}\ref{al2} is the only one; generic arrangements give a second family. To
show that some of these arrangements are irreducible, and hence, by Theorem~\ref{thm:arr}, that their
complements are not Alexander formal, we use a simple observation. As $M(\A)\cong U(\A)\times\C^\ast$,
the Poincar\'e polynomial of $M(\A)$ is $1+t$ times that of $U(\A)$; for a product of two arrangements
it is therefore divisible by $(1+t)^2$. So $\A$ is irreducible whenever $\chi(U(\A))\ne0$.

\begin{proposition}
\label{prop:gluing}
Let $\Gamma_1$ and $\Gamma_2$ be finite graphs whose intersection is empty, a 
single vertex, or a single edge, and let $\Gamma=\Gamma_1\cup\Gamma_2$. 
Then $M(\A_\Gamma)\cong M(\A_{\Gamma_1})\times N$, where $N$
is a retract of $M(\A_{\Gamma_2})$. Consequently, if $M(\A_{\Gamma_1})$ and 
$M(\A_{\Gamma_2})$ are strongly formal, then so is $M(\A_\Gamma)$.
\end{proposition}

\begin{proof}
If $\Gamma_1$ and $\Gamma_2$ are disjoint, then $M(\A_\Gamma)=M(\A_{\Gamma_1})\times M(\A_{\Gamma_2})$,
and we take $N=M(\A_{\Gamma_2})$. Otherwise, let $W$ be the vertex set of 
$\Gamma_1\cap\Gamma_2$, so that $W=\{u\}$, or $W=\{u,w\}$ with $uw$ an
edge, and write the points of $\C^{V(\Gamma)}$ as pairs $(y,z)$, with $y\in\C^{V(\Gamma_1)}$ and
$z\in\C^{V(\Gamma_2)\smallsetminus W}$. Every edge of $\Gamma$ lies in $\Gamma_1$ or in $\Gamma_2$, so
$(y,z)\in M(\A_\Gamma)$ if and only if $y\in M(\A_{\Gamma_1})$ and $(y|_W,z)\in M(\A_{\Gamma_2})$. The
conditions $x_a\ne x_b$ that define $M(\A_{\Gamma_2})$ are preserved by the maps $x\mapsto\alpha
x+\beta$ applied to all coordinates at once, for $\alpha\in\C^\ast$ and $\beta\in\C$. Let $N$ be the
set of $z$ such that $M(\A_{\Gamma_2})$ contains the point with $x_u=0$, respectively with $x_u=0$ and
$x_w=1$, and with the remaining coordinates~$z$. Then the map sending $(y,z)$ to
\[
\big(y,\,z-y_u\big),\quad\text{respectively}\quad\Big(y,\,\frac{z-y_u}{y_w-y_u}\Big),
\]
where $y_u$ is subtracted from every coordinate of $z$, is a diffeomorphism of $M(\A_\Gamma)$ onto
$M(\A_{\Gamma_1})\times N$; it is defined because $y_w\ne y_u$ on $M(\A_{\Gamma_1})$. With
$\Gamma_1\cap\Gamma_2$ in place of $\Gamma_1$, so that $\Gamma$ becomes $\Gamma_2$, 
the same map gives $M(\A_{\Gamma_2})\cong\C\times N$, respectively
$M(\A_{\Gamma_2})\cong\C\times\C^\ast\times N$, so $N$ is a retract of $M(\A_{\Gamma_2})$. The last
assertion follows from Lemma~\ref{lem:products}\ref{pr1} and~\ref{pr2}.
\end{proof}

A graph is \emph{chordal} if every cycle of length at least four has a chord. By a theorem of
Fulkerson and Gross \cite{Fulkerson-Gross-pjm65}, a finite graph is chordal if and only if its
vertices can be ordered so that the earlier neighbors of each vertex are pairwise adjacent.

\begin{theorem}
\label{thm:chordal}
Let $\Gamma$ be a finite chordal graph. Then $M(\A_\Gamma)$ is strongly formal if and only if $\Gamma$
does not contain $K_4$. In that case, $M(\A_\Gamma)$ is diffeomorphic to a product of copies of $\C$,
$\C^\ast$, and $\C\smallsetminus\{0,1\}$.
\end{theorem}

\begin{proof}
If $\Gamma$ contains $K_4$, then $M(\A_\Gamma)$ is not strongly formal, by
Corollary~\ref{cor:A3-loc}\ref{al2}. If not, order the vertices as above. The earlier neighbors of each
vertex are pairwise adjacent, so there are at most two of them, and if there are two, they are joined
by an edge. Adding the vertices one at a time therefore builds $\Gamma$ by taking disjoint unions with
$K_1$, gluing copies of $K_2$ along a vertex, and gluing copies of $K_3$ along an edge. At each step the
factor $N$ of Proposition~\ref{prop:gluing} is $\C$, $\C^\ast$, or $\C\smallsetminus\{0,1\}$,
respectively, so $M(\A_\Gamma)$ is a product of copies of these spaces. Each of them is contractible
or homotopy equivalent to a graph, so $M(\A_\Gamma)$ is strongly formal, by
Lemma~\ref{lem:products}\ref{pr2} and~\ref{pr4}.
\end{proof}

\begin{example}
\label{ex:2tree}
A \emph{$2$-tree} is a graph obtained from an edge by repeatedly adding a vertex joined to both ends of
an existing edge. In the order of construction, each vertex after the first two has exactly two
earlier neighbors, which are adjacent. So a $2$-tree $\Gamma$ on $m\ge3$ vertices is chordal and does
not contain $K_4$, and the proof of Theorem~\ref{thm:chordal} gives
$M(\A_\Gamma)\cong\C\times\C^\ast\times(\C\smallsetminus\{0,1\})^{m-2}$, a strongly formal space. Its
Poincar\'e polynomial is $(1+t)(1+2t)^{m-2}$, so $\chi(U(\A_\Gamma))=(-1)^m\ne0$, and $\A_\Gamma$ is
irreducible, with $2m-3$ hyperplanes; by Theorem~\ref{thm:arr}, $M(\A_\Gamma)$ is not Alexander formal.
Thus $2$-trees give irreducible arrangements of every rank whose complements are strongly formal and
not Alexander formal.
\end{example}

\begin{corollary}
\label{cor:generic}
Let $\A$ consist of $n$ hyperplanes in general position in $\C^\ell$, that is, any $\ell$ of them meet
only at $0$, with $n\ge \ell\ge2$. Then $M(\A)$ is strongly formal. If $n>\ell$, then $\A$ is irreducible, and
$M(\A)$ is not Alexander formal.
\end{corollary}

\begin{proof}
Taking one of the hyperplanes to infinity identifies $U(\A)$ with the complement of $n-1$ affine
hyperplanes in general position in $\C^{\ell-1}$, which is homotopy equivalent to the $(\ell-1)$-skeleton
of the torus $T^{n-1}$, by Hattori's theorem \cite{Hattori-tokyo75}. That skeleton is strongly formal,
by Example~\ref{ex:skeleta}, and hence so is $M(\A)\cong U(\A)\times\C^\ast$, by
Lemma~\ref{lem:products}\ref{pr2} and~\ref{pr4}. The Poincar\'e polynomial of $U(\A)$ is
$\sum_{q<\ell}\tbinom{n-1}{q}t^q$, so $\chi(U(\A))=(-1)^{\ell-1}\tbinom{n-2}{\ell-1}$, which is nonzero
for $n>\ell$. So $\A$ is then irreducible, and $M(\A)$ is not Alexander formal, by Theorem~\ref{thm:arr}.
\end{proof}

\noindent These are the pencils of Corollary~\ref{cor:arr-corners}\ref{co2} for $\ell=2$, and a second
infinite family of this kind for $\ell\ge3$.

\begin{question}
\label{q:graphic}
Is $M(\A_\Gamma)$ strongly formal for every finite graph $\Gamma$ that does not contain $K_4$?
\end{question}

\noindent The answer is affirmative for cycles, by Corollary~\ref{cor:generic}: removing an edge from
a cycle on $m$ vertices leaves a spanning tree, so its graphic arrangement is, up to a factor $\C$, an
arrangement of $m$ hyperplanes in general position in $\C^{m-1}$. With Proposition~\ref{prop:gluing},
it is therefore affirmative for every graph obtained from cycles and complete graphs on at most three
vertices by disjoint unions and by gluing along vertices and edges; this includes the chordal graphs
without $K_4$, the ``house'' (a $4$-cycle and a triangle sharing an edge), and all graphs without $K_4$
on at most four vertices. The smallest graphs for which the question is open have five vertices. 
Among them are $K_{2,3}$, whose only separating pair of vertices is not an edge, and the 
wheel with four spokes, which is $3$-connected.

\subsection*{Acknowledgments}
Parts of this work were carried out with the assistance of Claude
(Anthropic), an AI system; several successive versions of its
Sonnet and Opus models were used for exact computer checks of the explicit 
examples and constructions, for checking citations against the original 
sources, and for suggestions on the exposition.  All mathematical
arguments and computations have been independently verified by the
author, who takes full responsibility for the contents of this paper.


\renewcommand{\MR}[1]{\href{https://mathscinet.ams.org/mathscinet-getitem?mr=#1}{MR~#1}}

\end{document}